\documentclass[10pt]{amsart}
\usepackage{amsmath,amscd,amsthm,amsxtra,yhmath,stackrel}
\usepackage{epsfig,graphics,color,colortbl}
\usepackage{amssymb,latexsym}
\usepackage{mathrsfs}
\usepackage{bm}
\usepackage[poly,all]{xy}
\usepackage{hyperref}
\usepackage{tikz-cd}
\usepackage[draft]{todonotes}
\usepackage{upgreek}
\usepackage{tensor}
\usepackage{ytableau}

\newtheorem{thm}{\bf Theorem}[section]

\newtheorem{prop}[thm]{\bf Proposition}
\newtheorem{cor}[thm]{\bf Corollary}
\newtheorem{lem}[thm]{\bf Lemma}
\newtheorem{rem}[thm]{\bf Remark}
\newtheorem{ex}[thm]{\bf Example}

\numberwithin{equation}{section}

\newcommand{\mc}{\mathcal}
\newcommand{\mf}{\mathfrak}
\newcommand{\mb}{\bm}

\newcommand{\pf}{\noindent{\bfseries Proof. }}
\newcommand{\ov}{\overline}
\newcommand{\un}{\underline}

\newcommand{\U}{{\mc U}}
\newcommand{\W}{\mc{W}}
\newcommand{\cP}{\mathscr{P}}
\newcommand{\I}{\mathbb{I}}
\newcommand{\be}{{\bf e}}

\newcommand{\Z}{\mathbb{Z}}
\newcommand{\Q}{\mathbb{Q}}

\newcommand{\e}{\epsilon}
\newcommand{\de}{\delta}
\newcommand{\ude}{\updelta}
\newcommand{\te}{\widetilde{e}}
\newcommand{\tf}{\widetilde{f}}
\newcommand{\td}{\widetilde}

\newcommand{\La}{\Lambda}
\newcommand{\la}{\lambda}

\newcommand{\hf}{\frac{1}{2}}

\newcommand{\ot}{\otimes}

\newcommand{\si}{(-1)^{\e_i}}
\newcommand{\bq}{{\bf q}}

\usepackage{stmaryrd}
\usepackage{verbatim}
\usepackage{mathtools}
\usepackage[title]{appendix}
\usepackage{stackengine}
\def\dring#1{\stackinset{c}{}{t}{-2.5pt}
  {$
  \mkern2mu\scriptscriptstyle{
  {\scalebox{.4}{\begin{picture}(1,1)
	\put(0,0){\makebox(0,0)[c]{$\circ$}}
	\end{picture}}}
  }
  \mkern-6mu{
  \scalebox{.4}{\begin{picture}(1,1)
	\put(0,0){\makebox(0,0)[c]{$\circ$}}
	\end{picture}}
  }
  $}{$#1$}}
  
\newcommand{\mfrac}[2]{\frac{\textstyle #1}{\textstyle #2}}
\newcommand{\ket}[1]{\left|#1\right\rangle}
\newcommand{\ib}{\bold{i}}

\begin{document}
\title[Representations of quantum affine orthosymplectic superalgebras]
{Oscillator representations of quantum affine orthosymplectic superalgebras II}


\author{JAE-HOON KWON}

\address{Department of Mathematical Sciences and RIM, Seoul National University, Seoul 08826, Korea}
\email{jaehoonkw@snu.ac.kr}

\author{SIN-MYUNG LEE}

\address{School of Mathematics, Korea Institute for Advanced Study, Seoul 02455, Korea}
\email{sinmyunglee@kias.re.kr}

\author{MASATO OKADO}

\address{Osaka Central Advanced Mathematical Institute \&
Department of Mathematics, Osaka Metropolitan University,
3-3-138 Sugimoto, Sumiyoshi-ku, Osaka, 558-8585, Japan}
\email{okado@omu.ac.jp}

\thanks{J.-H. Kwon is supported by the National Research Foundation of Korea(NRF) grant funded by the Korea government(MSIT) (No.2020R1A5A1016126). S.-M. Lee is supported by KIAS Individual Grant (MG095001, MG095002) at Korea Institute for Advanced Study. M. Okado is supported by JSPS KAKENHI Grant Number JP25K06957. This work was partly supported by MEXT Promotion of Distinctive Joint Research Center Program JPMXP0723833165 and Osaka Metropolitan University Strategic Research Promotion Project (Development of International Research Hubs).}

\begin{abstract}
As a continuation of our study on \emph{quantum affine spin-oscillator correspondence}, we introduce a category of $q$-oscillator representations of quantum affine superalgebras of type $XB$ ($X=B,D$) and construct irreducible objects in a uniform way.
Under the truncation functor inspired by the theory of super duality, they are related to finite-dimensional irreducible representations of $U_q(D^{(2)}_{n+1})$ ($X=B$) or $U_q(B^{(1)}_{n})$ ($X=D$).
We also considerably streamline our previous arguments, which clarifies how our construction can be understood as a quantum affinization of classical super duality.
\end{abstract}

\maketitle
\setcounter{tocdepth}{1}
\tableofcontents
\section{Introduction}

\subsection{} 
The theory of \emph{super duality} relates various maximal parabolic BGG categories of Lie (super)algebras by means of interpolating Lie superalgebras \cite{CL,CLW}.
Typically, one relates a module category $\mathcal{O}_{\overline{\boldsymbol{\epsilon}}_n}$ to its bosonic partner $\mathcal{O}_{\underline{\boldsymbol{\epsilon}}_n}$ of a possibly different Lie (super)algebra, or vice versa,
by introducing their super analogue $\mathcal{O}_{\boldsymbol{\epsilon}_n}$ of a (bigger) Lie superalgebra.
It interpolates the first two categories through exact functors that map an irreducible module to an irreducible one or zero, as depicted in the following diagram:
\begin{equation}\label{diag:cls-super-duality}
\xymatrixcolsep{3pc}\xymatrixrowsep{0.3pc}\xymatrix{
 & \mathcal{O}_{\boldsymbol{\epsilon}_n}  \ar@{->}_{\mf{tr}_{\un{\bm\e}}}[dl]\ar@{->}^{\mf{tr}_{\ov{\bm\e}}}[dr] &  \\
 \mathcal{O}_{\underline{\boldsymbol{\epsilon}}_n} & &  \mathcal{O}_{\overline{\boldsymbol{\epsilon}}_n}
}
\end{equation}
The key result is that they induce equivalences of highest weight categories at infinite rank.
This provides a powerful and conceptual way to understand, for example, the theory of finite-dimensional representations of Lie superalgebra $\mathfrak{gl}_{M|N}$
by comparing it directly with a well-studied BGG category of Lie algebra $\mathfrak{gl}_{M+N}$.

The module categories appearing in super duality are not in general closed under taking tensor product.  
On the other hand, the truncation functor itself preserves tensor products by definition, 
and hence one may also connect monoidal structures through  super duality by restricting oneself to subcategories closed under the tensor product.

One such example is the spin-oscillator correspondence \cite{CLW}.
The \emph{oscillator representation}, being a bosonic analogue of the spin representation of $\mathfrak{so}_{2m+1}$ (resp. $\mathfrak{so}_{2m}$), is an infinite-dimensional irreducible representation of the orthosymplectic Lie superalgebra $\mathfrak{osp}_{1|2m}$ (resp. $\mathfrak{sp}_{2m}$).
Surprisingly, its tensor power is semisimple and generates a semisimple monoidal category which has formally the same fusion rule as the category of finite-dimensional represenations of $\mathfrak{so}_{2m+1}$ (resp. $\mathfrak{so}_{2m}$).
This numerical coincidence also follows from the existence of Howe dualities $(\mathfrak{so}_{2m+1},G)$ and $(\mathfrak{osp}_{1|2m},G)$ sharing the same Lie group $G$,
while the theory of super duality enhances it to a larger non-semisimple categories.

\subsection{} 
Let $\mathfrak{g}$ be a Kac--Moody algebra of affine type (without derivation and at level 0) and $U_q(\mathfrak{g})$ the associated quantum group, or a quantum affine algebra.
The Hopf algebra $U_q(\mathfrak{g})$ has intriguing finite-dimensional representation theory
whose major interest lies in understanding the structure of the tensor product of irreducible representations.
The role of (re)normalized $R$-matrix becomes even significant nowadays, especially in the context of monoidal categorification of cluster algebras \cite{KKOP1}.

A finite-dimensional $U_q(\mathfrak{g})$-module $V$ can be seen as an affinization of its restriction to a subalgebra $U_q (\mathring{\mathfrak{g}})$, 
where  $\mathring{\mathfrak{g}}$ is the finite-dimensional simple Lie algebra obtained by removing the 0-th node from the Dynkin diagram of $\mathfrak{g}$.
In other words, $V$ is defined on a direct sum of finite-dimensional irreducible $U_q (\mathring{\mathfrak{g}})$-modules by defining an action of the additional Chevalley generators $e_0 , f_0$. 
The $R$-matrix on the tensor product $V\otimes W$ can then be described by a matrix, indexed by the irreducible components of $V\otimes W$ as a $U_q(\mathring{\mathfrak{g}})$-module, whose entries are rational functions in one variable.
This \emph{spectral decomposition}, especially the poles of the matrix elements, is known to contain crucial information on the $U_q(\mathfrak{g})$-module structure of $V\otimes W$.

This viewpoint can be pushed further to a more general setting as follows.
Given a category $\mathcal{O}$ of $U_q(\mathring{\mathfrak{g}})$-modules, its affinization $\widehat{\mathcal{O}}$ is the category of $U_q(\mathfrak{g})$-modules whose restriction to $U_q(\mathring{\mathfrak{g}})$ belong to $\mathcal{O}$, see \cite{MY} for a systematic study when $\mathcal{O}$ is the BGG category.
If $\mathcal{O}$ is closed under the tensor product, then $\widehat{\mathcal{O}}$ is monoidal and 
even the $R$-matrix can be defined if $\mathcal{O}$ is a subcategory of the BGG category for $U_q(\mathring{\mathfrak{g}})$.
Hence one may try to apply the fruitful theory on finite-dimensional representations to this setting, at least when $\mathcal{O}$ is semisimple.

In this sense, we established a \emph{quantum affinization} of the spin-oscillator correspondence for type $D$ and $C$ in our previous work \cite{KLO24}.
The key evidence was the compelling coincidence of the spectral decomposition of the $R$-matrix on the tensor product of $q$-spin representations of $U_q (D^{(1)}_n)$ and that of $q$-oscillator representations of $U_q (C^{(1)}_n)$ observed in \cite{KwO21}, 
which allowed us to construct irreducible representations of these quantum affine algebras in a uniform way, which are also interpolated under truncation functors.
This may be viewed as the first step towards \emph{quantum affine super duality} whose complete incarnation is still to be developed (see also \cite{KL20}).

In this paper we deal with the remaining case in \cite[Table~1]{KwO21} when $\mathring{\mathfrak{g}}$ is of type $B$. 
More precisely, we define the category $\,^{\mathfrak{x}}\widehat{\mc O}^{\,\mf b}_{\rm osc,\un{\bm\e}}$ of $q$-oscillator representations of type $D^{(2)}_{m+1}$ ($\mathfrak{x}=\mathfrak{b}$) or type $A^{(2)\dagger}_{2m}$ ($\mathfrak{x}=\mathfrak{d}$)
as two affinizations of the one  $\mathcal{O}^{\,\mf b}_{\rm osc,\un{\bm\e}}$ for $U_q (\mathfrak{so}_{2m+1})$.
We compare them with $\,^\mathfrak{x}\widehat{\mc O}^{\,\mf b}_{\rm osc,\ov{\bm\e}}$, the categories of finite-dimensional representations of type $D^{(2)}_{m}$ or $B^{(1)}_{m}$, respectively, via their super analogue $\,^\mathfrak{x}\widehat{\mc O}^{\,\mf b}_{\rm osc,{\bm\e}}$.
Thus, we obtain the following diagram similar to \eqref{diag:cls-super-duality}
\begin{equation}\label{diag:qa-super-duality}
\xymatrixcolsep{3pc}\xymatrixrowsep{0.3pc}\xymatrix{
 & \,^\mathfrak{x}\widehat{\mc O}^{\,\mf b}_{\rm osc,{\bm\e}}  \ar@{->}_{\mf{tr}_{\un{\bm\e}}}[dl]\ar@{->}^{\mf{tr}_{\ov{\bm\e}}}[dr] &  \\
 \,^\mathfrak{x}\widehat{\mc O}^{\,\mf b}_{\rm osc,\un{\bm\e}} & &  \,^\mathfrak{x}\widehat{\mc O}^{\,\mf b}_{\rm osc,\ov{\bm\e}}
}
\end{equation}
where irreducible objects are uniformly constructed so that they match well under the truncation functors $\mf{tr}_{\un{\bm\e}},\mf{tr}_{\ov{\bm\e}}$.

While the argument largely parallels that of the previous work \cite{KLO24}, one notable difference is the use of the twistor automorphism on the covering quantum group \cite{CFLW} to consider a classical limit of $q$-oscillator representations (see \emph{Step 2} below for more details).
This was already implemented in \cite{KwO21} for the non-super case $\mathcal{O}^{\mathfrak{b}}_{\mathrm{osc},\underline{\boldsymbol{\epsilon}}}$, but we are able to extend the theory of twistors naturally to Lie superalgebras with isotropic odd simple roots which seems to be a new result.

Moreover, we streamline our argument in two key ways.
First, we no longer rely on Howe duality to deduce a decomposition of tensor powers of the oscillator representations (Proposition~\ref{prop: decomp of tensor power of W}), which we instead deduce directly from classical super duality.
Second, we prove that the spectral decomposition of $R$-matrices is preserved under the truncation functor once it is verified for generating objects. 
This also makes rather technical computations in \cite[Section~7.2]{KLO24} redundant, although it is still an interesting problem in itself which remains unsolved for fundamental representations of quantum affine algebras in general.

\subsection{}
We shall give an exposition of our construction of \eqref{diag:qa-super-duality} for $\mathfrak{x}=\mathfrak{b}$.
We expect that it will also provide a strategy for establishing the remaining cases of the quantum affine spin-oscillator correspondence proposed in \cite[Section~8.3]{KLO24}.

Let $\mathcal{U}_{\underline{\boldsymbol{\epsilon}}}=U_q(D^{(2)}_{m+2})$ be the quantum affine algebra of type $D^{(2)}_{m+2}$ (resp. $\mathcal{U}_{\overline{\boldsymbol{\epsilon}}}=U_{\widetilde{q}}(D^{(2)}_{m+1})$ of type $D^{(2)}_{m+1}$ with quantum parameter $\widetilde{q}=-q^{-1}$)
and $\mathcal{W}_{\underline{\boldsymbol{\epsilon}}}(x)$ its $q$-oscillator representation (resp. $\mathcal{W}_{\overline{\boldsymbol{\epsilon}}}(x)$ its $q$-spin representation, that is, the $m$-th fundamental representation) with spectral parameter $x$.
Recall that the category of finite-dimensional representations of $\mathcal{U}_{\overline{\boldsymbol{\epsilon}}}$ is generated by $\mathcal{W}_{\overline{\boldsymbol{\epsilon}}}(x)$ as an abelian monoidal category.

\subsubsection*{Step 1.} We introduce a super analogue $\mathcal{W}_{\boldsymbol{\epsilon}}(x)$ of the generating object as a module over a quantum affine superalgebra,
and truncation functors $\mathfrak{tr}_{\epsilon}$ ($\epsilon=\overline{\boldsymbol{\epsilon}},\underline{\boldsymbol{\epsilon}}$) so that $\mathfrak{tr}_{\epsilon}(\mathcal{W}_{\boldsymbol{\epsilon}}(x))\cong \mathcal{W}_{\epsilon}(x)$.

As $\mathcal{W}_{\underline{\boldsymbol{\epsilon}}}(x)$ (resp. $\mathcal{W}_{\overline{\boldsymbol{\epsilon}}}(x)$) is identified with the $q$-symmetric algebra of an $(m+1)$-dimensional space (resp. $q$-exterior algebra of an $m$-dimensional space),
we define $\mathcal{W}_{\boldsymbol{\epsilon}}(x)$ to be the $q$-(super)symmetric algebra of a superspace with a $\mathbb{Z}/2\mathbb{Z}$-graded basis parametrized by $\boldsymbol{\epsilon}=(0,1,0,\dots,1,0)$
together with the action of a \emph{generalized quantum group} $\mathcal{U}_{\boldsymbol{\epsilon}}=\mathcal{U}_{BB}(\boldsymbol{\epsilon})$.
This algebra is a variant of the quantum affine superalgebra introduced in \cite{KOS}, see Proposition~\ref{prop:tau iso} for a precise relationship.
Then we prove the existence of an algebra homomorphism from $\mathcal{U}_{\epsilon}$ to $\mathcal{U}_{\boldsymbol{\epsilon}}$ (Proposition~\ref{prop:reduction}).
While its proof is far less obvious than the embeddings of Lie (super)algebras, the formulas are natural and imply that the truncation functor can be defined just as in the classical super duality.

\subsubsection*{Step 2.} We investigate a direct sum decomposition of tensor powers of $\mathcal{W}_{\epsilon}(x)$ as $\mathring{\mathcal{U}}_\epsilon$-modules ($\epsilon=\boldsymbol{\epsilon},\underline{\boldsymbol{\epsilon}},\overline{\boldsymbol{\epsilon}}$)
and prove that its irreducible components are compatible with the truncation functor.

Equivalently, we consider the semisimple monoidal category $\mathcal{O}^{\mathfrak{b}}_{\mathrm{osc},\epsilon}$ generated by the $\mathring{\mathcal{U}}_\epsilon$-module $\mathcal{W}_{\epsilon}(x)$, where $\mathring{\mathcal{U}}_\epsilon$ is a finite type subalgebra of $\mathcal{U}_\epsilon$.
The truncation functors are defined naturally between these categories,
and we lift results from classical super duality to the quantum case through classical limit to prove the compatibility of irreducible objects under the truncation.
But we should remark that the $q$-oscillator representations of type $B$ are not compatible with the limit $q\rightarrow 1$ (Remark~\ref{rem:classical limit}) unlike the cases in our previous work \cite{KLO24}.
We settle this obstruction by adopting the twistor on a covering quantum group \cite{CFLW} to switch them to modules over $U_q (\mathfrak{osp}_{1|2m})$ which have well-defined classical limit.
The covering quantum group was originally defined only when the underlying Lie superalgebra does not have isotropic odd simple roots.
We show that it can be naturally extended to the generalized quantum group $\mathring{\mathcal{U}}_{\boldsymbol{\epsilon}}$.

\subsubsection*{Step 3.} We prove that the truncation functor between the affinizations $\,^\mathfrak{b}\widehat{\mc O}^{\,\mf b}_{\rm osc,\epsilon}$
preserves the spectral decomposition of normalized $R$-matrices on a tensor product of fundamental type representations, and hence irreducible objects obtained by fusion construction.

Due to a pioneering work \cite{KKKO15}, the use of a normalized $R$-matrix in the study of the tensor product of \emph{real} irreducible objects has expanded considerably in recent years,
and it is even applicable to our categories $\,^\mathfrak{b}\widehat{\mc O}^{\,\mf b}_{\rm osc,\e}$ of possibly infinite-dimensional modules.
Here we perform the fusion construction of irreducible objects as the image of a composition of $R$-matrices.
To this end, we introduce fundamental $q$-oscillator representations for $\epsilon=\boldsymbol{\epsilon},\underline{\boldsymbol{\epsilon}}$, which indeed correspond to finite-dimensional fundamental representations of $\mathcal{U}_{\overline{\boldsymbol{\epsilon}}}$ through the diagram \eqref{diag:qa-super-duality}.
By \emph{Step 2}, the spectral decomposition of the $R$-matrix on the tensor product of them is preserved under the truncation.
In particular, this implies that the property of being real is preserved (Lemma~\ref{lem:unicity of R}).
Therefore, this enables us to give a uniform construction of a large family of irreducible objects in $\,^\mathfrak{b}\widehat{\mc O}^{\,\mf b}_{\rm osc,\e}$ which is compatible with the truncation and contains all the finite-dimensional irreducible representations when $\epsilon=\overline{\boldsymbol{\epsilon}}$.

\subsection*{Organization}
This article is organized as follows.
In Section~\ref{sec:definitions} we recall the notion and basic properties of generalized quantum groups and introduce their $q$-oscillator representations.
We define truncation functors in Section~\ref{sec:truncation}, concluding \textit{Step 1}.
The next two sections are devoted to \textit{Step 2}.
Covering quantum groups and twistor automorphism are introduced in Section~\ref{sec:covering QG} to define twisted classical limit.
In Section~\ref{sec:osc rep} we establish necessary facts on oscillator representations of Lie superalgebras and prove their $q$-analogues.
For the last \textit{Step 3}, in Section~\ref{sec:aff osc category type b}, we define the categories $^{\mathfrak{x}}\widehat{\mathcal{O}}^{\mathfrak{b}}_{\mathrm{osc},\epsilon}$
and construct irreducible objects in a uniform manner in $\epsilon$.

\subsection*{Acknowledgement} 
Part of this work was done while the first two authors were visiting
OCAMI in Osaka Metropolitan University, March 2026. They would like
to thank OCAMI for the invitation and warm hospitality during their visit.
The third-named author would like to thank the Department of Mathematical Sciences of Seoul National University for their hospitality during his stay in September, 2026.

We acknowledge the use of AI for double-checking several proofs in the article.

\section{Quantum affine superalgebras of type $B\!B$ and $D\!B$}\label{sec:definitions}
\subsection{Generalized quantum groups {$\U_{X\!B}(\e)$} of affine type $X\!B$ ($X=D, B$)}\label{subsec:def}
We assume that $n$ is a positive integer such that {$n\geq 4$}.
We denote by $\mathbb{Z}_+$ the set of nonnegative integers. 
Let $\Bbbk=\mathbb{C}(q^\hf)$ and $\Bbbk^{\times}=\Bbbk \setminus \left\{0\right\}$, where $q^\hf$ is an indeterminate. We put 
$[m]=\frac{q^m-q^{-m}}{q-q^{-1}}$ $(m\in \Z_{+})$,
$[m]!=[m][m-1]\cdots [1]\quad (m\geq 1),\quad [0]!=1$, and 
$
\left[\begin{smallmatrix}
m\\[2pt]
k
\end{smallmatrix}\right]
= \frac{[m][m-1]\cdots [m-k+1]}{[k]!}$ $(0\leq k\leq m)$.

Let us introduce the {\em generalized quantum groups of affine type $B\!B$ and $D\!B$}, which we denote by {$\U_{X\!B}(\e)$} ($X=B, D$) (based on \cite{KOS,Ma}).
We largely follow the notations in \cite[Section 2.1]{KLO24} (slightly different depending on the type of $\U_{X\!B}(\e)$):
\begin{itemize}
\item[$\bullet$] $\e=(\e_1,\dots,\e_n)$ : a sequence with $\e_i\in \{0,1\}$ ($i=1,\dots, n$),
 
\item[$\bullet$] $\mathbb{I}= \{1<2<\cdots <n\}$, 

\item[$\bullet$] $P ={\Z\La \oplus } \Z\de_1\oplus\cdots\oplus \Z\de_n$: a $\Z$-lattice with a symmetric bilinear form $(\,\cdot\,|\,\cdot\,)$ satisfying 
\begin{equation*}
(\de_i|\de_j)=(-1)^{\e_i}\delta_{ij},\quad 
{(\de_i|\La)={\hf}},
\end{equation*}

\item[$\bullet$] $I=\{0,1,\ldots,n\}$, $\ring{I}=I\setminus\{0\}$,

\item[$\bullet$] $\alpha_i=\de_{i+1}-\de_i \in P$ $(i\in I\setminus\{0,n\})$, $\alpha_n=-\de_{n}$

\item[$\bullet$] $\alpha_0=\de_{1}$ ($X=B$) or $\de_{1}+\de_{2}$ ($X=D$),

\item[$\bullet$] $q_i=\si q^{\si}$ $(i\in \mathbb{I})$ and $\td{q}=-q^{-1}$

\item[$\bullet$] ${\bf v}={{\bf i}q^{-\hf}}$ (${\bf i}=\sqrt{-1}$),

\item[$\bullet$] ${\bq}(\,\cdot\,,\,\cdot\,)$: a symmetric biadditive function from $P\times P$ to $\Bbbk^{\times}$ given by
{
\begin{equation*}
\begin{split}
&\bq(\mu,\nu) = {{\bf v}}^{{-}\sum_{j\in \I}(\ell'\mu_j+\ell\nu_j)}\prod_{i\in \mathbb{I}}q_i^{\mu_i \nu_i} \quad (\mu,\nu\in P),
\end{split}
\end{equation*}}
for $\mu=\ell\La+\sum_{i\in \I}\mu_i\de_i$ and $\nu=\ell'\La+\sum_{i\in \I}\nu_i\de_i$.
\end{itemize}

Then we first define {$\mathcal{U}_{B\!B}(\epsilon)$} to be the associative $\Bbbk$-algebra generated by 
$e_i,f_i,k_\mu$  $(i\in I, \mu\in P)$ satisfying the following relations:

{\allowdisplaybreaks
\begin{gather}
k_\mu=1\quad (\mu=0), \quad k_{\mu +\mu'}=k_{\mu}k_{\mu'} \quad (\mu, \mu' \in P),\label{eq:Weyl-rel-1-e} \notag \\ 
k_\mu e_i k_{-\mu}=\bq(\mu,\alpha_i)e_i, \notag\quad 
k_\mu f_i k_{-\mu}=\bq(\mu,\alpha_i)^{-1}f_i\quad (i\in I, \mu\in P),  \\ 
e_if_j-f_je_i= \delta_{ij} \displaystyle \frac{k_i-k^{-1} _i}{q-q^{-1}} \quad\text{if  $(i,j)\ne(0,0),(n,n)$}, \label{eq:Weyl-rel-B-1} \\
[e_i ,f_i]= \displaystyle \frac{k_i-k^{-1} _i}{{\bf v}-{\bf v}^{-1}} \quad\text{if }i=0,n, \label{eq:Weyl-rel-B-2} \\
e^2 _i =f^2 _i =0 \quad \text{if $(\alpha_i|\alpha_i)=0$}, \notag\\
e_ie_j-e_je_i=f_if_j-f_jf_i=0 \quad\text{if $(\alpha_i|\alpha_j)=0$ $(i\neq j)$},\notag\\
e_0^3e_1-(1+(-1)^{\epsilon_1}[2])(e_0^2e_1e_0-e_0e_1e_0^2)-e_1e_0^3=(e \rightarrow f)=0, \label{eq:Serre rel for 0 in B} \\
e^2_i e_j - (-1)^{\epsilon_i}[2]e_i e_j e_i +e_j e^2 _i =(e \rightarrow f)=0
\quad \text{if $i\in I\setminus\{0,n\}$, $|i-j|=1$, $\e_i=\e_{i+1}$}, \notag \\
e_n^3e_{n-1}-(1+(-1)^{\epsilon_n}[2])(e_n^2e_{n-1}e_n-e_ne_{n-1}e_n^2)-e_{n-1}e_n^3=(e \rightarrow f)=0, \notag \\
e_{i}e_{i-1}e_{i}e_{i+1}-e_{i}e_{i+1}e_{i}e_{i-1}+(-1)^{\epsilon_i}[2]e_{i}e_{i-1}e_{i+1}e_{i} \label{eq:odd Serre of type A} \\
\hspace{0.1cm}-e_{i-1}e_{i}e_{i+1}e_{i}+e_{i+1}e_{i}e_{i-1}e_{i}=(e \rightarrow f)=0
\quad\text{if $i\in \{1,\dots,n-1\}$, $\e_i\neq \e_{i+1}$}. \notag
\end{gather}
Here $k_i=k_{\alpha_i}$ for $i\in I$, and $(e\rightarrow f)$ denotes the relation which is obtained by replacing $e_i$ with $f_i$.

Next, we define {$\mathcal{U}_{D\!B}(\epsilon)$} to be the associative $\Bbbk$-algebra generated by 
$e_i,f_i,k_\mu$  $(i\in I, \mu\in P)$ satisfying the same relations as in $\mathcal{U}_{B\!B}(\epsilon)$ except that \eqref{eq:Weyl-rel-B-1} is satisfied for $(i,j)\neq (n,n)$ and \eqref{eq:Weyl-rel-B-2} is satisfied for $i=n$, 
and \eqref{eq:Serre rel for 0 in B} is replaced with
\begin{gather*}
e^{2}_0 e_2 - (-1)^{\epsilon_1}[2]e_0 e_2 e_0 +e_2 e^{2}_0=(e\rightarrow f)=0 \quad \text{if }\epsilon_1= \epsilon_2,  \\
e_0e_1e_2-e_1e_0e_2+(-1)^{\epsilon_2}[2](e_1e_2e_0-e_0e_2e_1)
+e_2e_0e_1-e_2e_1e_0=(e\rightarrow f)=0 \quad \text{if $\e_1\neq \e_2$}, \\
e_0e_2e_3e_2-e_3e_2e_0e_2+(-1)^{\epsilon_3}[2]e_2e_3e_0e_2
-e_2e_0e_2e_3+e_2e_3e_2e_0=(e \rightarrow f)=0
\quad \text{if $\e_2\neq \e_3$}. 
\end{gather*}

It has a Hopf algebra structure with the comultiplication $\Delta$ and the antipode $S$
\begin{equation}\label{eq:comult-e}
\begin{split}
\Delta(k_\mu)&=k_\mu\otimes k_\mu, \\ 
\Delta(e_i)&= 1\ot e_i + e_i\ot k_i^{-1}, \\
\Delta(f_i)&= f_i\ot 1 + k_i\ot f_i , \\  
S(k_\mu)=k_{-\mu},& \ \ S(e_i)=-e_ik_i, \ \  S(f_i)=-k_i^{-1}f_i,\\
\end{split}
\end{equation}
for $\mu\in P$ and  $i\in I$.

For $X=B$, we have $\U_{B\!B}(0^n)\cong U_q(D_{n+1}^{(2)})$ and $\U_{B\!B}(1^n)\cong U_{\td{q}}(D_{n+1}^{(2)})$ once we rescale $e_0$ and $e_n$ by multiplying
the factor
\begin{equation}\label{eq:scalar c}
c:=\frac{{\bf v}-{\bf v}^{-1}}{q^{\frac12}-q^{-\frac12}}.
\end{equation}
For $X=D$, we have $\U_{D\!B}(0^n)\cong U_q(B_{n}^{(1)})$ and $\U_{D\!B}(1^n)\cong U_{\td{q}}(B_{n}^{(1)})$ (with the same rescaling of $e_n$).
For example, the associated Dynkin diagrams are as follows:\vskip 5mm

\begin{center}\setlength{\unitlength}{0.16in} \smallskip
\raisebox{4ex}{$\U_{B\!B}(00011)$ :}
\begin{picture}(18,3)
\put(2,2){\makebox(0,0)[c]{$\bigcirc$}}
\put(5,2){\makebox(0,0)[c]{$\bigcirc$}}
\put(8,2){\makebox(0,0)[c]{$\bigcirc$}}
\put(11,1.95){\makebox(0,0)[c]{$\bigotimes$}}
\put(14,2){\makebox(0,0)[c]{$\bigcirc$}}
\put(17,2){\makebox(0,0)[c]{$\bigcirc$}}

\put(5.35,2){\line(1,0){2.1}} 
\put(8.4,2){\line(1,0){2.1}} 
\put(11.4,2){\line(1,0){2.1}}
\put(14.4,1.9){\line(1,0){2}}
\put(14.4,2.1){\line(1,0){2}}
\put(16.3,2){\makebox(0,0)[c]{\Large $>$}}

\put(2.6,1.9){\line(1,0){2}}
\put(2.6,2.1){\line(1,0){2}}
\put(2.7,2){\makebox(0,0)[c]{\Large $<$}}

\put(2,1){\makebox(0,0)[c]{\tiny $0$}}
\put(5,1){\makebox(0,0)[c]{\tiny $1$}}
\put(8.2,1){\makebox(0,0)[c]{\tiny $2$}}
\put(11,1){\makebox(0,0)[c]{\tiny $3$}}
\put(14,1){\makebox(0,0)[c]{\tiny $4$}}
\put(17,1){\makebox(0,0)[c]{\tiny $5$}}
\end{picture}
\end{center}

\begin{center}\setlength{\unitlength}{0.16in} \smallskip
\raisebox{4ex}{$\U_{B\!B}(01010)$ :}
\begin{picture}(18,3)
\put(2,2){\makebox(0,0)[c]{$\bigcirc$}}
\put(5,2){\makebox(0,0)[c]{$\bigotimes$}}
\put(8,2){\makebox(0,0)[c]{$\bigotimes$}}
\put(11,1.95){\makebox(0,0)[c]{$\bigotimes$}}
\put(14,2){\makebox(0,0)[c]{$\bigotimes$}}
\put(17,2){\makebox(0,0)[c]{$\bigcirc$}}

\put(5.35,2){\line(1,0){2.1}} 
\put(8.4,2){\line(1,0){2.1}} 
\put(11.4,2){\line(1,0){2.1}}
\put(14.4,1.9){\line(1,0){2}}
\put(14.4,2.1){\line(1,0){2}}
\put(16.3,2){\makebox(0,0)[c]{\Large $>$}}

\put(2.6,1.9){\line(1,0){2}}
\put(2.6,2.1){\line(1,0){2}}
\put(2.7,2){\makebox(0,0)[c]{\Large $<$}}

\put(2,1){\makebox(0,0)[c]{\tiny $0$}}
\put(5,1){\makebox(0,0)[c]{\tiny $1$}}
\put(8.2,1){\makebox(0,0)[c]{\tiny $2$}}
\put(11,1){\makebox(0,0)[c]{\tiny $3$}}
\put(14,1){\makebox(0,0)[c]{\tiny $4$}}
\put(17,1){\makebox(0,0)[c]{\tiny $5$}}
\end{picture}
\end{center}

\begin{center}\setlength{\unitlength}{0.16in} \smallskip
\raisebox{4ex}{$\U_{D\!B}(00011)$ :}
\begin{picture}(18,5)
\put(6,0){\makebox(0,0)[c]{$\bigcirc$}}
\put(6,4){\makebox(0,0)[c]{$\bigcirc$}}
\put(8,2){\makebox(0,0)[c]{$\bigcirc$}}
\put(11,1.95){\makebox(0,0)[c]{$\bigotimes$}}
\put(14,2){\makebox(0,0)[c]{$\bigcirc$}}
\put(17,2){\makebox(0,0)[c]{$\bigcirc$}}

\put(14.4,1.9){\line(1,0){2}}
\put(14.4,2.1){\line(1,0){2}}
\put(16.3,2){\makebox(0,0)[c]{\Large $>$}}

\put(6.35,0.3){\line(1,1){1.35}} \put(6.35,3.7){\line(1,-1){1.35}}
\put(8.4,2){\line(1,0){2.1}} \put(11.4,2){\line(1,0){2.1}}


\put(6,5){\makebox(0,0)[c]{\tiny $0$}}
\put(6,-1){\makebox(0,0)[c]{\tiny $1$}}
\put(8.2,1){\makebox(0,0)[c]{\tiny $2$}}
\put(11,1){\makebox(0,0)[c]{\tiny $3$}}
\put(14,1){\makebox(0,0)[c]{\tiny $4$}}
\put(17,1){\makebox(0,0)[c]{\tiny $5$}}
\end{picture}
\end{center}\vskip 5mm

\begin{center}\setlength{\unitlength}{0.16in} \smallskip
\raisebox{4ex}{$\U_{D\!B}(10101)$ :}
\begin{picture}(18,5)
\put(6,0){\makebox(0,0)[c]{$\bigotimes$}}
\put(6,4){\makebox(0,0)[c]{$\bigotimes$}}
\put(8,2){\makebox(0,0)[c]{$\bigotimes$}}
\put(11,1.95){\makebox(0,0)[c]{$\bigotimes$}}
\put(14,2){\makebox(0,0)[c]{$\bigotimes$}}
\put(17,2){\makebox(0,0)[c]{$\bigcirc$}}

\put(14.4,1.9){\line(1,0){2}}
\put(14.4,2.1){\line(1,0){2}}
\put(16.3,2){\makebox(0,0)[c]{\Large $>$}}

\put(6.35,0.3){\line(1,1){1.35}} \put(6.35,3.7){\line(1,-1){1.35}}
\put(8.4,2){\line(1,0){2.1}} \put(11.4,2){\line(1,0){2.1}}
\put(6,0.5){\line(0,1){3.1}} 


\put(6,5){\makebox(0,0)[c]{\tiny $0$}}
\put(6,-1){\makebox(0,0)[c]{\tiny $1$}}
\put(8.2,1){\makebox(0,0)[c]{\tiny $2$}}
\put(11,1){\makebox(0,0)[c]{\tiny $3$}}
\put(14,1){\makebox(0,0)[c]{\tiny $4$}}
\put(17,1){\makebox(0,0)[c]{\tiny $5$}}
\end{picture}
\end{center}

\medskip

We also use the notation
\begin{itemize}
\item[$\bullet$] $\ring{\U}_{X\!B}(\e)$: the subalgebra of ${\U}_{X\!B}(\e)$ generated by $e_i$, $f_i$, $k_\mu$ ($i\in \ring{I}$, $\mu\in P$),

\item[$\bullet$] $\dring{\U}_{X\!B}(\e)$: the subalgebra of ${\U}_{X\!B}(\e)$ generated by $e_i$, $f_i$, $k_\mu$ ($i\in \ring{I}\setminus \{n\}$, $\mu\in P$).
\end{itemize}
Note that they are independent of $X$ and can be identified with the quantum groups associated to Lie superalgebras of type $B$ and $A$, respectively.

For a {$\U_{X\!B}(\e)$}-module $V$ and $\la\in P$, we define the $\la$-weight space to be 
\begin{equation}\label{eq:weight space}
 V_\la = \{\,v\in V\,|\,k_{\mu} v= {\bq}(\la,\mu) v \ \ (\mu\in P) \,\},
\end{equation}
and write ${\rm wt}(v)=\la$ for $v\in V_\la$.
Let ${\rm wt}(V)=\{\,\mu\in P\,|\,V_\mu\neq 0\,\}$.

\subsection{Quantum affine superalgebras {$U_{X\!B}(\e)$} of type $X\!B$}\label{subsec:Yamane QSA}

Let $U_{B\!B}(\epsilon)$ be the associative $\Bbbk$-algebra with 1 generated by $E_{i}$, $F_{i}$, $K_{\mu}$ $(\mu\in P, i\in I)$ satisfying the following relations:
\allowdisplaybreaks{
\begin{gather}
  K_{\mu}=1 \quad(\mu=0),\quad K_{\mu+\mu'}=K_{\mu}K_{\mu'} \quad (\mu, \mu' \in P),\notag \\
  K_{\mu}E_{i}K_{\mu}^{-1}=q^{(\mu|\alpha_{i})}E_{i},\quad K_{\mu}F_{i}K_{\mu}^{-1}=q^{-(\mu|\alpha_{i})}F_{i} \quad (i\in I, \mu\in P),\notag \\
  E_{i}F_{j}-(-1)^{p(i)p(j)}F_{j}E_{i}=(-1)^{\epsilon_{i}}\delta_{ij}\frac{K_{i}-K_{i}^{-1}}{q-q^{-1}} \quad \text{if  $(i,j)\ne(0,0),(n,n)$},\label{eq:Weyl-rel-B-1 Yamane}\\
{E_{i}F_{i}-F_{i}E_{i}= \frac{K_{i}-K_{i}^{-1}}{{\bf v}-{\bf v}^{-1}} \quad \text{if $i=0,n$}},\label{eq:Weyl-rel-B-2 Yamane}\\  
  E_{i}^{2}=F_{i}^{2}=0\quad \text{if $(\alpha_i|\alpha_i)=0$},\notag \\
  E_{i}E_{j}-(-1)^{p(i)p(j)}E_{j}E_{i}=(E\rightarrow F)=0\quad\text{if $(\alpha_i|\alpha_j)=0$},\notag \\
{E_0^3E_1-(1+{[2]})(E_0^2E_1E_0-E_0E_1E_0^2)-E_1E_0^3=(E \rightarrow F)=0}, \label{eq:Serre rel for 0 in B Yamane} \\
  E_{i}^{2}E_{j}-[2]E_{i}E_{j}E_{i}+E_{j}E_{i}^{2}=(E\rightarrow F)=0\quad \text{if $i,j\in I\setminus\{0,n\}$, $|i-j|=1$, $\e_i=\e_{i+1}$},\notag \\
{E_n^3E_{n-1}-(1+[2])(E_n^2E_{n-1}E_n-E_nE_{n-1}E_n^2)-E_{n-1}E_n^3=(E \rightarrow F)=0}, \notag \\
\left[E_{i},\left[\left[E_{i-1},E_{i}\right]_{(-1)^{p(i-1)}q},E_{i+1}\right]_{(-1)^{(p(i-1)+p(i))p(i+1)}q^{-1}}\right]_{(-1)^{p(i-1)+p(i)+p(i+1)}} \label{eq:odd Serre of type A Yamane}\\
  \hfill=(E\rightarrow F)=0\quad\text{if $i\in \{1,\dots,n-1\}$, $\e_i\neq \e_{i+1}$} \notag,
\end{gather}}
\!\!where $[X,Y]_{t}=XY-tYX$ and $p(i)=0$ (resp.\,$1$) for $i\in I$ with $(\alpha_i|\alpha_i)\neq 0$ (resp. $i\in I$ with $(\alpha_i|\alpha_i)=0$). It is a $\Z_2$-graded algebra with ${\rm deg}(K_\mu)=0$ for $\mu\in P$ and ${\rm deg}(E_i)={\rm deg}(F_i)=\ov{p(i)}$ for $i\in I$.

Also, let ${U}_{D\!B}(\epsilon)$ be the associative $\Bbbk$-algebra generated by 
$E_i,F_i,K_\mu$  $(i\in I, \mu\in P)$ satisfying the same relations as in ${U}_{B\!B}(\epsilon)$ except that \eqref{eq:Weyl-rel-B-1 Yamane} is satisfied for all $(i,j)\neq (n,n)$ (where we put $\epsilon_0 = \epsilon_1$) and \eqref{eq:Weyl-rel-B-2 Yamane} is satisfied for $i=n$,
and \eqref{eq:Serre rel for 0 in B Yamane} is replaced with
\begin{gather*}
  E^{2}_0 E_2 -[2] E_0 E_2 E_0 + E_2 E^{2}_0 = (E\rightarrow F)=0\quad\text{if }\epsilon_1=\epsilon_2,\\
  E_{0}E_{1}E_{2}-E_{1}E_{0}E_{2}+(-1)^{\epsilon_{2}+\epsilon_{3}}[2](E_{1}E_{2}E_{0}-E_{0}E_{2}E_{1})+E_{2}E_{0}E_{1}-E_{2}E_{1}E_{0}\\ =(E\rightarrow F)=0\quad\text{ if \ensuremath{\epsilon_{1}\neq\epsilon_{2}}},\\
  \left[E_{2},\left[\left[E_{0},E_{2}\right]_{(-1)^{p(0)}q},E_{3}\right]_{(-1)^{(p(0)+1)p(3)}q^{-1}}\right]_{(-1)^{p(0)+p(3)+1}}=(E\rightarrow F)=0\quad\text{ if }\epsilon_{2}\neq\epsilon_{3},
\end{gather*}
where we set $p(0)=p(1)$, $\deg(E_0)=\deg(F_0)=\overline{p(0)}$ and the others are as above.
The algebra {$U_{X\!B}(\e)$ ($X=D, B$)} is the quantum affine superalgebra of type (XB) in \cite[Section 1.4 and Proposition 6.7.1]{Ya99} after rescaling $E_0$ and $E_n$ by multiplying $c$ in \eqref{eq:scalar c} when the Weyl relation is given by \eqref{eq:Weyl-rel-B-2 Yamane}.

As in the case of $\U_D(\e)$ \cite[Section 2.2]{KLO24}, $\U_{X\!B}(\e)$ and $U_{X\!B}(\e)$ are related as follows. Let
\begin{itemize}
\item[$\bullet$] $P_{\rm af} = \bigoplus_{j\in \I}\Z\mb{\delta}_j \oplus \Z\mb{\delta} \oplus \Z \boldsymbol{\Lambda}$ with a symmetric bilinear form
\begin{equation*}
(\mb{\delta}_i|\mb{\delta}_j)=(-1)^{\e_i}\delta_{ij},\quad (\boldsymbol{\Lambda}|\mb{\delta}_i)=\hf,\quad (\mb{\delta}_i|\mb{\delta})=(\mb{\delta}|\mb{\delta})=(\boldsymbol{\Lambda}|\mb{\delta})=0,
\end{equation*}

\item[$\bullet$] $\mb{\alpha}_i=\mb{\de}_{i+1}-\mb{\de}_i$ $(i\in I\setminus\{0,n\})$, $\mb{\alpha}_n=-\mb{\de}_{n}$,

\item[$\bullet$] $\mb{\alpha}_0=\mb{\de}_{1}+\mb{\de}$ ($X=B$) or $\mb{\de}_{1}+\mb{\de}_2+\mb{\de}$ ($X=D$)
\end{itemize}
Then we may regard $P\subset P_{\rm af}$ by identifying $\de_i$ with $\mb{\de}_i$ ($i\in \I$), and $\La$ with $\boldsymbol{\Lambda}$.

Let {${\rm cl} : P_{\rm af} \longrightarrow P$}
be the projection along $\boldsymbol{\delta}$ onto $P$ under this embedding.
We take a $\Bbbk$-algebra $\U_{X\!B}(\e)_{\rm af}$, defined in the same way as $\mathcal{U}_B (\epsilon)$ but $P$ and $\bq(\,\cdot\,,\,\cdot\,)$ are now replaced by $P_{\rm af}$ and $\bq(\mu,\nu)=\bq({\rm cl}(\mu),{\rm cl}(\nu))$ for $\mu,\nu\in P_{\rm af}$.

{
Let $\Sigma$ be the commutative $\Bbbk$-bialgebra generated by $\sigma_i$ for each $i\in \I$ such that $\epsilon_i=1$ and ${\bm \sigma}$ subject to $\sigma_j^2={\bm\sigma}^4=1$ for $i,j\in \I$, with the comultiplication $\Delta(\sigma_i)=\sigma_i\otimes \sigma_i$ and $\Delta({\bm\sigma})={\bm\sigma}\otimes {\bm\sigma}$ for $i\in \I$.
As a placeholder we also set $\sigma_i=1$ for $i\in\mathbb{I}$ with $\epsilon_i=0$.
}
We define {${\U_{X\!B}(\e)}_{\rm af}[\sigma]$}
to be the semidirect product of ${\U_{X\!B}(\e)}_{\rm af}$ and $\Sigma$, where $\Sigma$ acts on ${\U_{X\!B}(\e)}_{\rm af}$ by
\begin{gather*}
\sigma_j \cdot k_\mu =k_\mu,\quad \sigma_j \cdot e_i=(-1)^{\e_j(\mb{\de}_j|\mb{\alpha}_i)}e_i,\quad \sigma_j \cdot f_i=(-1)^{\e_j(\mb{\de}_j|\mb{\alpha}_i)}f_i, \\
{\bm \sigma} \cdot k_\mu =k_\mu,\quad {\bm \sigma} \cdot e_i=\ib^{-2(\boldsymbol{\Lambda}|\mb{\alpha}_i)} e_i,\quad {\bm \sigma} \cdot f_i=\ib^{2(\boldsymbol{\Lambda}|\mb{\alpha}_i)}f_i,
\end{gather*}
for $j\in \I$, $\mu\in P_{\rm af}$ and $i\in I$.  
We define ${U_{X\!B}(\e)}_{\rm af}[\sigma]$ in the same manner.

\begin{prop}\label{prop:tau iso}
Suppose that $\e=(\e_1,\dots,\e_n)=(0,1,0,\dots,0,1,0)$ for $X=B$ or $\e=(\e_1,\dots,\e_n)=(1,0,1,0,\dots,0,1,0)$ for $X=D$.
There exists an isomorphism of $\Bbbk$-algebras
$\tau : \mathcal{U}_{XB}(\epsilon)_{\mathrm{af}}[\sigma] \longrightarrow U_{XB}(\epsilon)_{\mathrm{af}}[\sigma]$
such that 
{\allowdisplaybreaks
\begin{equation*}
\begin{split}
& \tau(\sigma_a) =\sigma_a,\quad {\tau(k_{\boldsymbol{\delta}_i})=K_{\boldsymbol{\delta}_i}\sigma_i,\quad \tau(k_{\La})=K_{\La}{\bm\sigma}},\quad \tau(k_\delta)=K_\delta,\\
& \tau(e_i)= 
\begin{cases}
E_i & \text{for $i=n$ or $i=0$ with $X=B$},\\
E_0 \sigma_{\le 1}\varsigma_1 & \text{for $i=0$ with $X=D$},\\
E_i \sigma_{\le i}\varsigma_i^{i} & \text{for $1\leq i\leq n-1$},
\end{cases}\\
& \tau(f_i)= 
\begin{cases}
F_i & \text{for $i=n$ or $i=0$ with $X=B$},\\
-F_0 \sigma_{\le 1} & \text{for $i=0$ with $X=D$},\\
(-1)^{i}F_i \sigma_{\le i}\varsigma_i^{i+1} & \text{for $1\leq i\leq n-1$},
\end{cases}
\end{split}
\end{equation*}}
where $\sigma_{\leq i}=\sigma_{1}\sigma_{2}\cdots\sigma_{i}$ and $\varsigma_i=\sigma_i\sigma_{i+1}$ for $1\le i\le n-1$.
\end{prop}
}

\pf Let us consider the case of $X=B$. The proof of $X=D$ is similar. Since $\dring{\U}_{B\!B}(\e)$ is of type $A$, $\tau$ is well-defined on $\dring{\U}_{B\!B}(\e)$ \cite[Section 2.3]{KL20}. So it suffices to consider the relations where $e_i, f_i, k_i^{\pm 1}$ ($i=0,n$) are involved. Also we consider only the case of $i=0$ since the case of $i=n$ is symmetric.

First, consider the Serre relation where $e_0$ and $e_j$ with $(\alpha_i|\alpha_j)\neq 0$ are involved. 
Let 
\begin{equation*}
\begin{split}
\mc{S}_{n,n-1}(\e)
&= {e_n^3e_{n-1}-(1+(-1)^{\epsilon_n}[2])(e_n^2e_{n-1}e_n-e_ne_{n-1}e_n^2)-e_{n-1}e_n^3}, \\
&= e_n^3e_{n-1}-(1+[2])(e_n^2e_{n-1}e_n-e_ne_{n-1}e_n^2)-e_{n-1}e_n^3 \quad (\e_n=0),\\
S_{n,n-1}(\e)
&= E_n^3E_{n-1}-(1+[2])(E_n^2E_{n-1}e_n{-}E_nE_{n-1}E_n^2){-}E_{n-1}E_n^3
\end{split}
\end{equation*}
be the corresponding relations in $\U_{B\!B}(\e)$ and $U_{B\!B}(\e)$ respectively.
Then it is clear that $\tau(\mc{S}_{n,n-1}(\e))={S}_{n,n-1}(\e)\sigma'$ for some $\sigma'\in \Sigma$, since $\tau(e_n)=E_n$ and $\sigma_aE_n=E_n\sigma_a$ for all $a$. The relation for $e_0$ and $e_j$ with $(\alpha_i|\alpha_j)=0$ is also clear.
\smallskip

Next, consider the Weyl relation. It is straightforward to see
\[
\tau(e_nf_k-f_ke_n)=
\begin{cases}
\frac{\tau(k_n)-\tau(k^{-1}_n)}{\mathbf{v}-\mathbf{v}^{-1}}=E_n F_n -F_n E_n & \text{if $k=n$}, \\
0=(E_nF_k - F_k E_n)\sigma' & \text{if $k\neq n$}, 
\end{cases}
\]
where $\sigma'\in \Sigma$ (depending on $k$). 
\qed

\begin{rem}\label{rem: iso GQG to YQG}
{\rm 
There also exists an isomorphism for any $\e$ with $\e_1=\e_n=0$ when $X=B$, and for $\e$ with $\e_n=0$ when $X=D$, which can be proved in a similar manner (cf.~\cite[Proposition 2.1]{KLO24}). Hence, if we restrict $\tau$ to $\ring{\U}_{X\!B}(\e)[\sigma]$ and $\ring{U}_{X\!B}(\e)[\sigma]$, respectively, then there exists an isomorphism for $\e$ with $\e_n=0$.
} 
\end{rem}

\smallskip

\subsection{Quasi $R$-matrix for $\U_{X\!B}(\e)$}\label{subsec:quasi-R}
Let $\mathcal{U}_{X\! B}(\epsilon)^+$ (resp. $\mathcal{U}_{X\! B}(\epsilon)^-$) be the subalgebra of $\mathcal{U}_{X\! B}(\epsilon)$
generated by $e_i$ (resp. $f_i$) for $i\in I$. Let {$Q_{\pm}:=\pm\sum_{i\in I}\mathbb{Z}_{+}\boldsymbol{\alpha}_i$}. There is a natural $Q_\pm$-grading on $\mathcal{U}_{X\! B}(\epsilon)^{\pm}$.
The following result is an analogue of the construction of Lusztig \cite[Chapter~1]{Lu93} which can be proved as in the super type $A$ case \cite[Sections~2.2--2.4]{KL20}.

\begin{prop}\mbox{}

\begin{itemize}
\item[(1)] For $i \in I$, there exists a unique $\Bbbk$-linear map $_i r:\mathcal{U}_{X\! B}(\epsilon)^-\rightarrow \mathcal{U}_{X\! B}(\epsilon)^-$ such that $_i r(f_j)=\delta_{ij}$ for $j\in I$, and
$_i r(xy)= \,_i r(x)y +\mathbf{q}(\alpha_i, \beta)x \,_i r(y)$  for $x\in \mathcal{U}_{X\! B}(\epsilon)^{-}_{\beta}$ and $y\in \mathcal{U}_{X\! B}(\epsilon)^-$.

\item[(2)] There exists a unique non-degenerate symmetric bilinear form $(\,\cdot\,,\,\cdot\,)$ on $\mathcal{U}_{X\! B}(\epsilon)^-$
such that $(1,1)=1$ and $(f_i x,y)=(x,\,_i r(y))$ for $i\in I$.

\end{itemize}
\end{prop}

For each $\beta\in Q_{+}$, take a basis $B_\beta$ of $\mathcal{U}_{X\! B}(\epsilon)^{-}_{-\beta}$ and its dual $B^{*}_\beta=\{\,b^* \,|\, b\in B_\beta\,\}$ with respect to the bilinear form above.
We also denote by $b^+ \in \mathcal{U}_{X\! B}(\epsilon)^{+}_{\beta}$ the one obtained from $b\in B_\beta$ by replacing each $f_i$ with $e_i$ for all $i\in I$.
Let $\overline{\,\cdot\,}$ be the $\mathbb{C}$-linear involution on $\mathcal{U}_{X\! B}(\epsilon)$ defined by $\overline{e_i}=e_i$, $\overline{f_i}=f_i$, $\overline{k}_\mu=k_{-\mu}$ ($i\in I$, $\mu\in P$) and $\overline{q^\hf}=q^{-\hf}$.

\begin{thm}\mbox{}
\begin{itemize}
\item[(1)] For each $\beta \in Q_+$, there exists a unique $\Theta_\beta \in \mathcal{U}_{X\! B}(\epsilon)^{+}_{\beta}\otimes \mathcal{U}_{X\! B}(\epsilon)^{-}_{-\beta}$
such that $\Theta_0 = 1\otimes 1$ and $\Theta \Delta(u)=\overline{\Delta(\overline{u})}\Theta$  for all $u\in \mathcal{U}_{X\! B}(\epsilon)$ and $\Theta = \sum_{\beta \in Q^+}\Theta_\beta$.
Moreover, $\Theta_\beta$ is given by
\[
\Theta_\beta = (q-q^{-1})^{\mathrm{ht}(\beta)}\sum_{b\in B_\beta} b^+ \otimes b^*,
\]
where $\mathrm{ht}(\beta)=\sum a_i$ if $\beta=\sum_{i\in I} a_i \boldsymbol{\alpha}_i$.

\item[(2)] If we let $\overline{\Theta}_\beta=(\overline{\,\cdot\,}\otimes \overline{\,\cdot\,})(\Theta_\beta)$ and $\overline{\Theta}=\sum\overline{\Theta}_\beta$, then we have $\Theta \overline{\Theta}=\overline{\Theta}\Theta = 1$.

\end{itemize}
\end{thm}

The summation $\Theta$ is given in a suitable completion of $\mathcal{U}_{X\!B}(\epsilon)^+ \otimes \mathcal{U}_{X\!B}(\epsilon)^-$ (see \cite[Chapter~4]{Lu93} for detailed arguments).
The element $\Theta$ is called the \emph{quasi $R$-matrix} for $\mathcal{U}_{X\!B}(\epsilon)$.

\subsection{A $\U_{X\!B}(\e)$-module $\W(x)$}\label{subsec:W}
Let 
\begin{itemize}
 \item[$\bullet$] $\Z_+^n(\epsilon)=\{\,{\bf m}=(m_1,\ldots,m_n)\in \mathbb{Z}^{n}_{+}\,|\, \epsilon_i=1\Rightarrow m_i\in\{0,1\}\,\}$, ${\bf 0}=(0,\dots,0)$,

 \item[$\bullet$] $\{\,\be_i\,|\,1\le i\le n\,\}$: the standard basis of $\Z^n$,
 
 \item[$\bullet$] $\W = \bigoplus_{{\bf m}\in \Z^n_+(\e)}\Bbbk|{\bf m}\rangle$: a $\Bbbk$-space spanned by $\ket{{\bf m}}$.
 
\end{itemize}

\begin{prop}\label{prop:osc W for e}
Let $x\in\Bbbk^\times$ be given.
\begin{itemize}
\item[(1)] If $X=B$, then the following action of the generators defines a $\U_{B\!B}(\e)$-module structure on $\W$:
{\allowdisplaybreaks 
\begin{align*}
&e_0 |{\bf m}\rangle =x|{\bf m}+{\bf e}_1 \rangle, \\
&f_0 |{\bf m}\rangle =x^{-1}[m_1]|{\bf m}-{\bf e}_1 \rangle,\\
&k_0 |{\bf m}\rangle ={\bf v}^{-1}q_1^{m_1}|{\bf m}\rangle, \\
&e_i |{\bf m}\rangle = [m_{i}]|{\bf m}-{\bf e}_{i}+{\bf e}_{i+1} \rangle,\notag \\
&f_i |{\bf m}\rangle = [m_{i+1}]|{\bf m}+{\bf e}_{i}-{\bf e}_{i+1}\rangle, \notag \\
&k_i |{\bf m}\rangle =q^{-m_i} _i q^{m_{i+1}} _{i+1}  |{\bf m}\rangle,\notag \\
&e_n |{\bf m}\rangle =[m_n]|{\bf m}-{\bf e}_n \rangle, \notag\\
&f_n |{\bf m}\rangle =|{\bf m}+{\bf e}_n \rangle, \notag\\
&k_n |{\bf m}\rangle ={\bf v}q_n^{-m_n} |{\bf m}\rangle,\notag\\
&k_\La |{\bf m}\rangle ={\bf v}^{-\sum_im_i} |{\bf m}\rangle,\notag\\
\end{align*}}
for $\ket{\mathbf{m}}\in \W$ and $1\le i\le n-1$.

\item[(2)] If $X=D$ and $\e_1=1$, then the above action of the generators defines a $\U_{X\!B}(\e)$-module structure on $\W$ except that the actions of $e_0, f_0$ and $k_0$ are replaced by 
\begin{align*}
&e_{0} |{\bf m}\rangle = x|{\bf m}+{\bf e}_{1}+{\bf e}_{2} \rangle,\\
&f_{0} |{\bf m}\rangle = x^{-1}[m_2]|{\bf m}-{\bf e}_{1}-{\bf e}_{2} \rangle,\\
&k_{0} |{\bf m}\rangle =q_1^{m_1-1}q_2^{m_2}|{\bf m} \rangle, \\
&k_\La |{\bf m}\rangle ={\bf v}^{-\sum_im_i} |{\bf m}\rangle,\notag\\
\end{align*}
for $\ket{\mathbf{m}}\in \W$ and $1\le i\le n-1$. 
\end{itemize}
We denote the $\U_{XB}(\e)$-modules in (1) and (2) by $\W(x)$.
\end{prop}
\pf It can be checked in a straightforward manner.
\qed

\begin{rem}\label{rem:wt of m}
{\rm
We have $k_{\de_i}\ket{\bf m}={\bf v}^{-1}q_i^{m_i}\ket{\bf m}$ ($i\in \I$).
The $\U_{X\!B}(\e)$-module $\W(x)$ has a weight space decomposition in the sense of \eqref{eq:weight space}, where we have ${\rm wt}(\ket{{\bf m}})=\La +\sum_{i\in \I}m_i\de_i$.
In particular, $\W(x)$ is an irreducible highest weight $\ring{\U}_{X\!B}(\e)$-module with a highest weight vector $\ket{{\bf 0}}$ of weight $\La$.
}
\end{rem}

\begin{rem}\label{rem:classical limit}
{\rm
We will take a classical limit of $\W(x)$ as a $\ring{\U}_{X\!B}(\e)$-module to study a category of $q$-oscillator representations in Section~\ref{sec:osc rep}. 
In order to do that we should first view $\W(x)$ as a $\ring{U}_{X\!B}(\e)$-module via the isomorphism $\tau$ as in \cite{KLO24}.
However, the resulting $\ring{U}_{X\!B}(\e)$-module is still incompatible with classical limit, especially with respect to a subalgebra $U_{q^{1/2}}(\mf{sl}_2)=\langle\,E_n,F_n,K_n^{\pm 1}\,\rangle$. 
Indeed, we have 
\begin{equation*}
\begin{split}
(E_nF_n-F_nE_n)\ket{\bf m}
&=\frac{K_n-K_n^{-1}}{{\bf v}-{\bf v}^{-1}}\ket{\bf m}
= \frac{{\bf v}q^{-m_n}-{\bf v}^{-1}q^{m_n}}{{\bf v}-{\bf v}^{-1}}\ket{\bf m} \\
&= \frac{{\ib}q^{-\hf-m_n}+{\ib}q^{\hf+m_n}}{{\ib}q^{-\hf}+{\ib}q^\hf}\ket{\bf m} 
= \frac{q^{-\hf-m_n}+q^{\hf+m_n}}{q^{-\hf}+q^\hf}\ket{\bf m},
\end{split}
\end{equation*}
for $\ket{\bf m}\in \W$ (here $\e_n=0$). Hence $E_n, F_n$, and $\frac{K_n-K_n^{-1}}{{\bf v}-{\bf v}^{-1}}$ do not give rise to an $\mf{sl}_2$-action at $q^\hf=1$. 
This obstruction was first observed in \cite{KwO21} when $\e=(0^n)$ and resolved by further replacing the quantum $\mathfrak{sl}_{2}$ with the quantum supergroup associated with $\mf{osp}_{1|2}$
through the twistor automorphism of a covering quantum group \cite{CFLW}, which will be extended to general $\epsilon$ in Section~\ref{sec:covering QG} below.
}
\end{rem}

\section{Truncation functors}\label{sec:truncation}
Recall that we are dealing with two algebras $\mathcal{U}_{XB}(\epsilon)$ for $X=B,D$. 
To keep our discussion uniform, we will assume the following notations. First for $X=B$, we let 
\begin{itemize}
\item[$\bullet$] $\bm{\e}=(\e_1,\dots,\e_n)=(0,1,0\dots,1,0)$ with $n=2m+1$ ($m\ge 2$),

\item[$\bullet$] $\U_{\bm{\e}}=\U_{B\!B}(\bm{\e})$ with $I_{\bm{\e}}=I$,

\item[$\bullet$] $\un{\bm{\e}}=(0^{m+1})$, $\ov{\bm\e}=(1^m)$,

\item[$\bullet$] $I_{\un{\bm\e}}=\{\,0,1,\dots,m+1\,\}$, $I_{\ov{\bm\e}}=\{\,0,1,\dots,m\,\}$,

\item[$\bullet$] $\U_{\un{\bm\e}}=U_{q}(D_{m+2}^{(2)})=\langle\, e_i,f_i, k_i^{\pm 1}\,|\, i\in I_{\un{\bm\e}}\, \rangle$, 

\item[$\bullet$] $\U_{\ov{\bm\e}}=U_{\td{q}}(D_{m+1}^{(2)})=\langle\, e_i,f_i, k_i^{\pm 1}\,|\, i\in I_{\ov{\bm\e}}\, \rangle$.

\end{itemize}
When $X=D$, we instead set
\begin{itemize}
\item[$\bullet$] $\bm{\e}=(\e_1,\dots,\e_n)=(1,0\dots,1,0)$ with $n=2m$ ($m\ge 2$),

\item[$\bullet$] $\U_{\bm{\e}}=\U_{D\!B}(\bm{\e})$ with $I_{\bm{\e}}=I$,

\item[$\bullet$] $\un{\bm{\e}}=(0^m)$, $\ov{\bm\e}=(1^m)$,

\item[$\bullet$] $I_{\un{\bm\e}}=\{\,0,1,\dots,m\,\}=I_{\ov{\bm\e}}$,

\item[$\bullet$] $\U_{\un{\bm\e}}=U_{q}(A_{2m}^{(2)\dagger})=\langle\, e_i,f_i, k_i^{\pm 1}\,|\, i\in I_{\un{\bm\e}}\, \rangle$, 

\item[$\bullet$] $\U_{\ov{\bm\e}}=U_{\td{q}}(B_{m}^{(1)})=\langle\, e_i,f_i, k_i^{\pm 1}\,|\, i\in I_{\ov{\bm\e}}\, \rangle$,

\end{itemize}
where $A_{2m}^{(2)\dagger}$ signifies the opposite labeling of nodes of the Dynkin diagram of affine type $A_{2m}^{(2)}$.
In either case, we denote by $\mathcal{W}_{\epsilon}(x)$ the $\mathcal{U}_\epsilon$-module in Proposition~\ref{prop:osc W for e},
with the exception of $\mathcal{W}_{\underline{\boldsymbol{\epsilon}}}(x)\coloneqq \mf{tr}_{\underline{\boldsymbol{\epsilon}}}(\W_{\bm\e}(x))$ when $X=D$, where $\mathfrak{tr}_{\underline{\boldsymbol{\epsilon}}}$ is the truncation functor defined in Theorem~\ref{thm:truncation functor} below.

We shall construct functors sending $\U_{\bm\e}$-modules to $\U_{\e}$-modules ($\e=\un{\bm\e}, \ov{\bm\e}$) as in \cite{KL20,KL22,KLO24}. 
Indeed, it suffices to construct a homomorphism of $\Bbbk$-algebras from $\U_{\e}$ to $\U_{\bm\e}$.

We define  
$\check{e}_j, \check{f}_j, \check{k}_j \in \U_{\bm\e}$ for $j\in I_{\un{\bm\e}}$ and $\hat{e}_i, \hat{f}_i, \hat{k}_i\in \U_{\bm\e}$ for $i\in I_{\ov{\bm\e}}$ as follows. Let $\eta=\pm 1$, and recall $c=\frac{\mathbf{v}-\mathbf{v}^{-1}}{q^{1/2}-q^{-1/2}}$ \eqref{eq:scalar c}.
First assume $X=B$. We put
{\allowdisplaybreaks
\begin{equation}\label{eq:check{e} for b}
\begin{split}
 \check{e}_j&=
\begin{cases}
ce_0 & \text{if $j=0$},\\ 
[e_{2j-1},e_{2j}]_{-q^\eta} & \text{if $1\le j\le m$},\\
ce_{2m+1} & \text{if $j=m+1$},
\end{cases}\quad 
 \check{f}_j=
\begin{cases}
f_0 & \text{if $j=0$},\\ 
[f_{2j},f_{2j-1}]_{-q^{-\eta}} & \text{if $1\le j\le m$},\\
f_{2m+1} & \text{if $j=m+1$},
\end{cases}\\
 \check{k}_j&=
\begin{cases}
 k_0 & \text{if $j=0$},\\ 
 k_{2j-1}k_{2j} & \text{if $1\le j\le m$},\\
 k_{2m+1} & \text{if $j=m+1$},
\end{cases}
\end{split}
\end{equation}}
\begin{equation}\label{eq:hat{e} for b}
\begin{split}
\hat{e}_i=[e_{2i},e_{2i+1}]_{q^{\eta}},\quad
\hat{f}_i=[f_{2i+1},f_{2i}]_{q^{-\eta}},\quad
\hat{k}_i=k_{2i}k_{2i+1}.
\end{split}
\end{equation}
Next when $X=D$ (\textit{cf.} \cite[Proposition~3.1]{KLO24}),
{\allowdisplaybreaks
\begin{equation}\label{eq:check{e} for d}
\begin{split}
 \check{e}_j&=
\begin{cases}
[e_1,e_0]_{-q^{-2\eta}} & \text{if $j=0$},\\ 
[e_{2j},e_{2j+1}]_{-q^\eta} & \text{if $0< j < m$},\\
ce_{2m} & \text{if $j=m$},
\end{cases}\quad 
 \check{f}_j=
\begin{cases}
\frac{1}{[2]^{2\delta_{j,0}}}[f_0,f_1]_{-q^{2\eta}} & \text{if $j=0$},\\ 
[f_{2j+1},f_{2j}]_{-q^{-\eta}} & \text{if $0< j < m$},\\
f_{2m} & \text{if $j=m$},
\end{cases}\\
 \check{k}_j&=
\begin{cases}
 k_{2j}k_{2j+1} & \text{if $0\le j < m$},\\
 k_{2m} & \text{if $j=m$},
\end{cases}
\end{split}
\end{equation}
\begin{equation}\label{eq:hat{e} for d}
\begin{split}
 \hat{e}_j&=
\begin{cases}
[e_0,e_2]_{q^{\eta}} & \text{if $j=0$},\\ 
[e_{2j-1},e_{2j}]_{q^\eta} & \text{if $0< j \le m$},
\end{cases}\quad 
 \hat{f}_j=
\begin{cases}
[f_2,f_0]_{q^{-\eta}} & \text{if $j=0$},\\ 
[f_{2j},f_{2j-1}]_{q^{-\eta}} & \text{if $0< j \le m$},
\end{cases}\\
 \hat{k}_j&=
\begin{cases}
 k_0 k_2 & \text{if $j=0$},\\
 k_{2j-1}k_{2j} & \text{if $0<j\le m$}.
\end{cases}
\end{split}
\end{equation}
}

\begin{prop}\label{prop:reduction}
Suppose $X=B$ or $D$. For $\e=\un{\bm\e},\ov{\bm\e}$, there exists a homomorphism of $\Bbbk$-algebras $\phi_{\e} : \U_{\e} \longrightarrow \U_{{\bm \e}}$ such that 
\begin{itemize}
\item[(1)] $\phi_{\un{\bm\e}}(e_j)=\check{e}_j$, $
\phi_{\un{\bm\e}}(f_j)=\check{f}_j$, and $
\phi_{\un{\bm\e}}(k_j)=\check{k}_j$ for $j\in I_{\un{\bm\e}}$, when $\e=\un{\bm \e}$,

\item[(2)] $\phi_{\ov{\bm\e}}(e_i)=\hat{e}_i$, $
\phi_{\ov{\bm\e}}(f_i)=\hat{f}_i$, and $\phi_{\ov{\bm\e}}(k_i)=\hat{k}_i$  for $i\in I_{\ov{\bm\e}}$, when $\e=\ov{\bm \e}$.
\end{itemize}
\end{prop}
\pf
We shall verify that $\phi_\epsilon (x_j)$ ($x=e,f,k$) satisfy the defining relations of $\mathcal{U}_{\bm\e}$. 
We give a proof of the most involved one:
\begin{equation}\label{eq:Serre-hat{e}-B}
(X=B)\quad \hat{e}_{0}^{3} \hat{e}_1 -[3]_{\mathbf{v}}(\hat{e}_{0}^2 \hat{e}_{1} \hat{e}_{0} -\hat{e}_{0} \hat{e}_{1} \hat{e}_{0}^{2})-\hat{e}_{1}\hat{e}_{0}^{3} =0
\end{equation}
in Appendix \ref{sec:tr-pf-Serre}, and omit the others.
One may refer to \cite{KLO24,KY} for the Serre relations of type $D$ and $A$ respectively, as the verification of each relation is local.
\qed
\smallskip

\begin{rem}{\rm
Let $\check{k}_\Lambda$ and $\hat{k}_\Lambda$ be the elements in $\U_{\un{\boldsymbol{\e}}}$ and $\U_{\ov{\boldsymbol{\e}}}$ associated to $\La$. Then we define $\phi_{\un{\bm\e}}(\hat{k}_\Lambda)={k}_\Lambda$ and $\phi_{\ov{\bm\e}}(\check{k}_\Lambda)={k}_\Lambda$.
} 
\end{rem}

Let 
\begin{itemize}
\item[$\bullet$] $\I_{\bm\e}=\I$, $\I_{\ov{\bm\e}}=\{2,4,\dots,2m\}$, $\I_{\un{\bm\e}}=\mathbb{I}_{\mathbf{\epsilon}}\setminus\mathbb{I}_{\mathbf{\overline{\epsilon}}}$,

\item[$\bullet$] $P_{\e}= \Z\La \oplus \bigoplus_{a\in \I_\e}\Z\de_a$, 
 $P_{\ge 0,\e}=\Z\La \oplus \bigoplus_{a\in \I_\e}\Z_+\de_a$\ \ for $\e={\bm\e},\un{\bm\e},\ov{\bm\e}$,
 
\item[$\bullet$] $\{\,\alpha_i^\e\,|\,i\in I_\e\,\}$: the set of simple roots for $\U_\e$\ \ for $\e={\bm\e},\un{\bm\e},\ov{\bm\e}$, 
where we may identify $\alpha_i^\e \in P_{\bm\e}$ ($\e=\un{\bm\e},\ov{\bm\e}$) as follows:
\begin{align*}
(X=B)\quad & \alpha^{\un{\bm\e}}_i=
\begin{cases}
 \de_1 & \text{if $i=0$},\\
 -\de_{2i-1}+\de_{2i+1} & \text{if $1\le i\le m$},\\
 -\de_{2m+1} & \text{if $i=m+1$},
\end{cases}
\quad
\alpha^{\ov{\bm\e}}_i=
\begin{cases}
 \de_2 & \text{if $i=0$},\\
 -\de_{2i}+\de_{2i+2} & \text{if $1\le i<m$},\\
 -\de_{2m} & \text{if $i=m$}.
\end{cases} 
\\
(X=D)\quad & \alpha^{\un{\bm\e}}_i=
\begin{cases}
 2\de_2 & \text{if $i=0$},\\
 -\de_{2i}+\de_{2i+2} & \text{if $1\le i < m$},\\
 -\de_{2m} & \text{if $i=m$},
\end{cases}
\quad
\alpha^{\ov{\bm\e}}_i=
\begin{cases}
 \de_1+\de_3 & \text{if $i=0$},\\
 -\de_{2i-1}+\de_{2i+1} & \text{if $1\le i<m$},\\
 -\de_{2m-1} & \text{if $i=m$}.
\end{cases}
\end{align*}

\item[$\bullet$] $\mc{C}(\e)$ ($\e={\bm\e},\un{\bm\e},\ov{\bm\e}$): the category of $\U_{\e}$-modules $V$ such that 
\begin{itemize}
 \item[(i)] $V=\bigoplus_{\la\in P_\e}V_\la$ with
$V_\la=\{\,v\,|\,k_\mu v = {\bq}(\la,\mu)v\ (\mu\in P_\e)\,\}$, 
 \item[(ii)] ${\rm wt}(V)\subset P_{\ge 0,\e}$,
 
\end{itemize}

\item[$\bullet$] $\mc{C}^\ell(\e)$: the subcategory of $\mc{C}(\e)$ with objects $V$ such that
${\rm wt}(V)\subset \ell\La + \sum_{a\in \I_\e}\Z_+\de_a$.
\end{itemize}

Note that $\mc{C}(\e)=\bigoplus_{\ell\in \mathbb{Z}}\mc{C}^\ell(\e)$ is closed under taking submodules, quotients and tensor products, where $\mc{C}^\ell(\e)\ot \mc{C}^{\ell'}(\e)\subset \mc{C}^{\ell+\ell'}(\e)$, and $\W(x)^{\otimes \ell}\in \mc{C}^{\ell}(\e)$ ($x\in \Bbbk^\times$, $\ell\ge 1$) (\emph{cf.} Remark \ref{rem:wt of m}).
Let $\mathrm{pr}:P_{\boldsymbol{\epsilon}}\rightarrow \bigoplus\mathbb{Z}\delta_a$ be the projection $\mathrm{pr}(\Lambda)=0$, $\mathrm{pr}(\delta_a)=\delta_a$.
Then we have the following, which is an analogue of \cite[Section 3]{KLO24}.
\begin{thm}\label{thm:truncation functor}
For $\e=\un{\bm\e},\ov{\bm\e}$, there exists an exact monoidal functor 
\begin{equation*}\label{eq:tr}
\xymatrixcolsep{2pc}\xymatrixrowsep{3pc}
\xymatrix{
\mf{tr}_\e=\mf{tr}^{\bm\e}_{\e} : 
\mc{C}({\bm\e})  \ar@{->}[r] &  \mc{C}(\e)} 
\end{equation*}
such that for $V,W\in \mc{C}({\bm\e})$ and $f \in {\rm Hom}_{\mc{C}({\bm\e})}(V,W)$
\begin{itemize}
\item[(1)] $\mf{tr}_{\e}(V) = 
\bigoplus_{\mu}V_\mu$ with $\mu\in {\rm wt}(V)$ such that $({\rm pr}(\mu)|\de_{a})=0 $ for all $a\not\in \mathbb{I}_\e$,

\item[(2)] $\mf{tr}_\e(f) \in {\rm Hom}_{\mc{C}(\e)}(\mf{tr}_\e(V),\mf{tr}_\e(W))$ is given by $\mf{tr}_\e(f)(v)=f(v)$ for $v\in \mf{tr}_\e(V)$,
\end{itemize}
where the $\U_{\e}$-action on $\mf{tr}_{\e}(V)$ is induced from $\phi_\epsilon$ in Proposition \ref{prop:reduction}. We call $\mf{tr}_{\e}$ a truncation functor.
\end{thm}
\pf It can be checked by the same arguments as in \cite[Propositions 4.4]{KY} thanks to the condition ${\rm wt}(V)\subset P_{\ge 0,{\bm\e}}$.
\qed\smallskip

\begin{rem}\label{rem:truncation of W}
{\rm 
Recall the $\mathcal{U}_\epsilon$-module $\mathcal{W}_{\epsilon}(x)$. We have
\begin{itemize}
\item[(1)] $\W_{\e}(x)\cong \mf{tr}_{\e}(\W_{\bm\e}(x))$ for $\e=\un{\bm\e},\ov{\bm\e}$,

\item[(2)] $\W_{\ov{\bm\e}}(x)$ is isomorphic to the $q$-spin representation (the $m$-th fundamental representation) of $U_{\td{q}}(D_{m+1}^{(2)})$ ($X=B$) or $U_{\td{q}}(B_{m}^{(1)})$ ($X=D$),

\item[(3)] $\W_{\un{\bm\e}}(x)$ is isomorphic to the $q$-oscillator representation of $U_q(D_{m+2}^{(2)})$ ($X=B$) or $U_{\td{q}}(A_{2m}^{(2)\dagger})$ ($X=D$) in \cite{KO}.

\end{itemize}
} 
\end{rem}

\begin{rem}{\rm 
Let $\e^c=(1-\e_1,\dots,1-\e_n)$.
We may consider the case of
\begin{equation*}\label{eq:e compliment}
({\bm\e}',\un{\bm\e}',\ov{\bm\e}'):=({\bm\e}^c,\un{\bm\e}^c,\ov{\bm\e}^c) 
\end{equation*}
with $\U_\e$ for $\e={\bm\e},\un{\bm\e},\ov{\bm\e}$ in Section \ref{sec:truncation} replaced by ${\bm\e}',\un{\bm\e}',\ov{\bm\e}'$.
For example, when $X=B$, the homomorphism in Proposition \ref{prop:reduction} is replaced by $\phi_{\e} : \U_{\e} \longrightarrow \U_{{\bm \e}'}$ for $\e=\un{\bm\e}',\ov{\bm\e}'$ such that 
$\phi_{\un{\bm\e}}(x_i)=\hat{x}_i$ and $\phi_{\ov{\bm\e}}(x_j)=\hat{x}_j$ for $i\in I_{\un{\bm\e}}$, $j\in I_{\ov{\bm\e}}$ and $x=e,f,k$, where
\begin{equation}\label{eq:hat{e}' for b}
\begin{split}
\hat{e}_i=c^{\delta_{i0}+\delta_{im}}[e_{2i},e_{2i+1}]_{-q^{\eta}},\quad \hat{f}_i=[f_{2i+1},f_{2i}]_{-q^{-\eta}},\quad \hat{k}_i=k_{2i}k_{2i+1},
\end{split}
\end{equation}
\begin{equation}\label{eq:check{e}' for b}
\begin{split}
 \check{e}_j&=
\begin{cases}
e_0 & \text{if $j=0$},\\ 
[e_{2j-1},e_{2j}]_{q^\eta} & \text{if $1\le j\le m$},\\
e_{2m+1} & \text{if $j=m+1$},
\end{cases}\ \
 \check{f}_j=
\begin{cases}
f_0 & \text{if $j=0$},\\ 
[f_{2j},f_{2j-1}]_{q^{-\eta}} & \text{if $1\le j\le m$},\\
f_{2m+1} & \text{if $j=m+1$},
\end{cases}\\
 \check{k}_j&=
\begin{cases}
k_0 & \text{if $j=0$},\\ 
k_{2j-1}k_{2j} & \text{if $1\le j\le m$},\\
k_{2m+1} & \text{if $j=m+1$},
\end{cases}
\end{split} 
\end{equation}
where $\eta=\pm1$ and $c$ is as in \eqref{eq:scalar c}.
} 
\end{rem}

\section{Covering quantum group}\label{sec:covering QG}
\subsection{Orthosymplectic Lie superalgebras}

Let $V=V_{\ov 0}\oplus V_{\ov 1}$ be a finite-dimensional complex superspace with a non-degenerate skew supersymmetric bilinear form $\mathtt{B}$. 
Recall that $\mf{spo}(V)$ is defined to be the Lie sub-superalgebra of $\mathfrak{gl}(V)$ with $\mathfrak{spo}(V)_{\varepsilon}=\{\,T\,|\,\texttt{B}(T(v),w)+(-1)^{\varepsilon|v|}\texttt{B}(v,T(w))=0 \, (v,w\in V) \}$ ($\varepsilon=\ov{0},\ov{1}$), where $|v|$ denote the degree of a homogeneous vector $v$. 
Similarly, $\mf{osp}(W)$ is defined with respect to a non-degenerate supersymmetric bilinear form on a superspace $W$. 
Note that $\mf{spo}(V)\cong\mf{osp}(\Pi V)$, where $\Pi$ is the parity change functor and $\mf{osp}(\Pi V)$ is defined with respect to a supersymmetric bilinear form ${\mathtt B}^{\Pi}$ given by ${\mathtt B}^\Pi(\Pi(v),\Pi(w))=(-1)^{|v|}{\mathtt B}(v,w)$ for homogeneous vectors $v,w$ (cf.\cite{CW}).

Let $\e$ be given.
We put
\begin{equation}\label{eq:osp(e) type B}
 \mf{osp}_B(\e)=
\begin{cases}
 \mf{spo}(V) & \text{if $\e_n=1$},\\
 \mf{osp}(V') & \text{if $\e_n=0$}.
\end{cases}
\end{equation}
where $V$ (resp. $V'$) is the superspace such that $(\dim V_{\ov 0},\dim V_{\ov 1})=(2M,2N+1)$ (resp. $(\dim V'_{\ov 0},\dim V'_{\ov 1})=(2M+1,2N)$) with $M=|\{\,i\,|\,\e_i=0\,\}|$ and $N=|\{\,i\,|\,\e_i=1\,\}|$.

Here we take a fundamental system for $\mf{spo}(V)$ (the set of positive simple roots for $\mf{spo}(V)$) to be the one associated with $\e$, namely  the one corresponding to $\upvarepsilon \updelta$ sequence \cite[Section~1.3]{CW} replacing $0$ and $1$ in $\e$ with $\updelta$ and $\upvarepsilon$, respectively. 
Similarly, the fundamental system of $\mf{osp}(V')$ is understood to be the one for $\mf{spo}(V')$, where the $\upvarepsilon \updelta$ sequence now corresponds to $\e^c=(1-\e_1,\dots,1-\e_n)$.
In particular, the fundamental system for $\mf{osp}_B(\e)$ does not have an odd non-isotropic simple root by definition \eqref{eq:osp(e) type B}.

We also put
\begin{equation}\label{eq:osp(e) type B super}
 \mf{osp}_{B^\bullet}(\e)=
\begin{cases}
 \mf{spo}(V') & \text{if $\e_n=0$},\\
 \mf{osp}(V) & \text{if $\e_n=1$}.
\end{cases}
\end{equation}
Here the fundamental system for $\mf{spo}(V)$ with $\e_n=0$ is also chosen to be the one corresponding to $\e$ following \cite[Section 1.3]{CW} so that the simple root $\alpha_n$ of $\mf{osp}_{B^\bullet}(\e)$ is non-isotropic odd.
For example, the associated Dynkin diagrams are  
\medskip

\begin{center}\setlength{\unitlength}{0.16in} \smallskip
\begin{picture}(18,3)
\put(0,2){\makebox(0,0)[c]{$\mf{osp}_{B}(00011)$ : }}
\put(5,2){\makebox(0,0)[c]{$\bigcirc$}}
\put(8,2){\makebox(0,0)[c]{$\bigcirc$}}
\put(11,1.95){\makebox(0,0)[c]{$\bigotimes$}}
\put(14,2){\makebox(0,0)[c]{$\bigcirc$}}
\put(17,2){\makebox(0,0)[c]{$\bigcirc$}}

\put(5.35,2){\line(1,0){2.1}} 
\put(8.4,2){\line(1,0){2.1}} 
\put(11.4,2){\line(1,0){2.1}}
\put(14.4,1.9){\line(1,0){2}}
\put(14.4,2.1){\line(1,0){2}}
\put(16.3,2){\makebox(0,0)[c]{\Large $>$}}

\put(5,1){\makebox(0,0)[c]{\tiny $1$}}
\put(8.2,1){\makebox(0,0)[c]{\tiny $2$}}
\put(11,1){\makebox(0,0)[c]{\tiny $3$}}
\put(14,1){\makebox(0,0)[c]{\tiny $4$}}
\put(17,1){\makebox(0,0)[c]{\tiny $5$}}
\end{picture}
\end{center}

\begin{center}\setlength{\unitlength}{0.16in} \smallskip
\begin{picture}(18,3)
\put(0,2){\makebox(0,0)[c]{$\mf{osp}_{B}(01010)$ : }}
\put(5,2){\makebox(0,0)[c]{$\bigotimes$}}
\put(8,2){\makebox(0,0)[c]{$\bigotimes$}}
\put(11,1.95){\makebox(0,0)[c]{$\bigotimes$}}
\put(14,2){\makebox(0,0)[c]{$\bigotimes$}}
\put(17,2){\makebox(0,0)[c]{$\bigcirc$}}

\put(5.35,2){\line(1,0){2.1}} 
\put(8.4,2){\line(1,0){2.1}} 
\put(11.4,2){\line(1,0){2.1}}
\put(14.4,1.9){\line(1,0){2}}
\put(14.4,2.1){\line(1,0){2}}
\put(16.3,2){\makebox(0,0)[c]{\Large $>$}}

\put(5,1){\makebox(0,0)[c]{\tiny $1$}}
\put(8.2,1){\makebox(0,0)[c]{\tiny $2$}}
\put(11,1){\makebox(0,0)[c]{\tiny $3$}}
\put(14,1){\makebox(0,0)[c]{\tiny $4$}}
\put(17,1){\makebox(0,0)[c]{\tiny $5$}}
\end{picture}
\end{center}

\begin{center}\setlength{\unitlength}{0.16in} \smallskip
\begin{picture}(18,3)
\put(0,2){\makebox(0,0)[c]{$\mf{osp}_{B^\bullet}(01010)$ : }}
\put(5,2){\makebox(0,0)[c]{$\bigotimes$}}
\put(8,2){\makebox(0,0)[c]{$\bigotimes$}}
\put(11,1.95){\makebox(0,0)[c]{$\bigotimes$}}
\put(14,2){\makebox(0,0)[c]{$\bigotimes$}}
\put(17,2){\makebox(0,0)[c]{\scalebox{2.3}{$\bullet$}}}

\put(5.35,2){\line(1,0){2.1}} 
\put(8.4,2){\line(1,0){2.1}} 
\put(11.4,2){\line(1,0){2.1}}
\put(14.4,1.9){\line(1,0){2}}
\put(14.4,2.1){\line(1,0){2}}
\put(16.3,2){\makebox(0,0)[c]{\Large $>$}}

\put(5,1){\makebox(0,0)[c]{\tiny $1$}}
\put(8.2,1){\makebox(0,0)[c]{\tiny $2$}}
\put(11,1){\makebox(0,0)[c]{\tiny $3$}}
\put(14,1){\makebox(0,0)[c]{\tiny $4$}}
\put(17,1){\makebox(0,0)[c]{\tiny $5$}}
\end{picture}
\end{center}
where \scalebox{.8}{$\bigotimes$} and \raisebox{-0.7ex}{\scalebox{2}{$\bullet$}} denote the isotropic and non-isotropic odd simple roots, respectively. In particular, when $\e$ is homogeneous, that is, $\e=(0^n)$ or $(1^n)$, we have
\begin{equation*}
\mf{osp}_{B}(\e)=\mf{so}_{2n+1},\quad \mf{osp}_{B^{\bullet}}(\e)=\mf{osp}_{1|2n}.
\end{equation*}

\subsection{Covering quantum group and twistor}\label{subsec: covering QG and twistor}
\subsubsection{The generalized quantum groups $\U_{B}(\e)$ and $\U_{B^\bullet}(\e)$}\label{subsec: GQG of finite types}
Given $\e=(\e_1,\dots,\e_n)$, let {$\mathcal{U}_{B}(\epsilon)$} be the associative $\Bbbk$-algebra generated by 
$e_i,f_i,k_\mu$  $(i\in \ring{I}, \mu\in P)$ with the same relations given in Section \ref{subsec:def}. We may identify it with $\ring{\U}_{XB}(\e)$, and also regard it as a quantized enveloping algebra of $\mf{osp}_B(\e)$ \eqref{eq:osp(e) type B}. 

We define {$\U_{B^\bullet}(\e)$} to be the associative $\Bbbk$-algebra generated by 
$e_i,f_i,k_\mu$  $(i\in \ring{I}, \mu\in P)$ satisfying the following relations:

\begin{gather*}
k_{\mu}=1 \quad(\mu=0),\quad k_{\mu+\mu'}=k_{\mu}k_{\mu'} \quad (\mu, \mu' \in P),\\
k_\mu e_i k_{-\mu}=\bq(\mu,\alpha_i)e_i,\quad 
k_\mu f_i k_{-\mu}=\bq(\mu,\alpha_i)^{-1}f_i\quad (i\in \ring{I}, \mu\in P), \\
[e_i ,f_j]= \delta_{ij} \displaystyle \frac{k_i-k^{-1} _i}{q-q^{-1}} \quad\text{if $(i,j)\ne(n,n)$}, \\
e_nf_n+f_ne_n= \displaystyle \frac{k_n-k^{-1} _n}{q^{\hf}-q^{-\hf}},\\
e^2 _i =f^2 _i =0 \quad \text{if $(\alpha_i|\alpha_i)=0$},\\
e_ie_j -e_je_i=(e\rightarrow f)=0 \quad \text{if $(\alpha_i|\alpha_j)=0$ ($i\neq j$)},\\
e^2 _i e_j - (-1)^{\epsilon_i}[2]e_i e_j e_i +e_j e^2 _i =(e \rightarrow f)=0
\quad \text{if $i\in \ring{I}\setminus\{n\}$, $|i-j|=1$, $\e_i=\e_{i+1}$}, \\
{e_n^3e_{n-1}+(1-(-1)^{\e_n}[2])(e_n^2e_{n-1}e_n+e_ne_{n-1}e_n^2)+e_{n-1}e_n^3=(e \rightarrow f)=0,} \\
e_{i}e_{i-1}e_{i}e_{i+1}-e_{i}e_{i+1}e_{i}e_{i-1}+(-1)^{\epsilon_i}[2]e_{i}e_{i-1}e_{i+1}e_{i} \\
\hspace{0.1cm}-e_{i-1}e_{i}e_{i+1}e_{i}+e_{i+1}e_{i}e_{i-1}e_{i}=(e \rightarrow f)=0
\quad\text{if $i\in \{2,\dots,n-1\}$, $\e_i\neq \e_{i+1}$}. \nonumber
\end{gather*}
It is a variant of the quantized enveloping algebra of $\mf{osp}_{B^\bullet}(\e)$ \eqref{eq:osp(e) type B super} \cite{Ya94}, which is indeed a superalgebra with odd generators $e_n$ and $f_n$.
In particular, $\U_{B^\bullet}(0^n)\cong U_q(B(0,n))=U_q(\mf{osp}_{1|2n})$ is the quantum supergroup of $\mf{osp}_{1|2n}$ which we recall soon.

Recall the isomorphism $\tau$ in Section~~\ref{subsec:Yamane QSA}.
Let $U_B(\e)$ be defined similarly using the relations of $U_{XB}(\epsilon)$.
Then it can be identified with a subalgebra $\mathring{U}_{X\!B}(\e)$,
and $\tau$ naturally restricts to an isomorphism between the extensions $U_B(\e)[\sigma]$ and $\U_B(\e)[\sigma]$ by the bialgebra $\Sigma$ (Remark~\ref{rem: iso GQG to YQG}).

We have an analogue for $\U_{B^\bullet}(\e)$. Let $U_{B^\bullet}(\e)$ be the quantized enveloping algebra of $\mf{osp}_{B^\bullet}(\e)$ given in \cite{Ya94}.
Namely, $U_{B^\bullet}(\e)$ is the $\Z_2$-graded $\Bbbk$-algebra generated by 
$E_i,F_i,K_\mu$  $(i\in \ring{I}, \mu\in P)$ satisfying the following relations:
\allowdisplaybreaks{
\begin{gather*}
  K_{\mu}=1 \quad(\mu=0),\quad K_{\mu+\mu'}=K_{\mu}K_{\mu'} \quad (\mu, \mu' \in P),\\
  K_{\mu}E_{i}K_{\mu}^{-1}=q^{(\mu|\alpha_{i})}E_{i},\quad K_{\mu}F_{i}K_{\mu}^{-1}=q^{-(\mu|\alpha_{i})}F_{i} \quad (i\in I, \mu\in P),\\
  E_{i}F_{j}-(-1)^{p(i)p(j)}F_{j}E_{i}=(-1)^{\epsilon_{i}}\delta_{ij}\frac{K_{i}-K_{i}^{-1}}{q-q^{-1}} \quad \text{if  $(i,j)\ne (n,n)$},\\
{E_{n}F_{n}+F_{n}E_{n}= \frac{K_n-K_n^{-1}}{q^\hf-q^{-\hf}}},\\  
  E_{i}^{2}=F_{i}^{2}=0\quad \text{if $(\alpha_i|\alpha_i)=0$},\\
  E_{i}E_{j}-(-1)^{p(i)p(j)}E_{j}E_{i}=(E\rightarrow F)=0\quad\text{if $(\alpha_i|\alpha_j)=0$ ($i\neq j$)},\\
  E_{i}^{2}E_{j}-[2]E_{i}E_{j}E_{i}+E_{j}E_{i}^{2}=(E\rightarrow F)=0\quad \text{if $i,j\in \mathring{I}\setminus\{n\}$, $|i-j|=1$, $\e_i=\e_{i+1}$},\\
{\left[E_n,\left[E_n,\left[E_n,E_{n-1}\right]_{(-1)^{p(n-1)}q}\right]_{(-1)^{1+p(n-1)}}\right]_{(-1)^{p(n-1)}q^{-1}}=(E\rightarrow F)=0},\\ 
\left[E_{i},\left[\left[E_{i-1},E_{i}\right]_{(-1)^{p(i-1)}q},E_{i+1}\right]_{(-1)^{(p(i-1)+1)p(i+1)}q^{-1}}\right]_{(-1)^{p(i-1)+p(i)+p(i+1)}}\\
  \hfill=(E\rightarrow F)=0\quad\text{if $i\in \{2,\dots,n-1\}$, $\e_i\neq \e_{i+1}$}
\end{gather*}}
where we recall $[x,y]_t = xy-tyx$.
The parity  $p(i)$ of $E_i$ and $F_i$ is equal to the one in Section \ref{subsec:Yamane QSA} for $i\in \ring{I}\setminus\{n\}$ but $p(n)=1$.

\begin{prop}\label{prop:tau iso bullet}
Suppose that $\e=(\e_1,\dots,\e_n)=(0,1,0,\dots,0,1,0)$ or $(1,0,\dots,0,1,0)$.
There exists an isomorphism of $\Bbbk$-algebras
$\tau : \U_{B^\bullet}(\e)[\sigma] \longrightarrow U_{B^\bullet}(\e)[\sigma]$
such that
{\allowdisplaybreaks
\begin{equation*}\label{eq:iso phi}
\begin{split}
& \tau(\sigma_a) =\sigma_a,\quad \tau(k_{\delta_i})=K_{\delta_i}\sigma_i,\quad {\tau(k_{\La})=K_{\La}{\bm\sigma}},\\
& \tau(e_i)= 
\begin{cases}
E_i \sigma_{\le n} & \text{for $i=n$},\\
E_i \sigma_{\le i}\varsigma_i^{i} & \text{for $1\leq i\leq n-1$},
\end{cases}\\
& \tau(f_i)= 
\begin{cases}
F_i \sigma_{\le n} & \text{for $i=n$},\\
(-1)^{i}F_i \sigma_{\le i}\varsigma_i^{i+1} & \text{for $1\leq i\leq n-1$},
\end{cases}
\end{split}
\end{equation*}}
where $\sigma_{\leq i}=\sigma_{1}\sigma_{2}\cdots\sigma_{i}$ and $\varsigma_i=\sigma_i\sigma_{i+1}$ for $1\le i\le n-1$.
\end{prop}
\pf The proof is similar to that of Proposition \ref{prop:tau iso}.
\qed

\subsubsection{Covering quantum groups of $\U_{B}(\e)$ and $\U_{B^\bullet}(\e)$}
We introduce the notion of a covering quantum group for $\U_{B}(\e)$ and $\U_{B^\bullet}(\e)$ and its twistor which are inspired by those of quantum Kac-Moody superalgebras without isotropic simple roots \cite{CFLW}.

Let $\e=(\e_1,\dots,\e_n)$ be given. We keep the notations in Section \ref{subsec:def},
 and let
\begin{itemize}
\item[$\bullet$] $(d_1,\dots,d_n)=(2,2,\ldots,2,1)$,

\item[$\bullet$] $\mathtt{q}_i={\bf i}^{d_i\e_i}q^{d_i/2}$ for $i\in \ring{I}=\{1,\dots,n\}$,

\item[$\bullet$] $(a_{ij})_{i,j\in \ring{I}}=((-1)^{\e_i}2(\alpha_i|\alpha_j)/d_i)_{i,j\in \ring{I}}$, equivalently $d_ia_{ij}=(-1)^{\e_i}2(\alpha_i|\alpha_j)$,

\item[$\bullet$] $R^\pi=R[\pi]/(\pi^2-1)$ for a subring $R$ of $\Bbbk$ with $1$ and an indeterminate $\pi$,

\item[$\bullet$] $\pi_i=\pi^{d_i}$ for $i\in \ring{I}$, that is, $(\pi_1,\dots,\pi_n)=(1,1,\dots,1,\pi)$,

\item[$\bullet$] $(\mathtt{p}(1),\dots,\mathtt{p}(n))=(0,0,\ldots,0,1)$.

\end{itemize}

\begin{rem}{\rm
The matrix $(a_{ij})_{i,j\in \ring{I}}$ is the Cartan matrix of type $B_n$ when $\e$ is homogeneous.
Note that $\mathtt{p}(i)$ should not be confused with the parity $p(i)$ of $E_i$ and $F_i$ for $i\in \ring{I}$ for $U_{B^\bullet}(\e)$. Also, note that
\begin{equation*}
\mathtt{q}_i=
\begin{cases}
q & \text{if $i\neq n$ with $\e_i=0$},\\
\td{q}^{-1} & \text{if $i\neq n$ with $\e_i=1$}, \\
q^\hf & \text{if $i=n$ with $\e_n=0$},\\
-{\bf v}^{-1} & \text{if $i=n$ with $\e_n=1$}. 
\end{cases}
\end{equation*}

} 
\end{rem}

We now introduce analogues of the covering quantum group for $\U_B(\e)$ and $\U_{B^\bullet}(\epsilon)$ in order to resolve the problem of taking a classical limit of $\U_B(\e)$-modules in Remark~\ref{rem:classical limit}. 
It turns out that these analogues can be defined by a natural extension of the construction in \cite{CFLW}, where only the case without isotropic odd simple roots is considered.

We first define the associated covering quantum group $\mathbf{U}$ to be the $\Q(q^{\hf},\ib)^\pi$-algebra with generators $e_i,f_i,J_i, k_{\mu}$ for 
$i\in \ring{I}$ and $\mu\in P$ subject to the following relations:
\begin{gather*}
k_{\mu}=1 \quad(\mu=0),\quad J_i k_\mu=k_\mu J_i, \quad k_{\mu+\mu'}=k_{\mu}k_{\mu'} \quad (\mu, \mu' \in P),\\
k_{\mu}e_ik_{-\mu}=\bq(\mu,\alpha_i)e_i,\quad k_\mu f_i k_{-\mu}=\bq(\mu,\alpha_i)^{-1}f_i \quad  (i\in \ring{I}, \mu\in P),\\
J_ie_j = \pi^{d_i a_{ij}} e_jJ_i,  \quad J_if_j =  \pi^{-d_i a_{ij}} f_jJ_i \quad  (i\in \ring{I}), \\
e_if_j-f_je_i=\delta_{ij}\frac{J_i k_i-k_i^{-1}}{q-q^{-1}}\quad\text{if $(i,j)\ne(n,n)$},\\
{e_nf_n-\pi f_ne_n=\frac{J_n k_n-k_n^{-1}}{\pi q^{1/2}-q^{-1/2}},}\\
{e^{2}_i=f^{2}_i=0\quad\text{if $(\alpha_i|\alpha_i)=0$},}\\
e_ie_j-e_je_i=(e\rightarrow f)=0\quad\text{if $(\alpha_i|\alpha_j)=0$},\\
{\sum_{l=0}^{1+|a_{ij}|}(-1)^l\pi^{\tfrac{l(l-1)\mathtt{p}(i)}{2}}
{1+|a_{ij}|\brack l}_{{\mathtt{q}_i,\pi_i}}e_i^{1+|a_{ij}|-l}e_je_i^l=(e\rightarrow f)=0}\\
\hspace{6cm}\text{if $i\neq j$ where $i\in \ring{I}\setminus\{n\}$ with $\epsilon_i=\epsilon_{i+1}$ or $i=n$},\\
e_{i}e_{i-1}e_{i}e_{i+1}-e_{i}e_{i+1}e_{i}e_{i-1}+(-1)^{\epsilon_i}[2]e_{i}e_{i-1}e_{i+1}e_{i} \\
\hspace{0.1cm}-e_{i-1}e_{i}e_{i+1}e_{i}+e_{i+1}e_{i}e_{i-1}e_{i}=(e \rightarrow f)=0
\quad\text{if $i\in \{2,\dots,n-1\}$ with $\e_i\neq \e_{i+1}$}.
\end{gather*}
Here we put $[m]_{z,\e}=\frac{(\e z)^{m}-z^{-m}}{\e z-z^{-1}}$ for $m\ge 0$ and define the binomial coefficient accordingly. Note that each $J_i$ (which is in fact $J_{d_i i}$ in the notation of \cite{CFLW}) is central as $d_i a_{ij}\in 2\mathbb{Z}$.

To define the twistor automorphism, we extend the algebra $\mathbf{U}$ to $\widehat{\mathbf{U}}$ by adjoining $T_i^{\pm 1},\Upsilon_i^{\pm 1}$ ($i\in \ring{I}$) and $T_\Lambda^{\pm 1}$ satisfying the following relations:
\begin{gather*}
T_iT_j=T_jT_i,\quad \Upsilon_i\Upsilon_j=\Upsilon_j\Upsilon_i,\\
T_i\Upsilon_j=\Upsilon_j T_i,\quad T_ik_\mu=k_\mu T_i, \quad \Upsilon_ik_\mu=k_\mu \Upsilon_i,\\
T_\Lambda T_i = T_i T_\Lambda,\quad T_\La\Upsilon_i=\Upsilon_i T_\La,\quad T_\La k_\mu=k_\mu T_\La,\\
T_ie_j=\ib^{d_ia_{ij}}e_jT_i,\quad T_if_j=\ib^{-d_ia_{ij}}f_jT_i,\\
{T_\Lambda e_j = \ib^{2(\La|\alpha_j)} e_j T_\Lambda,\quad
T_\Lambda f_j = \ib^{-2(\La|\alpha_j)} f_j T_\Lambda}\\
\Upsilon_ie_j=\ib^{\phi_{ij}}e_j\Upsilon_i,\quad \Upsilon_if_j=\ib^{-\phi_{ij}}f_j\Upsilon_i,
\end{gather*}
where 
\[
\phi_{ij}=\begin{cases}
d_ia_{ij}&\text{if }i>j,\\
d_i&\text{if }i=j,\\
0 &\text{if }i<j.
\end{cases}
\]

\begin{prop}\label{prop: twistor}
There is a $\Q(\ib)$-algebra automorphism $\widehat{\Psi}$ on $\widehat{\mathbf{U}}$
such that
\begin{alignat*}{4}
&e_i\mapsto \ib^{-d_i}\Upsilon_i^{-1}T_ie_i,&\quad
&f_i\mapsto f_i\Upsilon_i,&\quad
&k_i\mapsto T_i^{-1}k_i,& \quad 
&{k_\La\mapsto T_\La^{-1}k_\La},\\
&T_i\mapsto T_i,&\quad
&{T_\La\mapsto T_\La},&\quad
&\Upsilon_i\mapsto \Upsilon_i,&\quad 
&J_i \mapsto T_i^2 J_i,\\
&q^{\hf}\mapsto \ib^{-1}q^{\hf}={\bf v}^{-1},&\quad
&\pi\mapsto -\pi,
\end{alignat*}
for $i\in \ring{I}$.
\end{prop}
\pf
We shall check that $\widehat{\Psi}$ respects the following relations
\begin{align*}
\sum_{l=0}^{3} & (-1)^l \pi^{\tfrac{l(l-1)}{2}}{3\brack l}_{{\mathtt{q}_n,\pi}}e_n^{3-l}e_{n-1} e_n^l \\
&=e_n^3 e_{n-1} -[3]_{\mathtt{q}_n,\pi}e_n^2 e_{n-1}e_n +\pi[3]_{\mathtt{q}_n,\pi}e_n e_{n-1}e_n^2 -\pi e_{n-1}e_n^3=0,\\
e_{i}e_{i-1}e_{i}e_{i+1}&-e_{i}e_{i+1}e_{i}e_{i-1}+(-1)^{\epsilon_i}[2]e_{i}e_{i-1}e_{i+1}e_{i}-e_{i-1}e_{i}e_{i+1}e_{i}+e_{i+1}e_{i}e_{i-1}e_{i}=0,
\end{align*}
and the others are similar or already checked in \cite[Proposition~4.12]{CFLW}.
Indeed, the first one follows from
\begin{align*}
e_{n-1} \Upsilon^{-1}_n T_n &= \mathbf{i}^{\phi_{n,n-1}-d_n a_{n, n-1}} \Upsilon^{-1}_n T_n e_{n-1}=\Upsilon^{-1}_n T_n e_{n-1}, \\
e_{n} \Upsilon^{-1}_{n-1} T_{n-1} &= \mathbf{i}^{-d_{n-1} a_{n-1, n}} \Upsilon^{-1}_{n-1} T_{n-1} e_{n}=-\Upsilon^{-1}_{n-1} T_{n-1} e_{n}, \\
\widehat{\Psi}([3]_{\mathtt{q}_n,\pi})&=\widehat{\Psi}\left( \frac{\pi \mathtt{q}^{3}_n - \mathtt{q}^{-3}_n}{\pi\mathtt{q}_n-\mathtt{q}^{-1}_n}\right) = -[3]_{\mathtt{q}_n,\pi}.
\end{align*}
and the latter one from
\begin{align*}
e_{j}\Upsilon^{-1}_{j-1} T_{j-1} &= \mathbf{i}^{-d_{j-1} a_{j-1, j}} \Upsilon^{-1}_{j-1} T_{j-1} e_{j}=-\Upsilon^{-1}_{j-1} T_{j-1} e_{j}, \\
e_{j-1}\Upsilon^{-1}_{j} T_{j} &= \mathbf{i}^{\phi_{j,j-1}-d_{j} a_{j, j-1}} \Upsilon^{-1}_{j} T_{j} e_{j-1}=\Upsilon^{-1}_{j} T_{j} e_{j-1},
\end{align*}
for $j\leq n$ and 
\[
e_{j} \Upsilon^{-1}_{k} T_{k} =  \Upsilon^{-1}_{k} T_{k} e_j \text{ if }|j-k|>1, \quad \widehat{\Psi}([2])=-[2].
\]
\qed
\smallskip

For an idempotent $\varepsilon_\pm=\frac{1\pm\pi}{2}$, since $\pi \varepsilon_\pm =\pm \varepsilon_\pm$
we may identify a subalgebra $\varepsilon_\pm \widehat{\mathbf{U}}$ with a specialization $\widehat{\mathbf{U}}|_{\pi=\pm1}$.
In particular, the map $\widehat{\Psi}$, called a twistor \cite{CFLW}, induces an algebra isomorphism between $\widehat{\mathbf{U}}|_{\pi=1}$ and $\widehat{\mathbf{U}}|_{\pi=-1}$.

To compare these specializations with the generalized quantum groups, we impose  the following additional relations:
\[
J_i=1,\quad T^{2}_i=1\quad(i\in\mathring{I}),
\]
on $\widehat{\mathbf{U}}$. These relations are justified as they are always satisfied on the modules we shall consider (see \eqref{eq:T_i action} below and Remark~\ref{rem:wt of m}).
Since the twistor $\widehat{\Psi}$ preserves them, it descends to an automorphism on the resulting quotient, which we still denote by $\widehat{\Psi}$.

\begin{prop}
Under the above assumptions, the following hold:
\begin{itemize}
\item[(1)] The subalgebra of $\widehat{\mathbf{U}}|_{\pi=1}$ (resp. $\widehat{\mathbf{U}}|_{\pi=-1}$) generated by
$e_i$, $f_i$, and $k_\mu$ $(i\in\mathring{I},\,\mu\in P)$ is isomorphic to
$\mathcal{U}_B(\epsilon)$ (resp. $\mathcal{U}_{B^\bullet}(\epsilon)$), up to a rescaling of $e_n$.

\item[(2)] The image of the subalgebra
$\mathcal{U}_{B^\bullet}(\epsilon)\subset\widehat{\mathbf{U}}|_{\pi=-1}$
under $\widehat{\Psi}$ is isomorphic to
$\mathcal{U}_{B}(\epsilon)$ with the quantum parameter $\mathbf{v}^{-1}$.
\end{itemize}
\end{prop}
\pf It is straightforward to verify that, upon specializing $\pi=1$ (resp. $\pi=-1$), the defining relations for $e_i,f_i, k_\mu$ ($i\in\mathring{I},\mu\in P$) coincide with those of $\mathcal{U}_B (\epsilon)$ (resp. $\mathcal{U}_{B^\bullet} (\epsilon)$), up to a rescaling of $e_n$.
Since $\widehat{\Psi}(q^{1/2})=\mathbf{v}^{-1}$ and $\widehat{\Psi}(\pi)=-\pi$, the subalgebra $\mathcal{U}_{B^\bullet}(\epsilon)\subset\widehat{\mathbf{U}}|_{\pi=-1}$ is then isomorphically mapped by $\widehat{\Psi}$ onto $\mathcal{U}_{B}(\epsilon)$ but with the quantum parameter $\mathbf{v}^{-1}$.
\qed
\smallskip

We may define $\widehat{\bf U}[\sigma]$ as the semidirect product of $\widehat{\bf U}$ and $\Sigma$ as in Section \ref{subsec:Yamane QSA}, 
where $\sigma_i$ and ${\bm \sigma}$ commute with $T_j$, $\Upsilon_j$ ($j\in \ring{I}$) and $T_\La$.
Then $\widehat{\bf U}[\sigma]|_{\pi=1}=\varepsilon_+ \widehat{\bf U}[\sigma]$ (resp. $\widehat{\bf U}[\sigma]|_{\pi=-1}=\varepsilon_- \widehat{\bf U}[\sigma]$) contains $\mathcal{U}_{B}(\epsilon)[\sigma]$ (resp. $\mathcal{U}_{B^\bullet}(\epsilon)[\sigma]$),
and $\widehat{\Psi}$ naturally extends to an automorphism of $\widehat{\bf U}[\sigma]$ fixing $\Sigma$.

\subsection{Twisted classical limits of ${\U}_B(\e)$-modules}\label{subsec:twisted-cls-lim}

Let $V$ be a ${\U}_B(\e)$-module with weight space decomposition in the sense of \eqref{eq:weight space} such that
\begin{equation*}
{\rm wt}(V)\subset \bigcup_{k=1}^l \left(\ell_k\La + \sum_{a\in \I}\Z_{+}\de_a\right),
\end{equation*}
for some $\ell_1,\dots,\ell_k\in \Z_+$.
For $v\in V_\la$ with $\la=\ell\La +\sum_{a}m_a\de_a$, we define
\begin{gather}
T_i v = 
\begin{cases}
\ib^{2(m_{i+1}-m_i)}v & \text{if $1\le i<n$},\\
\ib^{-2m_n}v & \text{if $i=n$},
\end{cases}\quad 
T_\La v = \ib^{2(\La|\la)} v,\label{eq:T_i action}
\\
\Upsilon_i v = 
\begin{cases}
\ib^{-2m_i}v & \text{if $1\le i<n$},\\
\ib^{\sum_am_a -2m_n}v & \text{if $i=n$},
\end{cases}\quad
\sigma_i v= (-1)^{\e_im_i}v,\quad 
{\bm\sigma}v= \ib^{2(\La|\la)}v.
\end{gather}
Then it is straightforward to check that $V$ becomes a well-defined $\widehat{\bf U}[\sigma]|_{\pi=1}$-module.

Suppose that $\e=\bm{\e}=(\e_1,\dots,\e_n)=(\dots,0,1,0)$ and consider the embedding $\psi$ defined as the composition of the following maps:
\begin{equation}\label{eq: psi}
\xymatrixcolsep{2pc}\xymatrixrowsep{3pc}\xymatrix{
 U_{B^\bullet}(\e) \ \ar@{^{(}->}[r] &\ 
U_{B^\bullet}(\e)[\sigma] \ \ar@{->}[r]^{\substack{\quad \tau^{-1} \\ \cong}} &\
\U_{B^\bullet}(\e)[\sigma] \ \ar@{^{(}->}[r] &\ 
\widehat{\bf U}[\sigma]|_{\pi=-1} \ \ar@{->}[r]^{\substack{\widehat{\Psi} \\ \cong}} &\
\widehat{\bf U}[\sigma]|_{\pi=1}},
\end{equation}
where $\tau$ is the isomorphism in Proposition \ref{prop:tau iso bullet}. 
One should note that the pullback of a $\mathcal{U}_{B}(\epsilon)$-module $V$ by $\psi$ is a module over $\psi(U_{B^\bullet}(\epsilon))=U_{\psi^{-1}(q^{1/2})}(\mathfrak{osp}_{B^\bullet}(\epsilon))$,
and hence to obtain an honest $U_{B^\bullet}(\epsilon)$-action we have to first switch the quantum parameter $q^{1/2}$ of the given $\mathcal{U}_B (\epsilon)$-module $V$ into $\psi(q^{1/2})=\mathbf{i}^{-1}q^{1/2}={\bf v}^{-1}$.
We denote by $V^\psi=\{v^\psi\,|\,v\in V\}$ the resulting $U_{B^\bullet}(\epsilon)$-module.

\begin{lem}
For $\beta\in \sum_{i\in \ring{I}}\Z\alpha_i$ and $v\in V_\la$ with $\la=\ell\La +\sum_{a}m_a\de_a$, we have
$K_\beta v^\psi = (-1)^{2(\lambda|\beta)}q^{(\la|\beta)}v^\psi$.
\end{lem}
\pf By Proposition \ref{prop: twistor} and \eqref{eq:T_i action}, we have 
\begin{align*}
\psi(K_i) v 
&= T_ik_i\sigma_i\sigma_{i+1}v\\
&=(-1)^{m_i+m_{i+1}}(-q)^{(-1)^{\e_{i+1}}m_{i+1}-(-1)^{\e_{i}}m_{i}} v\\
&= q^{(\la|\alpha_i)}v
\end{align*}
for $i\neq n$ where $(\Lambda|\alpha_i)$=0, and 
\begin{align*}
\psi(K_n) v
&=T_n k_n v \\
&=(-1)^{m_n}(-q^{-1/2})^\ell (-q)^{-m_n}v\\
&=(-1)^\ell (q^{1/2})^{2(\la|\alpha_n)}v.
\end{align*}
\qed
\smallskip

Let $V$ be a highest weight module of $\mathcal{U}_{B}(\epsilon)$ with the highest weight $\mu\in  \ell\La + \sum_{a\in \I}\Z_+\de_a$.
Since any weight of V is also contained in the same set, by multiplying $E_n$ and $K_n$ by $(-1)^{\ell}$
we may assume that
\[
K_\beta v^\psi = q^{(\lambda|\beta)}v^\psi.
\]
Therefore, the weight space decomposition of $V^{\psi}$ as a $U_{q}(\mathfrak{osp}_{B^\bullet}(\epsilon))$-module 
coincides with the one of $V$ as a $\mathcal{U}_{B}(\epsilon)$-module.
It is also clear that $V^\psi$ is a highest weight module with the same highest weight $\mu$.
Now in the classical limit $\ov{V^{\psi}}$ of $V^\psi$  (see \cite[Section 2.3]{KLO24} for more details),
we have for a weight vector $v^\psi \in V^\psi$ of the weight $\lambda$ and the simple coroot $h_i$ of $\mathfrak{osp}_{B^\bullet}(\epsilon)$
\[
h_i v^\psi = \lim_{q^{1/2}\rightarrow 1} (-1)^{\epsilon_i}\frac{K_i - K^{-1}_i}{q^{\frac{d_i}{2}}-q^{-\frac{d_i}{2}}} v^\psi = (-1)^{\epsilon_i}\frac{2(\lambda|\alpha_i)}{d_i}v^\psi = \langle \lambda,h_i \rangle v^\psi.
\]
Summing up, we conclude with the following result for the $\mf{osp}_{B^\bullet}(\e)$-module $\ov{V^{\psi}}$ which we call a \emph{twisted classical limit} of a $\mathcal{U}_{B}(\epsilon)$-module $V$.

\begin{prop}\label{thm:classical limit of hw module}
Under the above hypothesis, $\ov{V^{\psi}}$ is a highest weight $\mf{osp}_{B^\bullet}(\e)$-module with the same highest weight and the same character as $V$.
\end{prop}
\smallskip

\section{$q$-analogues of oscillator representations of $\mf{osp}_{1|2n}$}\label{sec:osc rep}

\subsection{Oscillator $\mf{osp}_{1|2n}$-modules and finite-dimensional $\mf{so}_{2n+1}$-modules}\label{subsec:Howe duality}
Let
\begin{itemize}
\item[$\bullet$] $\mf{g}=\mf{osp}_{1|2n}=\mf{osp}_{B^\bullet}(0^n)$ or $\mf{so}_{2n+1}=\mf{osp}_{B}(1^n)$,

\item[$\bullet$] $\mathtt{P} = \bigoplus_{i=1}^n\Z\upvarpi_i$: the weight lattice of $\mf{g}$, where $\upvarpi_i$ is the $i$-th fundamental weight, 
 
\item[$\bullet$] $\ude_1=-\upvarpi_1$, $\ude_i=
  \upvarpi_{i-1}-\upvarpi_{i}$ for $1<i< n$, and $\ude_n=\upvarpi_{n-1}-2\upvarpi_n$,

\item[$\bullet$] $\upalpha_i$: the simple root for $\mf{g}$, where $\upalpha_i=\ude_{i+1}-\ude_{i}$ for $1\le i< n$, and $\upalpha_n=-\ude_n$.
\end{itemize}
Note that
$\upvarpi_i=
2\upvarpi_n + \ude_{n}+\dots+\ude_{i+1}$ $(1\le i<n)$,
and hence $\ude_1+\dots+\ude_n=-2\upvarpi_n$. 
For $\ell\ge 0$ and a partition (or Young diagram) $\la=(\la_i)_{i\geq 1}$ with length $\ell(\la)\le {\lfloor\ell/2 \rfloor}$, let
\begin{equation*}\label{eq:Lambda highest weight}
\begin{split}
\upvarpi_{(\la,\ell),\mf{g}} 
&=
\begin{cases}
-\ell \upvarpi_n +\sum_{i=1}^n\la_i\ude_{n-i+1} & \text{if $\mf{g}=\mf{osp}_{1|2n}$ and $\ell(\la)\le n$},\\
\ell \upvarpi_n +\sum_{i=1}^n\la'_i\ude_{n-i+1} & \text{if $\mf{g}=\mf{so}_{2n+1}$ and $\ell(\la')\le n$},\\
0 & \text{otherwise},
\end{cases}
\end{split}
\end{equation*}
where $\la'=(\la'_i)_{i\ge 1}$ denotes the conjugate of $\la$.
We also let
\begin{itemize}
\item[$\bullet$] $\mc{P}(\ell,\mf{g})=\{\,(\la,\ell)\,|\, \upvarpi_{(\la,\ell),\mf{g}} \neq 0 \,\}$ for $\ell\geq 1$ and $\mathcal{P}(0,\mathfrak{g})=\{(\emptyset,0)\}$,

\item[$\bullet$] $V((\la,\ell),{\mf{g}})$: the irreducible highest weight $\mf{g}$-module with highest weight $\upvarpi_{(\la,\ell),\mf{g}}$ for $\la\in \mc{P}(\ell,\mf{g})$,

\item[$\bullet$] $W_{\mf{g}}=V(((0),1),\mf{g})$.
\end{itemize}

We remark that the $\mathfrak{osp}_{1|2n}$-module $V((\la,\ell),\mf{osp}_{1|2n})$, called an {\em oscillator representation}, is infinite-dimensional while $V\left((\la,\ell),\mf{so}_{2n+1}\right)$ is finite-dimensional. 
In particular, $W_{\mf{osp}_{1|2n}}$ is isomorphic to the classical limit $\ov{\W^{\psi}_{(0^n)}}$ by Proposition~\ref{thm:classical limit of hw module}, and $W_{\mf{so}_{2n+1}}$ is the spin representation.

In this subsection, we prove the semisimplicity of $W_{\mf{g}}^{\ot \ell}$ and describe its decomposition into irreducible modules. 

We first assume $\mf{g}=\mf{so}_{2n+1}$, where the semisimplicity is clear. Let $W=W_{\mathfrak{so}_{2n+1}}$ and consider its crystal ${B}(W)$, which can be identified with $\Z_+^n(1^n)$ whose highest weight element is ${\bf 0}$ of the weight $\upvarpi_n$ (\emph{cf.}~\cite{Kas95,KashNaka}).
The crystal operators $\tf_i$ ($i\in \ring{I}$) is then given by $\tf_n {{\bf m}}={{\bf m}+{\bf e}_n}$ and $\tf_i {{\bf m}}={{\bf m}+{\bf e}_i-{\bf e}_{i+1}}$ for $1\le i <n$ and ${\bf m}\in \Z_+^n(1^n)=\{0,1\}^n$. 

The irreducible components of ${B}(W^{\ot \ell})\cong{B}(W)^{\ot \ell}$ are parametrized by highest weight elements, namely ${{\bf m}_1}\ot \dots \ot {{\bf m}_\ell}$ satisfying
\begin{equation}\label{eq:h.w. crystal element}
 \te_i({{\bf m}_1}\ot \dots \ot {{\bf m}_\ell})=0\quad (1\le i\le n),
\end{equation}
where the tensor product rule of crystals is given with respect to the comultiplication in \eqref{eq:comult-e}.

\begin{lem}\label{lem:h.w. vectors in spin tensors}
The highest weight elements in $B(W)^{\ot\ell}$ are in one-to-one correspondence with the sequences of Young diagrams
\[
\emptyset=\la^{(0)}\subset\la^{(1)}\subset\cdots\subset\la^{(\ell)}\subset(\lfloor\ell/2\rfloor^n)
\]
satisfying for all $1\leq i\leq \ell$
\begin{enumerate}
\item $\la^{(i)}_1\le\lfloor i/2\rfloor$,
\item $\la^{(i)}/\la^{(i-1)}$ is a vertical strip.
\end{enumerate}
Moreover, the weight of the corresponding element is given by $\ell \upvarpi_n +\sum_{i=1}^n\la^{(\ell)}_i\ude_{n-i+1}=\upvarpi_{(\la,\ell),\mf{so}_{2n+1}}$, where $\la=(\la^{(\ell)})'$.
\end{lem}
\pf
Suppose we are given a highest weight element ${\bf m}_1\ot\cdots\ot{\bf m}_{\ell}$ in $B(W)^{\ot\ell}$.
By the tensor product rule of crystals, this is equivalent to
\begin{itemize}
\item[(a)] ${\bf m}_1\ot\cdots\ot{\bf m}_{\ell-1}$ is a highest weight vector in $B(W)^{\ot\ell-1}$, and
\item[(b)] $\varepsilon_i({\bf m}_\ell)\le\langle h_i,\mathrm{wt}({\bf m}_1\ot\cdots\ot{\bf m}_{\ell-1})\rangle$ for $1\le i\le n$,
\end{itemize}
where $h_i$ is the coroot of $\alpha_i$, and $\varepsilon_i({\bf m}_\ell)=\max\{\,k\,|\,\te^k_i{\bf m}_\ell\neq 0\,\}$. 
We shall assign to such an element a sequence of Young diagrams inductively in $\ell$.
Note that when $\ell=1$, the only highest weight element in $B(W)$ is ${\bold 0}$ and it corresponds to $\la^{(1)}=\emptyset$. 

We may assume that ${\bf m}_1\ot\cdots\ot{\bf m}_{\ell-1}$ has been associated to a sequence of 
Young diagrams $\emptyset=\la^{(0)}\subset\la^{(1)}\subset\cdots\subset\la^{(\ell-1)}\subset(\lfloor(\ell-1)/2\rfloor^n)$ and 
$\mathrm{wt}({\bf m}_1\ot\cdots\ot{\bf m}_{\ell-1})=(\ell-1)\varpi_n+\sum_{i=1}^n{\la_i^{(\ell-1)}}\delta_{n-i+1}$.
In particular, we have
\[
\langle h_i,\mathrm{wt}({\bf m}_1\ot\cdots\ot{\bf m}_{\ell-1})\rangle=
\left\{\begin{array}{ll}
\ell-1-2\la^{(\ell-1)}_1&(i=n),\\
\la^{(\ell-1)}_{n-i}-\la^{(\ell-1)}_{n-i+1}\;&(i<n).
\end{array}\right.
\]
Let us define a Young diagram $\lambda^{(\ell)}$ by adding a box to the $i$-th row of $\lambda^{(\ell)}$ whenever $(\mathbf{m}_{\ell})_{n-i+1}$ is non-zero.
Then the condition (b) above is equivalent to that $\lambda^{(\ell)}\leq \lfloor\ell/2\rfloor$ and $\lambda^{(\ell)}/\lambda^{(\ell-1)}$ is a vertical strip.
It is also clear that the weight of the highest weight element is given as stated, and hence we obtain a sequence of Young diagrams with the desired properties.

Conversely, given a sequence $(\lambda^{(i)})_{i=1,\dots,\ell}$ we define $(\mathbf{m}_{i})_k = \lambda^{(i)}_{n-k+1} - \lambda^{(i-1)}_{n-k+1}$,
and the conditions (1) and (2) ensure that $\mathbf{m}_{1}\otimes \cdots \otimes \mathbf{m}_{\ell}$ is a highest weight element in $B(W)^{\otimes\ell}$.
This completes the proof of the Lemma.
\qed
\smallskip

Given a partition $\la$  with $\ell(\la)\le {\lfloor\ell/2 \rfloor}$, one can easily see from Lemma \ref{lem:h.w. vectors in spin tensors} that if $n \ge \ell(\lambda')$, then the multiplicity of $V((\la,\ell),{\mf{so}_{2n+1}})$ in $W^{\ot \ell}$ is non-zero and independent of $n$, which we denote by $d^{(\la,\ell)}$.
Let $SST_\ell (\mu)$ be the set of semistandard tableaux of shape $\mu$ with entries in $\{1,\dots,\ell\}$.
Then we have the following description of $d^{(\la,\ell)}$ in terms of semistandard tableaux.

\begin{cor}\label{rem: King tableaux}
The multiplicity $d^{(\la,\ell)}$ is equal to the number of  $T\in SST_{\ell}(\la)$ such that $t_k\ge 2k$ for all $k\ge 1$, where $t_k$ denotes the entry of $T$ located in the first column (from the left) and the $k$-th row (from the top) of $\la$.
\end{cor}
\pf Suppose that ${{\bf m}_1}\ot \cdots \ot {{\bf m}_\ell}\in B(W)^{\ot \ell}$ is a highest weight vector \eqref{eq:h.w. crystal element} associated with a sequence of Young diagrams $(\lambda^{(i)})_{i=1,\dots,\ell}$ as in Lemma \ref{lem:h.w. vectors in spin tensors}. 
Let $T$ be a tableau of shape $\la$ such that its subtableau $T^{(i)}$ corresponding to the skew shape $(\la^{(i)})'/(\la^{(i-1)})'$ is filled with $i$ for $1\le i\le \ell$. Then we have $T\in SST_{\ell}(\la)$ by the condition (2), and $t_k\ge 2k$ for all $k\ge 1$ by (1). 
Since the assignment ${{\bf m}_1}\ot \dots \ot {{\bf m}_\ell} \mapsto T$ is obviously reversible, the claim follows.
\qed

\begin{rem}\label{rem: decomposition number}
{\rm
One may notice that $d^{(\la,\ell)}$ is also equal to the number of  $T\in SST_{\ell-1}(\la)$ such that $t_k\ge 2k-1$ for $k\ge 1$ (reducing all the entries in $T$ by 1). 
If $\ell$ is odd, it is the number of symplectic tableaux introduced by King \cite{Kng}, and hence the dimension of the finite-dimensional irreducible $Sp_{\ell-1}$-module corresponding to $\la$. On the other hand, when $\ell$ is even, $d^{(\la,\ell)}$ is equal to the number of (odd) symplectic tableaux of shape $\la$ introduced by Proctor \cite{Pro}, which gives the dimension of an indecomposable module of the odd symplectic group $Sp_{\ell-1}$.
} 
\end{rem}

\begin{prop}\label{prop: decomp of tensor power of W}
Let $\mathfrak{g}=\mathfrak{so}_{2n+1}$ or $\mathfrak{osp}_{1|2n}$. For $\ell\ge 1$, $W_{\mf{g}}^{\ot \ell}$ is semisimple and
\begin{equation*}
W_{\mf{g}}^{\ot \ell} \cong 
\bigoplus_{(\la,\ell)\in \mc{P}(\ell,\mf{g})} V((\la,\ell),{\mf{g}})^{\oplus d^{(\la,\ell)}}.
\end{equation*}
\end{prop}
\pf
The case when $\mf{g}=\mf{so}_{2n+1}$ is done in Lemma \ref{lem:h.w. vectors in spin tensors}. 
To prove the decomposition for $\mf{g}=\mf{osp}_{1|2n}$, we apply the theory of super duality \cite{CLW}.

Let us consider a triple of Lie (super)algebras of infinite rank $(\td{\mf{g}}^{\mf{b}^\bullet},{\mf{g}}^{\mf{b}^\bullet},\ov{\mf{g}}^{\mf{b}^\bullet})$ \cite[Section 2.3]{CLW}, where ${\mf{g}}^{\mf{b}^\bullet}$ and $\ov{\mf{g}}^{\mf{b}^\bullet}$ are the Lie (super)algebra obtained by taking $n\rightarrow \infty$ in $\mf{osp}_{1|2n}$ and $\mf{so}_{2n+1}$, respectively. 
First, there are equivalences of categories $T: \td{\mc{O}} \longrightarrow {\mc{O}}$ and $\ov{T}: \td{\mc{O}} \longrightarrow \ov{\mc{O}}$ between certain parabolic BGG categories $\td{\mc{O}}$, ${\mc{O}}$ and $\ov{\mc{O}}$ for $\td{\mf{g}}^{\mf{b}^\bullet},{\mf{g}}^{\mf{b}^\bullet}$ and $\ov{\mf{g}}^{\mf{b}^\bullet}$, respectively \cite[Theorem 5.4]{CLW}. 
If we let $\td{W}$, $W$, and $\ov{W}$ denote the irreducible highest weight modules of $\td{\mf{g}}^{\mf{b}^\bullet},{\mf{g}}^{\mf{b}^\bullet}$ and $\ov{\mf{g}}^{\mf{b}^\bullet}$ of highest weight $-\La_0$ (see~\cite[Section 2.4]{CLW}), respectively,
then $T(\td{W})=W$ and $\ov{T}(\td{W})=\ov{W}$ by \cite[Theorem 4.6]{CLW}.
\begin{equation*}
\begin{tikzcd}[row sep=.8cm, column sep=.5cm]
& \widetilde{\mathcal{O}}\arrow[ld,"T"',"\simeq"] \arrow[rd,"\overline{T}","\simeq"']&\\
\mathcal{O} \arrow[d,"\mathfrak{tr}^{\infty}_n"']&& \overline{\mathcal{O}} \arrow[d,"\mathfrak{tr}^{\infty}_n"]\\
\mathfrak{osp}_{1|2n}\text{-Mod} && \mathfrak{so}_{2n+1}\text{-Mod}
\end{tikzcd}
\end{equation*}

Next, we may further truncate $\overline{V}\in\overline{\mathcal{O}}$ to a module $\mathfrak{tr}^{\infty}_n (\overline{V})$ over the rank $n$ Lie subalgebra $\ov{\mf{g}}^{\mf{b}^\bullet}_n = \mathfrak{so}_{2n+1}$ (as in Section~\ref{sec:truncation}),
and similarly $V\in\mathcal{O}$ to a $\mf{g}^{\mf{b}^\bullet}_n=\mathfrak{osp}_{1|2n}$-module $\mathfrak{tr}^{\infty}_n (V)$.
While it is not an equivalence, it maps an irreducible object to an irreducible object or zero, and is faithful on the irreducibles that are not mapped to zero \cite[Theorem 4.6]{CLW}.
In particular, $\mathfrak{tr}^{\infty}_n (\overline{W})$ and $\mathfrak{tr}^{\infty}_n (W)$ are isomorphic to $W_{\mathfrak{so}_{2n+1}}$ and $W_{\mathfrak{osp}_{1|2n}}$, respectively \cite[Sections 3.3 and 3.4]{CLW}.

Now we prove the decomposition. By Lemma~\ref{lem: h.w. module in tensor power}, $\ov{W}^{\ot\ell}$ is semisimple and its irreducible components are parametrized by $(\lambda,\ell)$ satisfying $\ell(\lambda)\leq \lfloor\ell/2\rfloor$ and have the multiplicity $d^{(\lambda,\ell)}$.
Since $T$ and $\overline{T}$ are equivalences, $\td{W}^{\ot\ell}$ (resp. $W^{\ot\ell}$) have formally the same decomposition into irreducibles $\widetilde{V}=\overline{T}^{-1}(\overline{V})$ (resp. $V=T(\widetilde{V})$) for each irreducible component $\overline{V}$ in $\ov{W}^{\ot\ell}$, sharing the same multiplicities.
Furthermore, by the property of $\mathfrak{tr}^{\infty}_n$ above,
the multiplicity of $\mathfrak{tr}^{\infty}_n(\overline{V})$ in $\mathfrak{tr}^{\infty}_n(\overline{W}^{\otimes\ell})$ is equal to that of $\overline{V}$ in $\overline{W}^{\otimes\ell}$ whenever $\mathfrak{tr}^{\infty}_n(\overline{V})$ is non-zero (hence irreducible).
The same holds for the irreducible components of $W^{\otimes\ell}$ and their truncations.
Therefore, it remains to recall from \cite[Lemma~3.2, Theorem~4.6]{CLW} that for an irreducible component $\ov{V}$ (and $V=T(\overline{T}^{-1}(\overline{V}))$) corresponds to $(\la,\ell)$, 
we have $\mf{tr}^{\infty}_n(\ov{V})\neq 0$ (resp. $\mf{tr}^{\infty}_n({V})\neq 0$) if and only if $\upvarpi_{(\la,\ell),\mf{so}_{2n+1}}\neq 0$ (resp. $\upvarpi_{(\la,\ell),\mf{osp}_{1|2n}}\neq 0$), 
in which case the highest weight of $\mf{tr}^{\infty}_n(\ov{V})\neq 0$ (resp. $\mf{tr}^{\infty}_n({V})\neq 0$) is exactly $\upvarpi_{(\la,\ell),\mf{so}_{2n+1}}$ (resp. $\upvarpi_{(\la,\ell),\mf{osp}_{1|2n}}$).

\qed\smallskip

\begin{rem}{\rm
Remark \ref{rem: King tableaux} and Proposition \ref{prop: decomp of tensor power of W} suggest the existence of an action of $(\mathfrak{g}, {Sp}_{\ell-1})$ on $W_{\mathfrak{g}}^{\otimes \ell}$ or its $q$-analogue, which would in particular be semisimple when $\ell$ is odd. 
However, we were unable to find such a result in the literature.
} 
\end{rem}

\begin{prop}\label{prop:tensor product decomp}
Let $(\mu,\ell)\in \mc{P}(\ell,\mf{g})$ and $(\nu,\ell')\in \mc{P}(\ell',\mf{g})$ be given. 
\begin{itemize}
\item[(1)] $V((\mu,\ell),{\mf{g}})\ot V((\nu,\ell'),{\mf{g}})$ is isomorphic to a direct sum of $V((\la,\ell+\ell'),{\mf{g}})$ for $(\la,\ell+\ell')\in \mc{P}(\ell+\ell',\mf{g})$.

\item[(2)] Let $c^{(\la,\ell+\ell')}_{(\mu,\ell),(\nu,\ell')}(\mf{g})$ denote its multiplicity. 
If $(\la,\ell+\ell')\in \mc{P}(\ell+\ell',\mf{osp}_{1|2n})\cap \mc{P}(\ell+\ell',\mf{so}_{2n+1})$, then we have  $c^{(\la,\ell+\ell')}_{(\mu,\ell),(\nu,\ell')}(\mathfrak{osp}_{1|2n})=c^{(\la,\ell+\ell')}_{(\mu,\ell),(\nu,\ell')}(\mf{so}_{2n+1})$.

\end{itemize}
\end{prop}
\pf It also follows from the arguments in the proof of Proposition \ref{prop: decomp of tensor power of W}. 
\qed
\smallskip

\begin{ex}{\rm 
The following decomposition is well-known:
\begin{equation*}
W_{\mf{so}_{2n+1}}\ot W_{\mf{so}_{2n+1}} \cong
\bigoplus_{i=0}^n V(((i),2),\mf{so}_{2n+1}), 
\end{equation*}
where 
\[
\upvarpi_{((i),2),\mf{so}_{2n+1}}=
\begin{cases}
\upvarpi_{n-i} & \text{for $1\le i<n$},\\
2\upvarpi_{n} & \text{if $i=0$},\\
0 & \text{if $i=n$}.
\end{cases}
\]
By Proposition \ref{prop:tensor product decomp}, we have
\begin{equation}
W_{\mf{osp}_{1|2n}}\ot W_{\mf{osp}_{1|2n}} \cong
\bigoplus_{i\ge 0} V(((i),2),\mf{osp}_{1|2n}). 
\end{equation}
}
\end{ex}

\subsection{Category ${\mc O}^{\,\mf b}_{{\rm osc},\e}$}\label{subsec:osc category type b}
Let us introduce a $q$-deformed and super analogue of the semisimple monoidal category of $\mf{g}$-modules generated by $V((\la,\ell))$ for $(\la,\ell)\in \mc{P}(\ell,\mf{g})$ ($\ell\ge 1$).
We keep the same notations as in \cite[Section~5.1]{KLO24} except that the role of ${\bm \e}$ and ${\bm\e}'$ are exchanged here.

From now on, we assume that $\e$ is one of ${\bm\e},\un{\bm\e},\ov{\bm\e}$ in Section \ref{sec:truncation}.
Let 
\begin{itemize}
\item[$\bullet$] $\mc{P}(\ell)=\{\,(\la,\ell)\,|\,\text{$\la$ is a partition with $\ell(\la)\le \lfloor \ell/2 \rfloor$}\,\}$, 

\item[$\bullet$] for $(\la,\ell)\in \mc{P}(\ell)$,
\begin{equation}\label{eq: highest weight}
\La_{(\la,\ell),\e}= 
\begin{cases}
\ell \La + \sum_{i\in \I_\e} m_i\de_{i} & \text{if $\la\in \cP_\e$},\\
 0 & \text{otherwise},
\end{cases}
\end{equation}

\item[$\bullet$] $\mc{P}(\ell)_\e=\{\,(\la,\ell)\in \mc{P}(\ell)\,|\,\La_{(\la,\ell),\e}\neq 0\,\}=\mc{P}(\ell)\cap\cP_\e$,

\item[$\bullet$] $\mc{V}_{\e}^{(\la,\ell)}$: the irreducible $\ring{\U}_\e$-module with highest weight $\La_{(\la,\ell),\e}$ for $(\la,\ell)\in \mc{P}(\ell)_\e$.
 
\end{itemize}

\begin{rem}\label{rem: decomposition for super}
{\rm 
For example, let $({\bm\e},\un{\bm\e},\ov{\bm\e})$ be the one for $X=B$ in Section~\ref{sec:truncation}.
Note that 
\begin{equation}\label{eq: parametrizations in over and under}
\begin{split}
\mc{P}(\ell)_{\ov{\bm \e}}=\mc{P}(\ell,\mf{so}_{2m+1}),\quad \La_{(\la,\ell),\ov{\bm \e}}=\upvarpi_{(\la,\ell),\mf{so}_{2m+1}},\\
\mc{P}(\ell)_{\un{\bm \e}}=\mc{P}(\ell,\mf{osp}_{1|2m}),\quad\La_{(\la,\ell),\un{\bm \e}}=\upvarpi_{(\la,\ell),\mf{osp}_{1|2m}}  
\end{split}
\end{equation}
by identifying $\mathsf{P}$ with the quotient of $P$ by $\Lambda=-\frac{1}{2}\sum (-1)^{\epsilon_i}\delta_i $.

The case of $\epsilon=\boldsymbol{\epsilon}$ amounts to the middle category $\widetilde{\mathcal{O}}$ in the proof of Proposition \ref{prop: decomp of tensor power of W}.
Namely, we may also consider the rank $n$ Lie subalgebra $\widetilde{\mathfrak{g}}^{\mathfrak{b}^\bullet}_n\cong \mathfrak{osp}_{B^\bullet}(\boldsymbol{\epsilon})$ and truncation from $\widetilde{\mathcal{O}}$ to a category of $\mathfrak{osp}_{B^\bullet}(\boldsymbol{\epsilon})$-modules.
Again from \cite[Lemma~3.2, Theorem~4.6]{CLW}, if $\td{V}$ is an irreducible submodule of $\td{W}^{\ot \ell}$ associated to $(\la,\ell)\in\mathcal{P}(\ell)$,
then $\mf{tr}^{\infty}_n(\td{V})\neq 0$ if and only if $\La_{(\la,\ell),{\bm \e}}\neq 0$, in which case its highest weight is $\La_{(\la,\ell),{\bm\e}}$.
Hence, if we let $W_{\mf{osp}_{B^\bullet}({\bm \e})}:=\mf{tr}^{\infty}_n(\td{W})(\neq 0)$ and $V((\lambda,\ell),\mathfrak{osp}_{B^\bullet}(\boldsymbol{\epsilon}))$ be the irreducible highest weight $\mathfrak{osp}_{B^\bullet}(\boldsymbol{\epsilon})$-module with the highest weight $\Lambda_{(\lambda,\ell),\boldsymbol{\epsilon}}$, 
then we have the same decomposition as in Proposition~\ref{prop: decomp of tensor power of W} for $\mathfrak{g}=\mathfrak{osp}_{B^\bullet}(\boldsymbol{\epsilon})$.
}
\end{rem}

We define $\mc{O}^{\mf b}_{\text{osc},\e}$ to be the category of $\ring{\U}_{\e}$-modules 
$V$ such that $V=\bigoplus_{\ell\ge 0}V_{\ell}$, where $V_\ell$ is a direct sum of $\mc{V}_{\e}^{(\la,\ell)}$'s for $(\la,\ell)\in \mc{P}(\ell)_{\e}$ with finite multiplicity for each $(\la,\ell)$, and $V_\ell=0$ for all sufficiently large $\ell$.
For example, let $\W_\e$ denote $\W_\e(x)$ ($x\in \Bbbk^\times$) regarded as a $\ring{\U}_\e$-module.
Then we have $\W_\e\in \mc{O}^{\mf b}_{\text{osc},\e}$ by Remark~\ref{rem:wt of m}.

The category $\mc{O}^{\mf b}_{\text{osc},\e}$ is obviously semisimple, and is closed under the tensor product by the $q$-analogue of the decomposition of $W^{\otimes \ell}_{\mathfrak{g}}$ , see Proposition~\ref{prop:decomp of W for e} below.
Moreover, we will show that truncation functors $\mathfrak{tr}_{\epsilon}$ ($\epsilon=\underline{\boldsymbol{\epsilon}},\overline{\boldsymbol{\epsilon}}$) send simples to simples or zero (Theorem~\ref{thm:trunc sends simple to simple}), which implies that they restrict to categories $\mathcal{O}^{\mathfrak{b}}_{\mathrm{osc},\epsilon}$.
While these results can be argued as in \cite[Proposition~6.15 -- Corollary~6.18]{KL22}, here we provide more direct and simpler proofs.

\begin{lem}\label{lem: h.w. module in tensor power}
For $\ell\ge 1$, 
\begin{itemize}
\item[(1)] $\W_{\e}^{\otimes\ell}$ is semisimple as a $\mathring{\U}_{\e}$-module,

\item[(2)] any irreducible $\ring{\U}_\e$-submodule in $\W_\e^{\ot \ell}$ is isomorphic to $\mc{V}^{(\la,\ell)}_\e$ for some $(\la,\ell)\in \mc{P}(\ell)_\e$.
\end{itemize}
\end{lem}

\pf (1) We assume $\e={\bm \e}$ since the proof for the other cases are similar. Let $\eta$ be an anti-involution on $\mathring{\mathcal{U}}_\epsilon$ given by
\[
\eta(e_i)=\begin{cases} (-q^2)^{\epsilon_i-\epsilon_{i+1}}q_i f_i k^{-1}_i& \\ \mathbf{v}f_n k^{-1}_n& \end{cases},\quad
\eta(f_i)=\begin{cases} (-q^2)^{\epsilon_{i+1}-\epsilon_{i}}q^{-1}_i k_i e_i& \text{if }i<n\\ \mathbf{v}^{-1} k_n e_n& \text{if }i=n\end{cases},\quad
\eta(k_\mu)=k_\mu.
\]
We also define a symmetric bilinear form on $\mathcal{W}_{\epsilon}$ by
\[
(\ket{\mathbf{m}},\ket{\mathbf{m}^\prime})=\delta_{\mathbf{m},\mathbf{m}^\prime}q^{-\sum\frac{m_i(m_i-1)}{2}}\prod[m_i]!,
\]
which can be easily seen to satisfy $(uv,w)=(v,\eta(u)w)$ for any $u\in\mathring{\mathcal{U}}_\epsilon$ and $v,w\in\mathcal{W}_{\epsilon}$.
Moreover, since $\eta$ satisfies $\Delta \circ \eta = (\eta\otimes\eta)\circ \Delta$, 
the same property holds for a similar form on tensor powers $\mathcal{W}^{\otimes\ell}_\epsilon$ defined by
\[
(v_1 \otimes \cdots \otimes v_\ell, w_1 \otimes \cdots \otimes w_\ell)=\prod (v_i, w_i).
\]

Let $A_\infty\subset\mathbb{Q}(q)$ denote the subring consisting of those regular at $q=\infty$,
and $\mathcal{L}^{\otimes \ell}_\infty$ the $A_\infty$-span of $\ket{\mathbf{m}_1}\otimes \cdots \otimes \ket{\mathbf{m}_\ell}$ for all $\mathbf{m}_i\in\mathbb{Z}^{n}_{+}(\epsilon)$
which is an $A_\infty$-lattice inside $\mathcal{W}^{\otimes\ell}_{\epsilon}$ compatible with the weight space decomposition.
Then we have $(\mathcal{L}^{\otimes \ell}_\infty,\mathcal{L}^{\otimes \ell}_\infty)\subset A_\infty$ 
and the image of the set $\left\{ \ket{\mathbf{m}_1}\otimes \cdots \otimes \ket {\mathbf{m}_\ell}\,|\,\mathbf{m}_i\in\mathbb{Z}^{n}_{+}(\epsilon)\right\}$ in $\mathcal{L}_\infty/q^{-1}\mathcal{L}_\infty$ forms an orthonormal basis with respect to the induced bilinear form which is hence positive definite.
This allows us to repeat the argument in the proof of \cite[Theorem~2.12]{BKK} to deduce the semisimplicity.

(2) By (1), $\mathcal{W}^{\otimes \ell}_{\epsilon}$ is a direct sum of irreducible highest weight modules.
Let us denote by $\overline{\mathcal{W}^{\otimes \ell}_{\epsilon}}$  the (twisted) classical limit (see Section~\ref{subsec:twisted-cls-lim}) with respect to this decomposition.
Recall that the classical limit of $\mathcal{W}_\epsilon$ is isomorphic to $W_{\mathfrak{g}}$ where $\mf{g}=\mf{osp}_{B^\bullet}({\bm \e}),\, \mf{osp}_{1|2m+2}$ ($X=B$) or $\mf{osp}_{1|2m}$ ($X=D$), and $\mf{so}_{2m+1}$ for $\e={\bm \e}, \un{\bm\e}$, and $\ov{\bm\e}$, respectively.
Hence, the character of $\overline{\mathcal{W}^{\otimes \ell}_{\epsilon}}$ is equal to $\overline{\mathcal{W}_{\epsilon}}^{\otimes \ell}=W^{\otimes \ell}_{\mathfrak{g}}$ whose decomposition is given in Proposition~\ref{prop: decomp of tensor power of W} ($\epsilon=\underline{\boldsymbol{\epsilon}},\overline{\boldsymbol{\epsilon}}$) and Remark~\ref{rem: decomposition for super} ($\epsilon=\boldsymbol{\epsilon}$).
Now that the classical limit of an irreducible submodule in $\mathcal{W}^{\otimes \ell}_{\epsilon}$ is a highest weight module contained in $\overline{\mathcal{W}^{\otimes \ell}_{\epsilon}}$,
its highest weight must be of the form $\Lambda_{(\lambda,\ell),\epsilon}$.
\qed

\begin{prop}\label{prop:decomp of W for e}
For $\e={\bm\e}, \un{\bm\e}, \ov{\bm\e}$, we have 
\begin{equation*} 
\W_{\e}^{\ot \ell} \cong \bigoplus_{\la\in \mc{P}(\ell)_{\e}} \left(\mc{V}_{\e}^{(\la,\ell)}\right)^{\oplus d^{(\la,\ell)}}.
\end{equation*}
\end{prop}

\pf 
By the lemma above, the classical limit $\overline{\mathcal{W}^{\otimes \ell}_{\epsilon}}$ is a direct sum of $\overline{\mc{V}_{\e}^{(\la,\ell)}}$, the classical limit of the irreducible highest weight module $\mc{V}_{\e}^{(\la,\ell)}$ and hence is also a highest weight module of highest weight $\Lambda_{(\lambda,\ell),\epsilon}$.
Since $\overline{\mathcal{W}^{\otimes \ell}_{\epsilon}}$ has the same character as $W^{\otimes \ell}_{\mathfrak{g}}$,
it is enough to prove that each $\overline{\mc{V}_{\e}^{(\la,\ell)}}$ is irreducible, which follows once $\overline{\mathcal{W}^{\otimes \ell}_{\epsilon}}$ is semisimple.
To this end, we shall show that the irreducible modules $V((\lambda,\ell),\mathfrak{g})$ have distinct central characters by invoking the linkage principle.

Assume $\epsilon=\boldsymbol{\epsilon}$ and $X=D$, for the others being parallel and easier.
Recall that the weight lattice $\mathtt{P}$ of $\mathfrak{osp}_{B^\bullet}(\boldsymbol{\epsilon})$ is given by (see Remark~\ref{rem: decomposition for super})
\[
\mathtt{P}=P/\mathbb{Z}\cdot(\Lambda + \frac{1}{2}\sum^{2m}_{i=1} (-1)^{\epsilon_i}\delta_i) \,\subset \, \bigoplus_{i\in\mathbb{I}} \mathbb{Z} \cdot \frac{1}{2}\delta_i = L_1 \oplus L_2, \quad L_k = \bigoplus_{\epsilon_i =k-1}\mathbb{Z} \cdot \frac{1}{2}\delta_i.
\]
The Weyl vector associated to a sequence $\epsilon$ is defined to be the half of the signed sum of positive roots, which equals (the image of) $\Lambda$ for $\epsilon=\boldsymbol{\epsilon}$.
The Weyl group $W$ of $\mathfrak{osp}_{B^\bullet}(\boldsymbol{\epsilon})$ is the product of two copies of the Weyl group $S_m \ltimes (\mathbb{Z}/2\mathbb{Z})^m$ of type $B_{m}$,
and the $k$-th copy acts on $\mathtt{P}\cap L_k$ in the usual type $B$ manner.

It is well-known \cite[Theorem~2.30]{CW} that for $\mu,\nu \in \mathtt{P}$, their associated central characters are the same if and only if $\mu$ and $\nu$ are linked, namely
\[
\nu +\Lambda = w\left(\mu +\Lambda - \sum_a c_a \beta_a \right)
\]
for some $w\in W$, $c_a \in \mathbb{C}$ and mutually orthogonal isotropic odd roots $\beta_a$ such that $(\mu+\Lambda | \beta_a)=0$.
For $(\lambda,\ell)\in \mathcal{P}(\ell)_{\boldsymbol{\epsilon}}$, the weight $\Lambda_{(\lambda,\ell),\boldsymbol{\epsilon}}+\Lambda \in \mathtt{P}$ is of the form $\sum_{i\in\mathbb{I}} b_i \delta_i$, where
\[
b_1 > b_3 > \cdots > b_{2m-1} > \frac{l+1}{2} > 0 > b_{2} > b_4 > \cdots > b_{2m} \geq -\frac{l+1}{2}.
\]
In particular, there is no isotropic odd root (which is of the form $\delta_i \pm \delta _{i+2k+1}$) orthogonal to it.
It is also clear that $\Lambda_{(\lambda,\ell),\boldsymbol{\epsilon}}+\Lambda$ is the only such weight in its Weyl group orbit, and hence $\Lambda_{(\lambda,\ell),\boldsymbol{\epsilon}}$ and $\Lambda_{(\lambda',\ell),\boldsymbol{\epsilon}}$ are not linked unless $\lambda=\lambda'$.
This concludes the proof.

\qed

\begin{thm}\label{thm:trunc sends simple to simple}
For $(\la,\ell)\in \mc{P}(\ell)_{\bm \e}$,
\begin{equation}\label{eq: truncation of simple}
\mf{tr}_{\e}\left(\mc{V}^{(\la,\ell)}_{\bm\e}\right)\cong
\begin{cases}
\mc{V}^{(\la,\ell)}_{\e}  & \text{if $(\la,\ell)\in \mc{P}(\ell)_{\e}$},\\
0 & \text{otherwise},
\end{cases}\quad (\e=\un{{\bm \e}}, \ov{\bm \e}).
\end{equation}
In particular, $\mf{tr}_{\e}$ is an exact monoidal functor from $\mc{O}^{\mf b}_{\text{osc},{\bm\e}}$ to $\mc{O}^{\mf b}_{\text{osc},\e}$.
\end{thm}

\pf By Proposition~\ref{prop:decomp of W for e}, $\mc{V}^{(\la,\ell)}_{\bm\e}$ is a submodule of $\W_{\bm\e}^{\ot \ell}$ for $(\la,\ell)\in \mc{P}(\ell)_{\bm \e}$.
Its truncation $\mathfrak{tr}_{\epsilon}(\mc{V}^{(\la,\ell)}_{\bm\e})$ is semisimple, as it is in the semisimple $\mathring{\mathcal{U}}_{\epsilon}$-module $\W_{\e}^{\ot \ell}$.
Moreover, as we have seen in the proof of Proposition~\ref{prop:tensor product decomp}, the classical limit $\mc{V}^{(\la,\ell)}_{\bm\e}$ is isomorphic to $V((\lambda,\ell),\mathfrak{osp}_{B^\bullet}(\boldsymbol{\epsilon}))$ 
and so $\mathfrak{tr}_{\epsilon}(\mc{V}^{(\la,\ell)}_{\bm\e})$ has the same character as $\mathfrak{tr}_{\epsilon}\left(V((\lambda,\ell),\mathfrak{osp}_{B^\bullet}(\boldsymbol{\epsilon}))\right)$.
By \cite[Theorem~4.6]{CLW}, the latter is non-zero if and only if $(\lambda,\ell)\in\mathcal{P}(\ell)_\epsilon $, in which case it is isomorphic to $V((\lambda,\ell),\mathfrak{osp}_{B^\bullet}(\epsilon))$ and hence $\mathfrak{tr}_{\epsilon}(\mc{V}^{(\la,\ell)}_{\bm\e})$ is isomorphic to $\mc{V}^{(\la,\ell)}_{\e}$.

\qed

\section{Category {$\,^{\mathfrak{x}}\widehat{\mc O}^{\,\mf b}_{\rm osc,\e}$}}\label{sec:aff osc category type b}

\subsection{Category {$\,^{\mathfrak{x}}\widehat{\mc O}^{\,\mf b}_{\rm osc,\e}$}}\label{subsec:cat-O-hat}
{From now on, we extend our base field $\Bbbk$ to its algebraic closure $\mathbf{k}$ inside $\bigcup_{k>0} \mathbb{C}(q^{1/k})$ to which all of the preceding discussions extend naturally.}

We keep the notations from Section~\ref{sec:truncation}. 
For $\mf{x}=\mf{b}, \mf{d}$, let $\,^{\mathfrak{x}}\widehat{\mc O}^{\,\mf b}_{\rm osc,\e}$ denote the category of $\U_\e$-modules (corresponding to the  cases $X=B, D$ respectively), whose restriction to $\ring{\U}_\e$ belongs to ${\mc O}^{\,\mf b}_{\rm osc,\e}$.
It is a full subcategory of $\mc{C}(\e)$, which is closed under taking submodules, quotients and tensor product since $\mc{O}^{\mf b}_{\text{osc},\e}$ is a semisimple monoidal category by Proposition~\ref{prop:decomp of W for e} (cf.~\cite[Section 5.2]{KLO24}).

For instance, if $\e=\ov{\bm \e}$, then any $\U_{\ov{\bm \e}}$-module in $\,^{\mathfrak{x}}\widehat{\mc O}^{\,\mf b}_{\mathrm{osc},\ov{\boldsymbol{\e}}}$ is finite-dimensional.
In fact, it is the category of finite-dimensional  $U_{\td{q}}(D_{m+1}^{(2)})$-modules ($X=B$) or $U_{\td{q}}(B_{m}^{(1)})$-modules ($X=D$) as it contains the $q$-spin representation and its spectral parameter shifts (see Section~\ref{subsec:fusion-construction}).
On the other hand, for $\e=\un{\bm \e}$, $\,^{\mathfrak{x}}\widehat{\mc O}^{\,\mf b}_{\mathrm{osc},\un{\bm \e}}$ is a subcategory of an affinization of the BGG category \cite{MY} of $U_q(B_m)$ for $U_q(D_{m+2}^{(2)})$-modules ($X=B$) or $U_q(A_{2m}^{(2)\dagger})$-modules ($X=D$), provided that the set of weights is finitely dominated.

By Theorems \ref{thm:truncation functor} and \ref{thm:trunc sends simple to simple}, the restrictions of the truncation functors $\mf{tr}_{\un{\bm\e}}$ and $\mf{tr}_{\ov{\bm\e}}$ to $\,^{\mathfrak{x}}\widehat{\mc O}^{\,\mf b}_{\rm osc,{\bm\e}}$ give rise to exact monoidal functors to $\,^{\mathfrak{x}}\widehat{\mc O}^{\,\mf b}_{\rm osc,\un{\bm \e}}$ and $\,^{\mathfrak{x}}\widehat{\mc O}^{\,\mf b}_{\rm osc,\ov{\bm \e}}$, respectively:
\begin{equation}\label{eq:triangle}
\xymatrixcolsep{3pc}\xymatrixrowsep{0.3pc}\xymatrix{
 & \,^\mathfrak{x}\widehat{\mc O}^{\,\mf b}_{\rm osc,{\bm\e}}  \ar@{->}_{\mf{tr}_{\un{\bm\e}}}[dl]\ar@{->}^{\mf{tr}_{\ov{\bm\e}}}[dr] &  \\
 \,^\mathfrak{x}\widehat{\mc O}^{\,\mf b}_{\rm osc,\un{\bm\e}} & &  \,^\mathfrak{x}\widehat{\mc O}^{\,\mf b}_{\rm osc,\ov{\bm\e}}.
}
\end{equation}

We shall construct, in a uniform way, a family of irreducible objects in ${}^{\mf{x}}\widehat{\mc O}^{\,\mf b}_{\rm osc,\e}$ for each $\epsilon=\boldsymbol{\epsilon},\underline{\boldsymbol{\epsilon}},\overline{\boldsymbol{\epsilon}}$, which corresponds to one another under the truncation functors.
In particular, when $\epsilon=\overline{\boldsymbol{\epsilon}}$, the resulting family recovers all the finite-dimensional $\mathcal{U}_{\overline{\boldsymbol{\epsilon}}}$-modules. 
The construction is based on the same strategy as in \cite{KLO24}, namely, by comparing the normalized $R$ matrices in the categories \eqref{eq:triangle} via truncations and then applying fusion construction.

\subsection{Normalized $R$-matrix and fusion construction}\label{subsec:normalized R}
Let us briefly recall the fusion construction of irreducible representations via $R$-matrices (\emph{cf.}~\cite{KKKO15}). We present it in a slightly more general setting to include infinite-dimensional representations.
Let $M_1,M_2\in \,^{\mathfrak{x}}\widehat{\mc O}^{\,\mf b}_{\rm osc,\e}$ be given. For an indeterminate $z_k$ ($k=1,2$), let $M_k(z_k)={\bf k}[z^{\pm 1}_k]\ot M_k$ be the ${\bf k}[z^{\pm 1}_k]\otimes \mathcal{U}_\epsilon$-module defined by
\[
e_i (g(z_k)\otimes v)=g(z_k)z^{\delta_{0,i}}_k\otimes (e_i v), \quad f_i (g(z_k)\otimes v)=g(z_k)z^{-\delta_{0,i}}_k\otimes (f_i v),
\]
for $i\in I$, $g(z)\in {\bf k}[z^{\pm 1}]$ and $v\in M_k$.
The weight space decomposition of $M_k$ naturally extends to a $P_{\mathrm{af}}$-weight space decomposition of $M_k(z_k)$ by setting $\boldsymbol{\mathrm{wt}}(z_k^{l} \otimes v)=l\boldsymbol{\delta}+\mathrm{wt}(v)$ ($l\in\Z$), where $\mathrm{wt}(v)\in P$ is identified with its image in $P_{\mathrm{af}}$ under the embedding in Section~\ref{subsec:Yamane QSA}.

Recall the quasi $R$-matrix $\Theta$ given in Section~\ref{subsec:quasi-R}. Let $\Pi_\mathbf{q}:M_1\otimes M_2\rightarrow M_1 \otimes M_2$ be defined by $\Pi_\mathbf{q}(v_1\otimes v_2)=\mathbf{q}(\mathrm{wt}(v_1),\mathrm{wt}(v_2))(v_1\otimes v_2)$ for $v_k\in M_k$.
Then we have an isomorphism of $\mathcal{U}_\epsilon$-modules, called the universal $R$-matrix, given by the composition
\begin{equation*}
\mathcal{R}^{\mathrm{univ}}_{M_1,M_2}: M_1(z_1)\widehat{\otimes} M_2(z_2)\overset{\Theta}{\longrightarrow} M_1(z_1)\widehat{\otimes} M_2(z_2)\overset{\Pi_\mathbf{q}^{-1}}{\longrightarrow} M_1(z_1)\widehat{\otimes} M_2(z_2)\overset{s}{\longrightarrow} M_2(z_2)\widetilde{\otimes} M_1(z_1),
\end{equation*}
where $M_1(z_1)\widehat{\otimes} M_2(z_2)$ and $M_2(z_2)\widetilde{\otimes} M_1(z_1)$ are completions as given in \cite[Section~7]{Kas02} and $s$ is the flip map.
It satisfies the following triangle identities:
\begin{equation}\label{eq:triangle-id}
\begin{split}
\mathcal{R}^{\mathrm{univ}}_{M_1\otimes M_2,N}
&=
\left(
\mathcal{R}^{\mathrm{univ}}_{M_1,N}
\otimes \mathrm{id}_{M_2}
\right)
\circ
\left(
\mathrm{id}_{M_1}
\otimes
\mathcal{R}^{\mathrm{univ}}_{M_2,N}
\right),\\
\mathcal{R}^{\mathrm{univ}}_{N,M_1\otimes M_2}
&=
\left(
\mathrm{id}_{M_1}
\otimes
\mathcal{R}^{\mathrm{univ}}_{N,M_2}
\right)
\circ
\left(
\mathcal{R}^{\mathrm{univ}}_{N,M_1}
\otimes \mathrm{id}_{M_2}
\right),
\end{split}
\end{equation}
for $N \in \,^{\mathfrak{x}}\widehat{\mc O}^{\,\mf b}_{\rm osc,\e}$.

Suppose that $M_1$ and $M_2$ are irreducible and
\begin{center}
($\ast$)\quad there exists $\lambda_i\in P$ such that $(M_i)_{\lambda_i}=\mathbf{k}v_i$ and $\mathrm{wt}(M_i)\subset \lambda_i - \sum_{k\in\mathring{I}}\mathbb{Z}_{+}\alpha_k$.
\end{center}
When $\e=\ov{\bm\e}$, we may choose the dominant extremal weight vector $v_i$ of $M_i$ and its weight $\lambda_i$ (\emph{cf.}~\cite[Section 8]{Kas02}).
By weight comparison we have $\mathcal{R}^{\mathrm{univ}}_{M_1,M_2}(v_1 \otimes v_2)= a_{M_1,M_2}(z) (v_2 \otimes v_1)$ for some $a_{M_1,M_2}(z)\in\mathbf{k}\llbracket z \rrbracket^\times$ ($z=z_1/z_2$).
Hence we may define a \emph{normalized $R$-matrix} by 
{\small
\begin{equation}\label{eq:normalized R matrix general}
\mathcal{R}_{M_1,M_2}(z):=\frac{\mathcal{R}^{\mathrm{univ}}_{M_1,M_2}}{a_{M_1,M_2}(z)}:M_1(z_1) \otimes M_2(z_2) \longrightarrow \mathbf{k}\llbracket z \rrbracket\otimes_{\mathbf{k}[(z_1/z_2)^{\pm1}]} (M_2 (z_2)\otimes M_1(z_1)),
\end{equation}}
so that $\mathcal{R}_{M_1,M_2}(z)(v_1 \otimes v_2)=v_2 \otimes v_1$.
Note that if $\mathbf{k}(z_1,z_2)\otimes_{\mathbf{k}[z^{\pm1}_1,z^{\pm1}_2]} (M_1(z_1)\otimes M_2(z_2))$ is generated by $v_1 \otimes v_2$ as a $\mathbf{k}(z_1,z_2)\otimes\mathcal{U}_\epsilon$-module,
then $\mathcal{R}_{M_1,M_2}(z)$ is the unique $\mathbf{k}[z^{\pm 1}_1,z^{\pm 1}_2]\otimes\mathcal{U}_\epsilon$-module homomorphism normalized by $v_1\otimes v_2  \mapsto v_2 \otimes v_1$, and we also have ${\rm Im}\mathcal{R}_{M_1,M_2}(z) \subset \mathbf{k}(z_1,z_2)\otimes_{\mathbf{k}[z^{\pm1}_1,z^{\pm1}_2]} (M_2 (z_2)\otimes M_1(z_1))$.
In particular, if $M_1,M_2,M_3\in \,^{\mathfrak{x}}\widehat{\mc O}^{\,\mf b}_{\rm osc,\e}$ are such irreducible modules, then the associated normalized $R$ matrices satisfy the following Yang--Baxter equation:
\begin{equation}\label{eq:YBE}
\mathcal{R}_{M_2,M_3}\left(\tfrac{z_2}{z_3}\right)
\mathcal{R}_{M_1,M_3}\left(\tfrac{z_1}{z_3}\right)
\mathcal{R}_{M_1,M_2}\left(\tfrac{z_1}{z_2}\right)
=
\mathcal{R}_{M_1,M_2}\left(\tfrac{z_1}{z_2}\right)
\mathcal{R}_{M_1,M_3}\left(\tfrac{z_1}{z_3}\right)
\mathcal{R}_{M_2,M_3}\left(\tfrac{z_2}{z_3}\right). 
\end{equation}

Let $M_1,\dots,M_t \in {}^{\mf{x}}\widehat{\mathcal{O}}^{\mathfrak{b}}_{\mathrm{osc},\epsilon}$ be irreducible $\mathcal{U}_\epsilon$-modules satisfying the condition ($\ast$) above so that they admit normalized $R$-matrices \eqref{eq:normalized R matrix general}
\[
\mathcal{R}_{M_i,M_j}(z_i/ z_j): M_i (z_i) \otimes M_j (z_j) \longrightarrow \mathbf{k}\llbracket z_i/ z_j \rrbracket \otimes_{{\bf k}[(z_i/z_j)^{\pm1}]} (M_j (z_j) \otimes M_i (z_i)).
\]
Then we may define the $\U_\e$-linear map
\begin{equation*}
\begin{split}
& \mathcal{R}_{M_1,\dots,M_t}(z_1,\dots,z_t):\\ 
& \quad M_1(z_1)\otimes\cdots\otimes M_t(z_t) \longrightarrow \mathbf{k}\llbracket\, z_i/z_j\,|\,i<j \,\rrbracket\otimes_{\mathbf{k}[(z_i/z_j)^{\pm 1}|i<j]}(M_t(z_t)\otimes\cdots\otimes M_1(z_1)) 
\end{split}
\end{equation*}
by taking any reduced expression of the longest element of the symmetric group $S_t$ and composing the normalized $R$-matrices accordingly, which is well-defined due to \eqref{eq:YBE}.

Let $c_1,\dots,c_t \in \mathbf{k}^\times$ be given.
We further suppose that there exists
$b_{M_i,M_j}(z)\in\mathbf{k}\llbracket z\rrbracket$ for $i\le j$ such that the maps
\begin{equation*}\label{eq:r matrix}
\mathbf{r}_{M_i,M_j}\left(z_i/z_j\right)
=
b_{M_i,M_j}\left(z_i/z_j\right)
\mathcal{R}_{M_i,M_j}\left(z_i/z_j\right)
\end{equation*}
satisfy the following:
  \begin{itemize}
    \item[(1)] $\mathbf{r}_{M_i,M_j}(z_i/z_j)$ admits a non-zero specialization at $(z_i,z_j)=(c_i,c_j)$, that is,
\begin{equation}\label{eq:renormalization of R matrix}
\mathbf{r}_{M_i,M_j}(c_i/c_j): M_i(c_i)\ot M_j(c_j) \longrightarrow M_j(c_j)\ot M_i(c_i)
\end{equation}
is a non-zero $\U_\e$-linear map called a {\em renormalized $R$-matrix},
    
    \item[(2)] $\mathbf{r}_{M_i,M_i}(1)\in \mathbf{k}^\times \mathrm{id}_{M_i(c_i)^{\ot 2}}$ for all $i$.
\end{itemize}
Replacing $\mathcal{R}_{M_i,M_j}(z_i/ z_j)$ in $\mathcal{R}_{M_1,\dots,M_t}(z_1,\dots,z_t)$ with $\mathbf{r}_{M_i,M_j}\left(z_i/z_j\right)$ and then specializing $z_i=c_i$ ($i=1,\dots,t$), we obtain a $\mathcal{U}_\epsilon$-linear map
\[
\mathbf{r}_{M_1,\dots,M_t}(c_1,\dots,c_t):M_1(c_1)\otimes\cdots\otimes M_t(c_t) \longrightarrow M_t(c_t)\otimes\cdots\otimes M_1(c_1).
\]

\begin{thm}[{\cite{KKKO15}}]\label{thm:fusion-construction}
The image of $\mathbf{r}_{M_1,\dots,M_t}(c_1,\dots,c_t)$ is irreducible if it is not zero.
\end{thm}

\subsection{Normalized $R$ matrix on $\W_{\e}(x)\ot \W_{\e}(y)$}\label{subsec:norm R matrix for osc B} 
We first need the following.
\begin{prop}\label{prop:generic-simple-spin}
The tensor product
$\mathcal{W}_{\epsilon}(x)\otimes\mathcal{W}_{\epsilon}(y)$ is an irreducible $\U_\e$-module
for generic $x,y\in {\bf k}^{\times}$. Moreover,  $\mathcal{W}_{\epsilon}(x)^{\ot 2}=\mathcal{W}_{\epsilon}(x)\otimes\mathcal{W}_{\epsilon}(x)$ is irreducible for any $x\in {\bf k}^\times$.
\end{prop}
\pf The case when $\e=\ov{\bm\e}$ is well-known \cite{Kas02}.
The case of $\epsilon=\underline{\boldsymbol{\epsilon}}$ immediately follows from the subsequent lemmas.
Then the case when $\epsilon=\boldsymbol{\epsilon}$ follows from the arguments in \cite[Proof of Theorem 7.6]{KLO24} using Theorem~\ref{thm:trunc sends simple to simple}.
\qed
\smallskip

Recall that we have the following decomposition as a $\mathring{\mathcal{U}}_{\epsilon}$-module: 
\begin{equation}\label{eq:spin2 decomposition}
\mathcal{W}_{\epsilon}\otimes\mathcal{W}_{\epsilon}\cong\bigoplus_{i\geq0}\mathcal{V}_{\epsilon}^{((i,0),2)}, 
\end{equation}
(independent of $X=B, D$) by  Proposition~\ref{prop:decomp of W for e} when $\ell=2$.

\begin{lem}
For $i\geq0$, the vector
\begin{equation*}\label{eq: maximal vectors}
v_{i}=\sum_{a=0}^{i}\left[(-1)^{a}\mathbf{v}^{-a}\prod_{k=1}^{a}q^{k-1}\frac{[i-k+1]}{[k]}\right]\ket{(i-a)\mathbf{e}_{n}}\otimes\ket{a\mathbf{e}_{n}} 
\end{equation*}
is a $\mathring{\mathcal{U}}_{\underline{\boldsymbol{\epsilon}}}$-highest
weight vector of the irreducible component $\mathcal{V}_{\underline{\boldsymbol{\epsilon}}}^{((i,0),2)}$
in $\mathcal{W}_{\underline{\boldsymbol{\epsilon}}}\otimes\mathcal{W}_{\underline{\boldsymbol{\epsilon}}}$.
\end{lem}
\pf It is straightforward to check that $v_i$ is annihilated by $e_j$ for $j\in \ring{I}$.
\qed

\begin{lem}\label{lem:E,F-on-hw-vectors} Let $v_i \in \mathcal{W}_{\underline{\boldsymbol{\epsilon}}}(x)\otimes\mathcal{W}_{\underline{\boldsymbol{\epsilon}}}(y)$  $(i\geq 0)$ be the $\mathring{\mathcal{U}}_{\un{\bm\epsilon}}$-highest weight vector of $\mathcal{V}_{\underline{\boldsymbol{\epsilon}}}^{((i,0),2)}$.
  \begin{enumerate}
    \item If $X=B$, then we have
      \begin{align*}
        f_0 f_1 \cdots f_{m} v_i  & = [i](x^{-1}+q^i y^{-1})v_{i-1},\\
        e_m \cdots e_1 e_0 v_i & = h_{i,0}(x,y) v_{i+1} + h_{i,1}(x,y) f_{m+1}v_{i} + h_{i,2} (x,y) f^{2}_{m+1}v_{i-1},
      \end{align*}
for some $h_{i,k}(z_1,z_2) \in \Bbbk[z_1,z_2]$ with $h_{i,0}(c,c) \neq 0$   for any $c\in\mathbf{k}^\times$.
    \item If $X=D$, then we have
       \begin{align*}
        f_0 f^{(2)}_1 \cdots f^{(2)}_{m-1} v_i  & = \frac{[i][i-1]}{[2]^2}(x^{-1}+q^{2i-1}y^{-1})v_{i-2} \quad (i\geq 2),\\
        f_0 f^{(2)}_1 \cdots f^{(2)}_{m-1} f_m v_1  & = \frac{x^{-1}+qy^{-1}}{[2]}v_0,\\
        e_m e^{(2)}_{m-1}\cdots e^{(2)}_1 e_0 v_0 & = \frac{1}{\mathbf{v}(q-1)[2]}((q^{-1}x+y)v_1 - (x+y)f_m v_0), \\
        e^{(2)}_{m-1}\cdots e^{(2)}_1 e_0 v_i & = \sum^{4}_{k=0} h^{\prime}_{i,k} (x,y)f^{k}_m v_{i+2-k},
      \end{align*}
for some $h^{\prime}_{i,k} (z_1,z_2) \in \Bbbk[z_1,z_2]$ with $h^{\prime}_{i,0}(c,c) \neq 0$ for any $c\in\mathbf{k}^\times$.
  \end{enumerate}
\end{lem}

\pf
The formulas except the last ones in (1) and (2) can be verified by a straightforward calculation. 
For the remaining ones, say when $X=B$, the key step is to prove that
  \begin{align*}
  e_{m+1}^2 e_m\cdots e_1 e_0 v_i \in \mathbf{k}^\times h_{i,2}(x,y)v_{i-1},
  \end{align*}
  from which the claim can be deduced by considering the coefficient of $v_i$ in $e_{m+1} e_m\cdots e_1 e_0 v_i$ and then that of $v_{i+1}$ in $e_m\cdots e_1 e_0 v_i$.
  We refer the reader to \cite[Appendix~C]{KLO24} for more details where similar formulas are obtained for the quantum affine algebra of type $D^{(1)}_{m}$.
\qed
\smallskip

Now consider the normalized $R$-matrix on $\mathcal{W}_{\epsilon}(z_{1})\otimes\mathcal{W}_{\epsilon}(z_{2})$ with respect to $\mathcal{V}_{\epsilon}^{((0,0),2)}$:
\begin{equation}\label{eq:normalized R matrix}
\begin{split}
& \mathcal{R}_{\epsilon}(z)\, :\,\mathcal{W}_{\epsilon}(z_{1})\otimes\mathcal{W}_{\epsilon}(z_{2})\longrightarrow\Bbbk(z_{1},z_{2})\otimes_{\Bbbk[z_{1}^{\pm1},z_{2}^{\pm1}]}\left(\mathcal{W}_{\epsilon}(z_{2})\otimes\mathcal{W}_{\epsilon}(z_{1})\right),\\[2pt]
 & \mathcal{R}_{\epsilon}(z)\vert_{\mathcal{V}_{\epsilon}^{((0,0),2)}}=\mathrm{id}_{\mathcal{V}_{\epsilon}^{((0,0),2)}},
\end{split}
\end{equation}
where $z=z_1/z_2$.
The irreducibility of $\mathcal{W}_{\epsilon}(z_{1})\otimes\mathcal{W}_{\epsilon}(z_{2})$ in Proposition \ref{prop:generic-simple-spin} implies that $\mathcal{R}_{\epsilon}(z)$ is unique.
Then we have the following compatibility of normalized $R$-matrices since $\mf{tr}_{\e}(\mc{V}_{\bm\e}^{((0,0),2)})=\mc{V}_{\e}^{((0,0),2)}$  by Theorem~\ref{thm:trunc sends simple to simple}.

\begin{lem}\label{lem:truncation of R}
For $\epsilon=\underline{\boldsymbol{\epsilon}},\,\overline{\boldsymbol{\epsilon}}$,
we have
\[
\mathfrak{tr}_{\epsilon}(\mathcal{R}_{\boldsymbol{\epsilon}}(z))=\mathcal{R}_{\epsilon}(z).
\]
\end{lem}

We can now give the spectral decomposition of $\mc{R}_\e(z)$.
\begin{thm}\label{thm:spectral decomp}
We have the following:
\begin{align*}
\mathcal{R}_{\epsilon}(z)
=
\begin{cases}
\displaystyle
\sum_{i\geq0}\prod_{j=1}^{i}
\mfrac{1+q^{j}z}{z+q^{j}}
\mathcal{P}_{(i,0),\epsilon}
& (X=B),\\[6pt]
\displaystyle
\sum_{i\geq0}\prod_{j=1}^{i}
\mfrac{1+zq^{4j-1}}{z+q^{4j-1}}
\mathcal{P}_{(2i,0),\epsilon}
+
\sum_{i\geq0}\prod_{j=1}^{i}
\mfrac{1+zq^{4j-3}}{z+q^{4j-3}}
\mathcal{P}_{(2i-1,0),\epsilon}
& (X=D),
\end{cases}
\end{align*}where $\mathcal{P}_{(i,0),\epsilon}$ is the projection onto the irreducible component $\mathcal{V}_{\epsilon}^{((i,0),2)}$ in \eqref{eq:spin2 decomposition}.
\end{thm}

\pf
The formula for $\epsilon=\boldsymbol{\underline{\epsilon}}$ follows easily from Lemma~\ref{lem:E,F-on-hw-vectors} and the fact that the $R$-matrix is $\mathcal{U}_\epsilon$-linear. It was originally obtained in \cite[Proposition~7]{KO} for $X=B$ and \cite[Proposition~B.1]{KwO21} for $X=D$ (with a different comultiplication).
Then by Lemma \ref{lem:truncation of R} and Theorem~\ref{thm:trunc sends simple to simple}, we obtain the formula for $\epsilon=\boldsymbol{\epsilon}$ and subsequently for $\epsilon=\boldsymbol{\overline{\epsilon}}$ by applying $\mathfrak{tr}_{\overline{\boldsymbol{\epsilon}}}^{\boldsymbol{\epsilon}}$.
\qed
\smallskip

\begin{rem}\label{rem:denominator}
{\rm
In the case of $\epsilon=\overline{\boldsymbol{\epsilon}}$, we see that $\mathcal{P}_{(i,0),\epsilon}=0$ if $i > m$, and hence the summation in the spectral decomposition of $\mc{R}_{\epsilon}(z)$ is finite. So we may regard the finite product $\prod_{j=1}^{m}(z+q^j)$ (resp. $\prod_{j=1}^{n} (z+q^{2j-1})$) when $X=B$ (resp. $X=D$) as the \emph{denominator} of the normalized $R$-matrix.
Indeed, $\mathcal{R}_{\epsilon}(z)$ coincides with the normalized $R$-matrix for $V(\upvarpi_m)_x\ot V(\upvarpi_m)_y$  under the identification 
\begin{equation}\label{eq:identifying with spin}
\mathcal{W}_{\overline{\boldsymbol{\epsilon}}}(x)=V(\upvarpi_m)_x, 
\end{equation}
which is the $q$-spin representation with spectral parameter $x$ of type $D^{(2)}_{m+1}$ (resp. $B^{(1)}_m$) when $X=B$ (resp. $X=D$) (see Remark~\ref{rem:truncation of W}).
Then our denominator (or even the spectral decomposition itself) agrees with the corresponding formula in \cite[Proposition~4.9.5]{KKMMNN} (resp. \cite[(5.9)]{O90}) once the conventions $q_m = q^{1/2}$ and $\widetilde{q}=-q^{-1}$ for the quantum parameter of $\mathcal{U}_{\overline{\boldsymbol{\epsilon}}}$ are taken into account. 
} 
\end{rem}

\subsection{Fundamental representations and fusion construction  in $\,^{\mathfrak{x}}\widehat{\mc O}^{\,\mf b}_{\rm osc,\e}$}\label{subsec:fusion-construction}\mbox{}

Let us first introduce representations in $^{\mathfrak{x}}\widehat{\mathcal O}^{,\mathfrak b}_{\mathrm{osc},\epsilon}$ that play a role analogous to that of fundamental representations in the finite-dimensional setting.

For a positive integer $l>0$ and $x\in {\bf k}^\times$, we define a $\mathcal{U}_{\epsilon}$-module $\mc{W}_{l,\e}(x)$ as follows:
{\small
\begin{gather*}
\mathcal{W}_{l,\boldsymbol{\epsilon}}(x)  = 
\begin{cases}
  {\rm Ker}\left(\mathcal{R}_{\boldsymbol{\epsilon}}(-q^{-l})\,:\,\mathcal{W}_{\boldsymbol{\epsilon}}(x)\otimes\mathcal{W}_{\boldsymbol{\epsilon}}(-xq^{l})\longrightarrow\mathcal{W}_{\boldsymbol{\epsilon}}(-xq^{l})\otimes\mathcal{W}_{\boldsymbol{\epsilon}}(x)\right) & (X=B),\\
 {\rm Ker}\left(\mathcal{R}_{\boldsymbol{\epsilon}}(-q^{-2l+1})\,:\,\mathcal{W}_{\boldsymbol{\epsilon}}(x)\otimes\mathcal{W}_{\boldsymbol{\epsilon}}(-xq^{2l-1})\longrightarrow\mathcal{W}_{\boldsymbol{\epsilon}}(-xq^{2l-1})\otimes\mathcal{W}_{\boldsymbol{\epsilon}}(x)\right) & (X=D),
\end{cases}\\[3pt]
\mathcal{W}_{l,\epsilon}(x) =\mathfrak{tr}_{\epsilon}^{\boldsymbol{\epsilon}}\left(\mathcal{W}_{l,\boldsymbol{\epsilon}}(x)\right)\quad (\epsilon=\underline{\boldsymbol{\epsilon}},\,\overline{\boldsymbol{\epsilon}}).   
\end{gather*}}
We also set $\mathcal{W}_{0,\epsilon}(x)= \mathcal{W}_{\epsilon}(x)$.

\begin{prop}\label{prop:irreducibility of W and fundamental}
We have the following:
\begin{itemize}
\item[(1)] For $l>0$ and $x\in {\bf k}^\times$, $\mathcal{W}_{l,\epsilon}(x)$ is either zero or irreducible. Moreover, it is zero if and only if $\epsilon=\overline{\boldsymbol{\epsilon}}$ and $l>m$.

\item[(2)] For $l_1, l_2\in\Z_+$, $\mathcal{W}_{l_1,\epsilon}(x)\otimes \mathcal{W}_{l_2,\epsilon}(y)$ is irreducible for generic $x, y\in {\bf k}^\times$.
\end{itemize}
\end{prop}

\pf (1)
Consider the following renormalization of $\mathcal{R}_{\boldsymbol{\epsilon}}(z)$:
\begin{gather*}
\mathcal{R}_{\boldsymbol{\epsilon}}^{\diamond}(z)=
\begin{cases}
\prod_{j=1}^{l}(z+q^{j})\,\mathcal{R}_{\boldsymbol{\epsilon}}(z) & (X=B),\\
\prod_{j=1}^{l}(z+q^{2j-1})\,\mathcal{R}_{\boldsymbol{\epsilon}}(z) & (X=D),
\end{cases}\\
\mathcal{R}_{\epsilon}^{\diamond}(z)=\mathfrak{tr}_{\epsilon}\left( \mathcal{R}_{\boldsymbol{\epsilon}}^{\diamond}(z)\right)\quad (\epsilon=\underline{\boldsymbol{\epsilon}},\overline{\boldsymbol{\epsilon}}).
\end{gather*}
Then it follows from Lemma \ref{lem:truncation of R} and Theorem \ref{thm:spectral decomp} that
{\small
\begin{equation*}
\mathcal{W}_{l,\epsilon}(x)=
\begin{cases}
{\rm Im}\left(\mathcal{R}_{\epsilon}^{\diamond}(-q^{l}):\mathcal{W}_{\epsilon}(-xq^{l})\otimes\mathcal{W}_{\epsilon}(x)\longrightarrow\mathcal{W}_{\epsilon}(x)\otimes\mathcal{W}_{\epsilon}(-xq^{l})\right) & (X=B),\\
{\rm Im}\left(\mathcal{R}_{\epsilon}^{\diamond}(-q^{2l-1}):\mathcal{W}_{\boldsymbol{\epsilon}}(-xq^{2l-1})\otimes\mathcal{W}_{\boldsymbol{\epsilon}}(x) \longrightarrow \mathcal{W}_{\boldsymbol{\epsilon}}(x)\otimes\mathcal{W}_{\boldsymbol{\epsilon}}(-xq^{2l-1})\right) & (X=D).  
\end{cases}
\end{equation*}}

Let $M_k=\mathcal{W}_{\epsilon}(1)$ ($k=1,2$) with $(c_1,c_2)=(-xq^l,x)$ and (resp. $(c_1,c_2)=(-xq^{2l-1},x)$) when $X=B$ (resp. $X=D$), and let ${\bf r}_{M_1,M_2}(z_1/z_2)=\mathcal{R}_{\boldsymbol{\epsilon}}^{\diamond}(z_1/z_2)$ and ${\bf r}_{M_k,M_k}(z_k/z_k)=\mc{R}_{\e}(1)={\rm id}_{M_k^{\ot 2}}$.
Then ${\bf r}_{M_i,M_j}(c_i/c_j)$ $(i\le j)$ is a renormalized $R$-matrix in the sense of \eqref{eq:renormalization of R matrix}. Hence by Theorem \ref{thm:fusion-construction}, $\mathcal{W}_{l,\epsilon}(x)$ is irreducible when it is non-zero.
Finally, we see from Theorem \ref{thm:spectral decomp} that as a $\ring{\U}_\e$-module
\begin{equation*}\label{eq:classical decomposition of fundamental rep}
\mathcal{W}_{l,\epsilon}(x)\cong
\begin{cases}
\bigoplus_{i\geq0}\mathcal{V}_{\epsilon}^{((l+i,0),2)} & (X=B),\\
\bigoplus_{i\geq0}\mathcal{V}_{\epsilon}^{((l+2i,0),2)} & (X=D),
\end{cases}
\end{equation*}
and hence it vanishes if and only if $\epsilon=\overline{\boldsymbol{\epsilon}}$ and $l>m$.

(2) It follows from the generic irreducibility of $\mathcal{W}_{\epsilon}(x)\otimes\mathcal{W}_{\epsilon}(y)$ (Proposition~\ref{prop:generic-simple-spin}), see \cite[Proposition~6.8]{KLO24}.
\qed
\smallskip

Suppose $\epsilon=\overline{\boldsymbol{\epsilon}}$ and $0<l<m$.
The maximal weight of $\mathcal{W}_{l,\overline{\boldsymbol{\epsilon}}}(x)$, namely the weight of its $\ell$-highest weight vector, is $\upvarpi_{((l,0),2),\mf{so}_{2m+1}}=\upvarpi_{m-l}$ (cf.~Remark \ref{rem: decomposition for super} and \eqref{eq:identifying with spin}).
Therefore, 
\begin{equation}\label{eq:identifying with fundamental}
\mathcal{W}_{l,\overline{\boldsymbol{\epsilon}}}(x)\cong V(\upvarpi_{m-l})_{x},
\end{equation}
which is the $(m-l)$-th fundamental representation of the quantum affine algebra $\mathcal{U}_{\overline{\boldsymbol{\epsilon}}}=U_{\widetilde{q}}(D^{(2)}_{m+1})$ (resp. $U_{\widetilde{q}}(B^{(1)}_{m})$) when $X=B$ (resp. $X=D$), up to a shift of spectral parameter. So it is natural to call $\mathcal{W}_{l,\epsilon}(x)$ a \emph{fundamental $q$-oscillator representation}  for $\epsilon=\boldsymbol{\epsilon},\underline{\boldsymbol{\epsilon}}$. 
\smallskip

Next, consider a normalized $R$-matrix on $\W_{l_1,\e}(z_1)\ot \W_{l_2,\e}(z_2)$ for $l_1, l_2\in \Z_+$.
Let $(\la_0,\ell)\in \mc{P}(\ell)$ be given by
\begin{equation}\label{eq:canonical maximal component}
(\la_0,\ell)=
\begin{cases}
((\max\{l_{1},l_{2}\},\min\{l_{1},l_{2}\},0,0),4) & \text{if $\min\{l_1,l_2\}>0$},\\
((\max\{l_{1},l_{2}\},0,0),3) & \text{if $\max\{l_1,l_2\}>\min\{l_1,l_2\}=0$},\\ 
((0,0),2) & \text{if $\max\{l_1,l_2\}=0$}.
\end{cases}
\end{equation}
Suppose that $(\la_0,\ell)\in \mc{P}(\ell)_\e$. Then we may consider the normalized $R$-matrix on $\W_{l_1,\e}(z_1)\ot \W_{l_2,\e}(z_2)$ with respect to $\mc{V}_{\e}^{(\la_{0},\ell)}$, that is,
\begin{equation*}
\begin{split}
&\mc{R}_{l_{1},l_{2},\e}(z):\W_{l_1,\e}(z_1)\ot \W_{l_2,\e}(z_2)\longrightarrow {\bf k}(z_{1},z_{2})\otimes_{\mathbf{k}[z^{\pm1}_1,z^{\pm1}_2]} \W_{l_2,\e}(z_2)\ot \W_{l_1,\e}(z_1),\\
&\mathcal{R}_{l_{1},l_{2},\e}(z)\vert_{\mc{V}_{\e}^{(\la_{0},\ell)}}={\rm id}_{\mc{V}_{\e}^{(\la_{0},\ell)}},
\end{split}
\end{equation*} 
with $z=z_1/z_2$.

\begin{lem}\label{lem:tr-pres-R}
We have $\mathfrak{tr}_{\epsilon}\left(\mathcal{R}_{l_{1},l_{2},\boldsymbol{\epsilon}}(z)\right)=\mathcal{R}_{l_{1},l_{2},\epsilon}(z)$ for $\epsilon=\underline{\boldsymbol{\epsilon}},\,\overline{\boldsymbol{\epsilon}}$.
.
\end{lem}

\pf It follows from Theorem \ref{thm:trunc sends simple to simple} and Proposition \ref{prop:irreducibility of W and fundamental}(2) by the same arguments as in Lemma \ref{lem:truncation of R}.
\qed
\smallskip

\begin{lem}\label{lem:unicity of R}
For $l\in\Z_+$, $\mathcal{R}_{l,l,\epsilon}(z)$ does not have a pole at $z$=1 and satisfies $\mathcal{R}_{l,l,\epsilon}(1)= \mathrm{id}_{\W_{l,\e}(z)^{\ot 2}}$.
\end{lem}
\pf Let us first compare the spectral decompositions of $\mathcal{R}_{l_{1},l_{2},\epsilon}(z)$ for $\e={\bm\e}, \un{\bm\e}, \ov{\bm\e}$.
Let $S_\e$ be the set of $(\la,\ell)\in \mc{P}(\ell)_{\e}$ such that $\mc{V}^{(\la,\ell)}_\e$ appears in $\W_{l_{1},\e}\ot \W_{l_{2},\e}$ and hence in $\W_{l_{2},\e}\ot \W_{l_{1},\e}$ (this set is independent of the omitted spectral parameters). Note that $S_\e\subset S_{\bm\e}$ for $\e=\un{\bm\e}, \ov{\bm\e}$.

For each $(\la,\ell)\in S_{\e}$, let us fix $\mathring{\mathcal{U}}_{\e}$-highest weight vectors $v_{\lambda,i,\e}$ and $v_{\lambda,i,\e}^{\prime}$ of $\mathcal{V}^{(\lambda,\ell)}_{\epsilon}$
in $\W_{l_{1},\e}\ot \W_{l_{2},\e}$ and $\W_{l_{2},\e}\ot \W_{l_{1},\e}$ respectively.
Here we put additional index $i$ to
emphasize that the multiplicity of $\mathcal{V}_{\overline{\boldsymbol{\epsilon}}}^{(\lambda,\ell)}$ can be greater than one in general.
Let 
\begin{equation}\label{eq:projection}
\mathcal{P}_{\lambda,i,j,\e}^{l_{1},l_{2}}:\mathcal{W}_{l_{1},\e}\otimes\mathcal{W}_{l_{2},\e}\longrightarrow\mathcal{W}_{l_{2},\e}\otimes\mathcal{W}_{l_{1},\e} 
\end{equation}
be the projection as a $\mathring{\mathcal{U}}_{\e}$-linear
map given by $\mathcal{P}_{\lambda,i,j,\e}^{l_{1},l_{2}}(v_{\mu,k,\e})=\delta_{\lambda\mu}\delta_{ik}v_{\lambda,j,\e}^{\prime}$.
We also use the notation $\mathcal{P}_{\lambda,i,j,\e}^{l_{1},l_{2}}=0$ if $(\lambda,\ell)\notin \mathcal{P}(\ell)_\epsilon$.
Then we may write 
\begin{equation}\label{eq:spectral decomp}
\mathcal{R}_{l_{1},l_{2},\e}(z)=\sum_{\lambda,i,j}\rho_{\lambda,i,j,\e}^{l_{1},l_{2}}(z)\mathcal{P}_{\lambda,i,j,\e}^{l_{1},l_{2}},
\end{equation}
for some $\rho^{l_1, l_2}_{\lambda,i,j,\e}(z)\in {\bf k}(z_1,z_2)$.

By Proposition \ref{prop:tensor product decomp}, the multiplicity of $\mc{V}^{(\la,\ell)}_\e$ in $\W_{l_1,\e}\ot \W_{l_2,\e}$ is independent of $\e$ once it is non-zero.
So by Theorem~\ref{thm:trunc sends simple to simple}, we may assume
$\mf{tr}_{\e}\left(\mc{P}_{\la,i,j,\bm{\e}}^{l_{1},l_{2}}\right)=\mc{P}_{\la,i,j,{\e}}^{l_{1},l_{2}}$ ($\e=\un{\bm\e}, \ov{\bm\e}$) for $(\la,\ell)\in S_{\bm\e}$.

On the other hand, $\rho^{l_1, l_2}_{\lambda,i,j,{\e}}(z)\in {\bf k}(z_1,z_2)$ does not depend on the rank of ${\e}$, once it is non-zero (which can be proved by applying the same arguments as in \cite[Section 4.2]{KL20}).
Furthermore, by Lemma~\ref{lem:tr-pres-R}, the non-zero coefficient $\rho_{\lambda,i,j,{\e}}^{l_{1},l_{2}}(z)$ is independent of $\e$. 
Summarizing, we may write
\begin{equation}\label{eq:normalized R for two fundamental}
\mathcal{R}_{l_{1},l_{2},\epsilon}(z)=\sum_{\lambda,i,j}\rho_{\lambda,i,j}^{l_{1},l_{2}}(z)\mathcal{P}_{\lambda,i,j,\epsilon}^{l_{1},l_{2}} \quad (\epsilon={\bm\e}, \un{{\bm\e}}, \ov{{\bm\e}}), 
\end{equation}
for some $\rho_{\lambda,i,j}^{l_{1},l_{2}}(z)\in {\bf k}(z_1,z_2)$, where the sum is over $(\la,\ell)\in\mc{P}(\ell)$.

It is well-known that $\mc{R}_{l,l,\ov{\bm\e}}(1)={\rm id}_{\W_{l,\ov{\bm\e}}(z)^{\ot 2}}$ \cite{KKKO15} (cf.~\eqref{eq:identifying with fundamental}), which is equivalent to $\rho_{\lambda,i,j}^{l,l}(1)=\delta_{i,j}$ for all $i,j$.
This together with \eqref{eq:normalized R for two fundamental} implies that
$\mathcal{R}_{l,l,\epsilon}(1)=\operatorname{id}_{\W_{l,\e}(z)^{\ot 2}}$ 
for $\epsilon=\boldsymbol{\epsilon},\underline{\boldsymbol{\epsilon}}$.
\qed
\smallskip

Now, we apply the fusion construction.
For $l_1,l_2 \in \Z_+$, let $\mc{Z}^{(l_1,l_2)}_\e$ be the set of poles in $\mc{R}_{l_1,\l_2,\e}(z)$.
Let $\mc{P}^+_\e$ be the set of $({\mb l},{\mb c})$ such that
\begin{itemize}
\item[(1)]  ${\mb l}=(l_{1},\dots,l_{t})\in\Z_+^{t}$ and $\mb{c}=(c_1,\dots,c_t)\in ({\bf k}^\times)^t$ for some $t\geq 1$,

\item[(2)] $c_{i}/c_{j}\not\in \mc{Z}^{(l_i,l_j)}_\e$ for all $i<j$ when $t\ge 2$.
\end{itemize}
By definition of $\mc{Z}^{(l_i,l_j)}_\e$ and Lemma \ref{lem:unicity of R}, the maps 
\begin{equation*}
{\bf r}_{\W_{l_i,\epsilon},\mathcal{W}_{l_j,\epsilon}}(z_i/z_j)\coloneqq \mathcal{R}_{l_i,l_j,\epsilon}(z_i/z_j)\quad (i\le j) 
\end{equation*}
satisfy the conditions for Theorem \ref{thm:fusion-construction} ($M_i=\W_{l_i,\epsilon}$).
Hence for $({\mb l},{\mb c})\in \mc{P}^+_\e$, we have 
a $\U_\e$-linear map
\[
\mathbf{r}_{\mathcal{W}_{l_1,\epsilon},\dots,\mathcal{W}_{l_t,\epsilon}}(c_1,\dots,c_t):\mathcal{W}_{l_1,\epsilon}(c_1)\otimes\cdots\mathcal{W}_{l_t,\epsilon}(c_t)\longrightarrow \mathcal{W}_{l_t,\epsilon}(c_t)\otimes\cdots\otimes\mathcal{W}_{l_1,\epsilon}(c_1),
\]
whose image we denote by $\W_{\e}({\mb l},{\mb c})$.
Summarizing, we have the following.
\begin{thm}
For $({\mb l},{\mb c})\in \mc{P}^+_\e$, we have the following.
\begin{itemize}
\item[(1)] $\W_{\e}({\mb l},{\mb c})$ is either zero or an irreducible $\U_\e$-module in $\,^{\mathfrak{x}}\widehat{\mc O}^{\,\mf b}_{\rm osc,\e}$,

\item[(2)] $\mf{tr}_{\e}\left(\W_{\bm \e}({\mb l},{\mb c})\right)=\W_{\e}({\mb l},{\mb c})$ for $\e=\un{\bm \e},\ov{\bm\e}$.
\end{itemize}
\end{thm}
\pf (1) follows from Theorem \ref{thm:fusion-construction}. (2) follows from Theorem \ref{thm:truncation functor} and Lemma \ref{lem:tr-pres-R}.
\qed
\smallskip
 
Let us make some further remarks on $\mathcal{Z}^{(l_1, l_2)}_{\epsilon}$, although an explicit formula for the spectral decomposition \eqref{eq:normalized R for two fundamental} is not known. 

\begin{prop}
For $\epsilon=\boldsymbol{\epsilon},\underline{\boldsymbol{\epsilon}}$, we have
\begin{equation*}
\mathcal{Z}^{(l_1,l_2)}_{\epsilon} = \bigcup_{M\ge 2}\mc{Z}^{(l_1,l_2)}_{(1^{M})}.
\end{equation*}
\end{prop}
\pf
The set $\mathcal{Z}^{(l_1,l_2)}_{\epsilon}$ is the union of the set of poles of $\rho^{l_1,l_2}_{\lambda,i,j}(z)$ for all $\lambda\in\mathcal{P}(\ell)_{\epsilon}$ and $i,j$ \eqref{eq:normalized R for two fundamental}.
Therefore, the claim follows from 
\[
\mathcal{P}(\ell)_\epsilon=\bigcup_{M\geq 2} \mathcal{P}(\ell)_{(1^M)}.
\]
\qed

\begin{rem}{\rm
Similarly we obtain $\mathcal{Z}^{(l_1,l_2)}_{(1^M)} \subset \mathcal{Z}^{(l_1,l_2)}_{(1^{M+1})}$. 
The precise description of the (finite) set $\mathcal{Z}^{(l_1,l_2)}_{(1^M)}$ can be read from the denominator formula for $\mathcal{R}_{l_1,l_2,(1^M)}(z)$ (see, for example, Theorem~6.7 for $X=B$ and \cite[Theorem 6.5]{OS} for $X=D$).
In particular, we have 
$
\mc{Z}^{(l_1,l_2)}_{(1^M)}\cap\left(\bigcup_{k>0}\mathbb{C}\llbracket q^{-1/k} \rrbracket \right) = \emptyset
$
which also follows from a general theory \cite[Proposition~9.3]{Kas02} without referring to the denominator formula
(note that $\mathcal{U}_{(1^M)}$ has the quantum parameter $-q^{-1}$).
Therefore, we also have
\begin{equation*}
\mc{Z}^{(l_1,l_2)}_{\e}\cap \bigcup_{k>0}\mathbb{C}\llbracket q^{-1/k} \rrbracket = \emptyset
\end{equation*}
for any $\epsilon=\boldsymbol{\epsilon},\underline{\boldsymbol{\epsilon}}$.
Finally, as $\mathcal{W}_{l_1,\epsilon}\otimes \mathcal{W}_{l_2,\epsilon}$ is of infinite length as a $\mathring{\mathcal{U}}_\epsilon$-module unless $\epsilon=\overline{\boldsymbol{\epsilon}}$,
the set $\mc{Z}^{(l_1,l_2)}_{\e}$ is infinite in general.
}
\end{rem}

\appendix
\section{Proof of \eqref{eq:Serre-hat{e}-B}}\label{sec:tr-pf-Serre}
In this appendix, we give a proof of the identity
\[
\hat{e}_{0}^{3}\hat{e}_{1}-[3]_{\mathbf{v}}(\hat{e}_{0}^{2}\hat{e}_{1}\hat{e}_{0}-\hat{e}_{0}\hat{e}_{1}\hat{e}_{0}^{2})-\hat{e}_{1}\hat{e}_{0}^{3} =0
\]
where $\hat{e}_{0}, \hat{e}_{1} \in \mathcal{U}_{BB}(\boldsymbol{\epsilon})$ ($\boldsymbol{\epsilon}=(0,1,0,1,\dots)$) are given by
\[
\hat{e}_{0}=e_0 e_1 -q^{\eta}e_1 e_0,\quad\hat{e}_{1}=e_2 e_3 -q^{\eta}e_3 e_2, \quad \eta=\pm 1.
\]
Since the proof is the same, we assume $\eta=1$.

We will use the following shorthand notation
\[
(i_{1}i_{2}\dots i_{l})\coloneqq e_{i_{1}}e_{i_{2}}\cdots e_{i_{l}}\in\mathcal{U}_{BB}(\epsilon).
\]
For example, the Serre relations in $\mathcal{U}_{BB}(\boldsymbol{\epsilon})$ we are going to invoke are written as
\begin{align}
(0001)-(q+1+q^{-1})(0010)+(q+1+q^{-1})(0100)-(1000) & =0\label{eq:even-Serre-0}\\
(1012)-(1210)+(q+q^{-1})(1021)-(0121)+(2101) & =0\label{eq:odd-Serre}\\
(11)=(22)=0,\quad(02)=(20).
\end{align}
Moreover, an overscore above a part of monomial $(i_1 \dots i_l)$ indicates which factor of $e_{i_1}\cdots e_{i_l}$ is replaced by the relation above.

Let us substitute the definition of $\hat{e}_0$ and $\hat{e}_1$ into the identity we shall prove.
Since $e_3$ commutes with $e_0$ and $e_1$ and so with $\hat{e}_0$, it suffices to check that
\begin{align*}
 & (1010102)-q^{-1}(1010012)+q^{-2}(1001012)-q^{-3}(0101012)\\
 & +(q-1+q^{-1})\begin{pmatrix}(1010210)-q^{-1}(1001210)+q^{-2}(0101210)\\
-q^{-1}(1010201)+q^{-2}(1001201)-q^{-3}(010\overline{1201})
\end{pmatrix}\\
 & -(q-1+q^{-1})\begin{pmatrix}(1021010)-q^{-1}(1021001)+q^{-2}(1020101)\\
-q^{-1}(0121010)+q^{-2}(0121001)-q^{-3}(0\overline{1201}01)
\end{pmatrix}\\
 & -(2101010)+q^{-1}(2101001)-q^{-2}(2100101)+q^{-3}(2010101)
\end{align*}
vanishes. Replacing each overscored factors above using (\ref{eq:odd-Serre}),
it is equal to 
\begin{align*}
 & (1010102)-q^{-1}(1010012)+q^{-2}(1001012)\\
 & +(q-1+q^{-1})((1010210)-q^{-1}(1001210)-q^{-1}(1010201)+q^{-2}(100\overline{1201}))\\
 & -(q-1+q^{-1})((1021010)-q^{-1}(1021001)+q^{-2}(1020101)-q^{-1}(0121010))\\
 & -(2101010)+q^{-1}(2101001)-q^{-2}(2100101)\\
 & -(q^{-2}-q^{-1})(0101210)+(q^{-2}-q^{-1})(0121001)-q^{-3}(0100121)+q^{-3}(0012101).
\end{align*}
We keep reducing it to
\begin{align*}
 & (1010102)-q^{-1}(1010012)\\
 & +(q-1+q^{-1})((10\overline{1021}0)-q^{-1}(1010201)-(\overline{1021}010)+q^{-1}(1021001)+q^{-1}(0121010))\\
 & -(2101010)+q^{-1}(2101001)-q^{-2}(2100101)\\
 & -(q^{-2}-q^{-1})(0101210)+(q^{-2}-q^{-1})(0121001)-q^{-3}(0100121)+q^{-3}(0012101)\\
 & +(q^{-1}-1)(1001210)-(q^{-3}+q^{-1})(1020101)-q^{-2}(1001201)+q^{-2}(1000121)\\
 =& -q^{-1}(1010012)\\
 & +(q-1+q^{-1})(-q^{-1}(1010201)+q^{-1}(\overline{1021}001))\\
 & +q^{-1}(2101001)-q^{-2}(2100101)\\
 & -(q^{-2}-q^{-1})(0101210)+(q^{-2}-q^{-1})(0121001)-q^{-3}(0100121)+q^{-3}(0012101)\\
 & +q^{-1}(1001210)-(q^{-3}+q^{-1})(1020101)-q^{-2}(1001201)+q^{-2}(1000121)\\
 & +(q^{-2}-q^{-1})(0121010)+(1012100)-(1210010)
\end{align*}
and then to, after rearrangement,
\begin{align*}
 & -q^{-1}(1010012)-q^{-3}(0100121)+q^{-2}(1000121)-q^{-2}(1001201)\\
 & +q^{-1}(1001210)-(q^{-2}-q^{-1})(0101210)-(q^{-2}+1)(1010201)\\
 & +q^{-3}(0012101)+(1012100)\\
 & -(q^{-3}+q^{-1})(1020101)+(q^{-2}-q^{-1})(0121010)+q^{-2}(0121001)\\
 & -q^{-1}(1021001)-(1210010)+q^{-1}(1210001)-q^{-2}(2100101).
\end{align*}
Note that the first term $-q^{-1}(1\overline{0100}12)$ equals $-q^{-1}(100\overline{1012})$
by (\ref{eq:even-Serre-0}) and $(11)=0$ and the last term $-q^{-2}(\overline{2101}001)$
so that we may reduce it further to
\begin{align*}
 & -q^{-3}(0100121)+(q^{-2}-q^{-1})(\overline{1000}121)\\
 & +(1001021)-(q^{-2}-q^{-1})(0101210)-(1010201)\\
 & +q^{-3}(0012101)+(1012100)\\
 & -q^{-3}(1020101)+(q^{-2}-q^{-1})(0121010)+q^{-3}(1021001)\\
 & -(1210010)+(q^{-1}-q^{-2})(121\overline{0001}).
\end{align*}
The middle line can be replaced with $q^{-3}(0010121)+(1210100)$
by (\ref{eq:odd-Serre}) (applied to $1210$ and $0121$ respectively)
and $(11)=0$. Now we use (\ref{eq:even-Serre-0}) to replace the overscored factors, arriving at
\begin{align*}
 & -(0100121)+(0010121)\\
 & +(1001021)-(q^{-2}-q^{-1})(0101210)-(1010201)\\
 & -q^{-3}(1020101)+(q^{-2}-q^{-1})(0121010)+q^{-3}(1021001)\\
 & -q^{-3}(1210010)+q^{-3}(1210100).
\end{align*}
Finally, since
\begin{align*}
(1\overline{0010}21)-(1010201) & =\frac{1}{q+1+q^{-1}}((1000121)-(1100021))\\
 & =\frac{1}{q+1+q^{-1}}((\overline{1000}121)-(0001121))\quad\text{as $(11)=0$ in either side}\\
 &=(0100121)-(0010121)
\end{align*}
and similarly
\[
-q^{-3}(1020101)+q^{-3}(1021001)=q^{-3}(1210010)-q^{-3}(1210100),
\]
it remains to see
\[
(q^{-2}-q^{-1})(-(0\overline{1012}10)+(0121010))=0.
\]

{\small


\begin{thebibliography}{O20}


\bibitem{BKK}
G. Benkart, S.-J. Kang, M. Kashiwara, {\em Crystal bases for the quantum superalgebra $U_q(\mathfrak{gl}(m,n))$}, J. Amer. Math. Soc. \textbf{13} (2000) 295--331.

\bibitem{CL}
S.-J. Cheng, N. Lam, {\em Irreducible characters of general linear superalgebra and super duality}, Commun. Math. Phys. \textbf{298} (2010) 645--672.

\bibitem{CLW}
S.-J. Cheng, N. Lam, W. Wang, {\em Super duality and irreducible characters of ortho-symplectic Lie superalgebras}, Invent. Math. \textbf{183} (2011) 189--224.

\bibitem{CW}
S.-J. Cheng, W. Wang, {\em Dualities and Representations of Lie Superalgebras}, Graduate Studies in Mathematics 144, Amer. Math. Soc., 2013.

\bibitem{CFLW}
S. Clark, Z. Fan, Y. Li, W. Wang, {\em Quantum supergroups III. Twistors}, Comm. Math. Phys. \textbf{332} (2014) 415--436.

\bibitem{KKKO15}
S.-J. Kang, M. Kashiwara, M. Kim, S.-j. Oh, {\em Simplicity of heads and socles of tensor products}, Compos. Math. \textbf{151} (2015) 377--396.

\bibitem{KKMMNN}
S.-J. Kang, M. Kashiwara, K. Misra, T. Miwa, T. Nakashima, A. Nakayashiki, {\em Perfect crystals of quantum affine Lie algebras}, Duke Math. J. \textbf{68} (1992) 499--607.

\bibitem{Kas95}
M. Kashiwara, {\em On crystal bases}, Representations of groups,
CMS Conf. Proc., 16, Amer. Math. Soc., Providence, RI, (1995) 155--197.

\bibitem{Kas02}
M. Kashiwara, {\em On level zero representations of quantum affine algebras}, Duke Math. J. \textbf{112} (2002) 117--175.

\bibitem{KKOP1}
M. Kashiwara, M. Kim, S.-j. Oh, E. Park, {\em Monoidal categorification and quantum affine algebras}, Compos. Math. \textbf{156} (2020) 1039--1077.

\bibitem{KashNaka}
M. Kashiwara, T. Nakashima, {\em Crystal graphs for representations of the $q$-analogue of classical Lie algebras}, J. Algebra \textbf{165} (1994) 295--345.

\bibitem{Kng}
R. C. King, {\em Weight multiplicities for the classical groups, Group theoretical methods in physics}, (Fourth Internat. Colloq., Nijmegen, 1975), pp. 490–499, Lecture Notes in Phys., Vol. 50, Springer-Verlag, Berlin-New York, 1976.

\bibitem{KO}
A. Kuniba, M. Okado, {\em Tetrahedron equation and quantum $R$ matrices for $q$-oscillator representations of $U_q(A^{(2)}_{2n})$, $U_q(C^{(1)}_n)$ and $U_q(D^{(2)}_{n+1})$}, Comm. Math. Phys. \textbf{334} (2015) 1219--1244.

\bibitem{KOS}
A. Kuniba, M. Okado, S. Sergeev, {\em Tetrahedron equation and generalized quantum groups}, J. Phys. A: Math. Theor. \textbf{48} (2015) 304001 (38p).

\bibitem{KL20}
J.-H. Kwon, S.-M. Lee, {\em Super duality for quantum affine algebras of type $A$}, Int. Math. Res. Not. \textbf{23} (2022) 18446--18525.

\bibitem{KL22}
J.-H. Kwon, S.-M. Lee, {\em Affinization of $q$-oscillator representations of $U_q(\mf{gl}_n)$}, Lett. Math. Phys. \textbf{113} (2023) 58.

\bibitem{KLO24}
J.-H. Kwon, S.-M. Lee, M. Okado,
{\em Oscillator representations of quantum affine orthosymplectic superalgebras}, Comm. Math. Phys. \textbf{405} (2024) 101.

\bibitem{KwO21}
J.-H. Kwon, M. Okado, {\em Higher level $q$-oscillator representations for $U_q(C^{(1)}_n)$, $U_q(C^{(2)}(n+1))$ and $U_q(B^{(1)}(0,n))$}, Comm. Math. Phys. \textbf{385} (2021) 1041--1082. 

\bibitem{KY}
J.-H. Kwon, J. Yu, {\em $R$-matrix for generalized quantum groups of type $A$}, J. Algebra, \textbf{566} (2021) 309--341.

\bibitem{Lu93}
G. Lusztig, {\em Introduction to quantum groups}, Progress in Math. \textbf{110}, Birkh\"{a}user, 1993.

\bibitem{Ma}
S. Machida, {\em Quantized superalgebras and generalized quantum groups}, M.S. Thesis (2017), Osaka City University, available at \url{https://www.omu.ac.jp/orp/ocami/assets/23_09.pdf}.

\bibitem{MY}
E. Mukhin, C. A. S. Young, {\em Affinization of category $\mathcal{O}$ for quantum groups}, Trans. Amer. Math. Soc. {\bf 366} (2014) 4815--4847.

\bibitem{OS}
S.-j. Oh, T. Scrimshaw, {\em Denominators of $R$-matrices, higher Dorey's rules and a generalization of T-systems for quantum affine algebras}, preprint (2025) arXiv:2510.10874.

\bibitem{O90}
M. Okado, {\em Quantum $R$ matrices related to the spin representations of $B_n$ and $D_n$}, Comm. Math. Phys. \textbf{134} (1990) 467--486.

\bibitem{Pro}
R.A. Proctor, {\em Odd symplectic groups}, Invent. Math. \textbf{92} (1988), 307--332.

\bibitem{Ya94}
H. Yamane, {\em Quantized enveloping algebras associated to simple Lie superalgebras and universal R-matrices},  Publ. Res. Inst. Math. Sci. \textbf{30} (1994) 15--84.

\bibitem{Ya99}
H. Yamane, {\em On defining relations of affine Lie superalgebras and affine quantized universal enveloping superalgebras}, Publ. Res. Inst. Math. Sci. \textbf{35} (1999) 321--390.


\end{thebibliography}
\end{document}